\documentclass[11pt, reqno]{amsart}

\usepackage{latexsym,amssymb,amsmath,eufrak}

\usepackage{amssymb,amsmath,eufrak}

\title{Regularization by free additive convolution, square and rectangular 
cases}

\author{Serban Belinschi,  
 Florent Benaych-Georges
and Alice Guionnet}
\address{Serban T. Belinschi, University of Saskatchewan and IMAR. Department of Mathematics and Statistics, University of Saskatchewan,  
106 Wiggins Road, Saskatoon, SK, S7N 5E6, 
Canada.} \email{belinsch@math.usask.ca}
\address{Florent Benaych-Georges, LPMA, UPMC Univ Paris 6, Case courier 188, 4, Place Jussieu, 75252 Paris Cedex 05, France.} \email{florent.benaych@gmail.com}
\address{Alice Guionnet, Ecole Normale Sup\'erieure de Lyon,
Unit\'e de Math\'ematiques pures et appliqu\'ees,
UMR 5669,
46 All\'ee d'Italie,
69364 Lyon Cedex 07, France.} \email{aguionne@umpa.ens-lyon.fr}

\date{\today}

\newcommand{\bck}{\backslash}

\newcommand{\ds}{\displaystyle}

\newcommand{\ninf}{\underset{n\to\infty}{\longrightarrow}}
\newcommand{\ssi}{if and only if }

\newcommand{\R}{\mathbb{R}}
\newcommand{\C}{\mathbb{C}}

\newcommand{\ud}{\mathrm{d}}

\newcommand{\trv}{Voiculescu transform }

\newcommand{\pro}{probability }

\newcommand{\f}{\frac}
\newcommand{\ff}{\frac{1}}
\newcommand{\lf}{\left}
\newcommand{\ri}{\right}

\newcommand{\st}{such that }
\newcommand{\la}{\lambda}

\newcommand{\vfi}{\varphi}
\newcommand{\ste}{\, ;\, }
\newcommand{\mc}{\mathcal }

\newcommand{\eps}{\varepsilon}
\newcommand{\arc}{{\scriptscriptstyle\boxplus_{\la}}}

\newcommand{\bxp}{{\scriptscriptstyle\boxplus}}

\def\e{{\epsilon}}
\def\ra{{\rightarrow}}
\numberwithin{equation}{section}

\newtheorem{Th}{Theorem}[section]
\newtheorem{propo}[Th]{Proposition}

\newtheorem{lem}[Th]{Lemma}

\newtheorem{rmq}[Th]{Remark}
\newtheorem{cor}[Th]{Corollary}

\newenvironment{pr}{\noindent {\bf Proof. }}{\ \ \ $\square$}

\def\ts

\d{
\lambda}

\long\def\symbolfootnote[#1]#2{\begingroup
\def\thefootnote{\fnsymbol{footnote}}\footnote[#1]{#2}\endgroup}

\begin{document}

\maketitle

\symbolfootnote[0]{{\it MSC 2000 subject classifications.}  primary 
46L54, 
60E10 
30A99, 
secondary 15A52.} 
 
\symbolfootnote[0]{{\it Key words.} Free convolution, rectangular free convolution, regularization, free infinitely divisible distributions, free probability theory, random matrices.}

\symbolfootnote[0]{This work was partially supported by Miller institute for Basic
 Research in Science, University of California Berkeley.}

\begin{abstract}
The free convolution $\bxp$ is the binary operation on the set of probability measures on the real line which allows to deduce, from the individual spectral distributions,  the spectral distribution of a sum of independent unitarily invariant square random matrices or of a sum of free operators in a non commutative probability space. In the same way,  the rectangular free convolution $\arc$ allows to deduce, from the individual singular 
 distributions,  the singular distribution  of a sum of independent unitarily invariant rectangular random matrices. In this paper, we consider 
the regularization properties of these free convolutions
 on the whole real line. More specifically, we try to find continuous semigroups $(\mu_t)$ of probability measures such that $\mu_0=\delta_0$ and such that for all $t>0$ and all probability measure $\nu$, $\mu_t \bxp\nu$ (or, in the rectangular context,  $\mu_t\arc \nu$) is absolutely continuous with respect to the Lebesgue measure, with a positive analytic density on the whole real line.  In the square case, for $\bxp$, we prove that in semigroups satisfying this property, no measure can have a finite second moment, and we give a sufficient condition on semigroups to satisfy this property, with examples. In the rectangular case, we prove that in most cases, for $\mu$ in a $\arc$-continuous semigroup, $\mu\arc\nu$ either has an atom at the origin or doesn't put any mass in a neighborhood of the origin, 
and
thus the expected property does not hold. However, we give sufficient conditions for analyticity of the density 
of $\mu\arc\nu$ except on a negligible set of points, as well as existence and continuity of a density everywhere. 
\end{abstract}

\tableofcontents

\section{Introduction}
It is a very natural question to study 
the spectrum of the sum of two matrices, being 
given the spectrum of each of them.
Such a question can of course have many different answers 
depending on the relations between the eigenspaces of the
two matrices. If they are the same,
but for instance the eigenvalues are independent
and say equidistributed inside each matrix, the spectrum 
will simply be given by the classical convolution. If, on the contrary the eigenspaces
are chosen as arbitrarily as possible 
with respect to each other,
which corresponds to 
conjugating one of the
matrices with an independent unitary 
matrix following the Haar measure, Voiculescu 
 \cite{Vo0}  proved
 that in the limit where the
size of the matrices goes to infinity while the spectral
measure of each 
of the two matrices converges
weakly, the outcome only depends on these limiting
measures and is given by their free convolution. 
More precisely, if we let $A_N,B_N$ be two sequences
of $N\times N$ Hermitian matrices with eigenvalues $(a^N_i)_{1\le i\le N}$
and $(b^N_i)_{1\le i\le N}$ respectively, such that the spectral measures
$L^N_A:=\frac{1}{N}\sum_{i=1}^N\delta_{a^N_i}$  and 
$L^N_B:=\frac{1}{N}\sum_{i=1}^N\delta_{b^N_i}$ converge 
to probability measures $\mu_A$ and $\mu_B$ as
$N$ goes to infinity, 
and if $U_N$ follows
the Haar measure on the unitary group and is independent of $A_N$ and $B_N$, 
then the spectral measure of 
$A_N+U_NB_NU_N^*$ converges towards the free convolution
 $\mu_A\bxp \mu_B$ of $\mu_A$ and $\mu_B$.

One of the authors, F. Benaych-Georges \cite{benaych.rectangular}, 
generalized  this convergence to the case
of rectangular matrices. In this case, $A_{N,M}$ and $B_{N,M}$ are
$N\times M$ matrices and we assume that
their singular values (for $N\leq M$, the singular values of an $N\times M$ matrix $C$ are the eigenvalues of $\sqrt{CC^*}$)
converge towards $\nu_A$ and $\nu_B$. We let, for $C=A$ or $B$,  $\mu_C$ 
be the symmetrization of $\nu_C$: $\mu_C(A)=\frac{1}{2}(\nu_C(A)+\nu_C(-A))$.
We consider $U_N$ and $V_M$ following
the Haar measure on the $N\times N$ and the $M\times M$
unitary matrices respectively. Then, F. Benaych-Georges
proved that, if $N/M$ converges to some $\lambda\in [0,1]$, 
then the symmetrization of the empirical measure of the 
singular values of  $A_{N,M}+ U_N B_{N,M}V_M$ converges
towards a probability measure $\mu_A\arc\mu_B$.

Free convolution naturally shows up 
in random matrix theory since important matrices
such as the Gaussian ensembles are invariant under 
conjuguation under the unitary group
and therefore can always be written as $U_NB_NU_N^*$ 
for some Haar distributed matrix $U_N$, independent of $B_N$,
or have asymptotically the same behaviour (for instance
matrices with independent equidistributed entries, see \cite{dykema}).

Convolution is a standard tool 
in classical analysis for regularizing functions or measures.
In this article, we study the regularizing properties of free (square 
and rectangular) additive convolution.
Because we wish to be able to regularize measures 
by  perturbations as small as 
desired, it is natural to regularize them
by processes $\mu_t, t\ge 0$ 
such that $\mu_t$ 
tends to $\delta_0$ as $t$ goes to zero.
To simplify, we shall consider more precisely
processes 
corresponding to infinitely divisible laws $\mu^{\boxplus t}$,
as constructed by Bercovici and Voiculescu  \cite{BercoviciVoiculescuIUMJ}
(see also 
 Nica and Speicher \cite{nicaspeicher}).

Such issues are naturally related 
with the possibility that the density vanishes,
since the density is then 
likely \cite{BA} 
to have some
infinite derivative at  the boundary of the support. 

Hence, we shall more precisely ask
the following question: can we find $\mu$
(likely a free infinitely divisible law) such that 

{\rm (H) } {\it For any probability measure $\nu$, $\mu
\bxp\nu$ (or in the rectangular case $\mu\arc\nu$) is absolutely continuous with respect
to the 
Lebesgue measure,  with a density which is analytic 
and which does not vanish on $\mathbb R$.}

Such questions already showed up in several papers.
In \cite{Vo1}, D. Voiculescu used 
the regularizing properties of 
free convolution by 
semi-circular laws to study free entropy for one variable.
In \cite{GZ,bianevoiculescu}, regularization by Cauchy laws is used
to smooth free diffusions in one case, 
and to study Wasserstein
metric and 
derive large deviations principles 
in the other. 
In fact, it was shown
in \cite{CDG2} that free convolution
is the natural concept to regularize
a free diffusion since 
the result will still be a free diffusion,
but with a different (and hopefully more regular)
drift. Free convolution by Cauchy laws
is well understood since it coincides 
with standard convolution with Cauchy laws.
In particular, Cauchy laws satisfy (H).
The drawback is that Cauchy laws do not possess any
moments, and thus do not allow a combinatorial approach
by moments. In section \ref{negan}, we show that it is in fact impossible
to find a $\boxplus$-infinitely divisible 
probability measure $\mu$ satisfying (H) and with finite variance.
In fact, we can then construct 
another probability measure $\nu$ such that the density
of $\mu\bxp\nu$ vanishes  at a point
where its derivative is infinite.
As a 
positive answer, we provide in section \ref{sectionanalysq}
sufficient conditions for a probability
measure $\mu$ to satisfy  (H). They require that
$\mu$ has  either none or  infinite first moment.

In the rectangular case,
we exhibit in section \ref{studyzero}
a sharp transition concerning the behaviour 
of the free rectangular convolution of two measures  at the origin.
If $\mu(\{0\})+\nu(\{0\})>1$, 
we prove that
$\mu\arc\nu(\{0\})$ is positive. This generalizes
a similar result of Bercovici and
Voiculescu in the square case (\cite{BVReg}). 
More surprisingly, when $\mu(\{0\})+\nu(\{0\})<1$,  
we show
the existence of a nonempty open neighbourhood of
the origin which does not intersect the support of
the density of $\mu\arc\nu
$,
for any infinitely divisible law $\mu$ and
any probability measure $\nu$. 

This phenomenon is related 
to the  repulsion at the
origin of 
the spectrum.  Such a repulsion 
was also shown to
hold at the finite matrices level
 in the square case by Haagerup \cite{haagerupunp}
 (by adding 
a form of Cauchy matrices) and by Sniady \cite{sniady} (by
adding Gaussian matrices). Our result is 
less strong since it is 
clear that the rectangular case carries naturally
a repulsion of the origin (as can be 
seen on the Pastur-Marchenko laws) and it
holds only asymptotically.
However, we find rather amazing 
that it holds for any infinitely divisible law
$\mu$ and any probability measure $\nu$.

This interesting phenomenon shows 
 that 
(H) cannot hold in the rectangular case. 
We thus show a weaker result
in Corollary \ref{final.corollary.rectangular}
and Proposition \ref{really-last}
by giving sufficient conditions
for analyticity of the density
of $\mu\arc\nu$  except on a discrete set
of points, as well as existence and
continuity of a density everywhere.

{\bf Acknowledgments.} The authors would like to thank  Professor Hari
Bercovici for useful discussions and suggestions, and  his encouragements
during the work on this paper.

\section{Prerequisites in complex analysis and free probability background}
\subsection{Complex analysis}
Let $\mathbb D:=\{z\in\mathbb C\colon|z|<1\}$.
We denote by $\sphericalangle\lim_{z\to w}f(z)$, $\sphericalangle f(w)$,
or
$$\lim_{\stackrel{z\longrightarrow w}{\sphericalangle}}f(z)$$
the limit of $f$ at $w\in\partial \mathbb  D$ along points inside any angle with
vertex
at $w$ and included
in $\mathbb D$, and name it the {\it nontangential limit of $f$ at $w$}.
{\it Unconditionnal limit} 
at $w$ 
corresponds to limit taken
over any path in the domain $D$ 
which ends at $w$.

We will denote $\C^+=\{z\in\C: \Im(z)>0\}$,
$\C^-=\{z\in\C:\Im(z)<0\}$ and $\R^+=[0,+\infty)$. The notion of nontangential limit extends
naturally to maps defined on a half-plane. Moreover, since the rational
transformation
$z\mapsto\frac{z-i}{z+i}$ of the extended complex plane $\mathbb
C\cup\{\infty\}$
carries the upper half-plane onto the unit disc, most properties of
analytic functions defined in
the unit disc transfer naturally to functions defined on $\mathbb C^+$.

We shall use in our paper several results describing the boundary
behaviour of
analytic functions defined in the unit disc or the upper half-plane.
Let us first cite a theorem due to Lindel\"of, (theorem 2.20(i) in
\cite{CL}).
\begin{Th}\label{lindel}
[Lindel\"of] Let $f$ be a meromorphic function on the upper half plane
such that there are at least
three points of $\C\cup\{\infty\}$ which are not attained by $f$ on the
upper half plane.
Consider $a\in \R\cup \{\infty\}$ such that there is a path
$\gamma : [0,1)\to \C^+$ with limit $a$ at $1$ such that
$$l:=\ds\lim_{t\to 1}f(\gamma(t))$$ exists in $\C\cup\{\infty\}$.
Then the nontangential limit of $f$  at $a$ exists and equals $l$.\end{Th}

Recall that if $f$
is a function from a subset $D$ of $\C\cup\{\infty\}$ into
$\C\cup\{\infty\}$, for all $z$ in
the boundary of $D$, the {\it cluster set} $C(f,z)$ of $f$ at $z$ is the
set of limits in
$\C\cup\{\infty\}$ of images, by $f$, of sequences of 
points in 
$D$ which tend to
$z$.  For 
any subset
$D'$ of $D$, $C_{D'}(f,z)$ denotes the cluster set of the restriction of
$f$ to $D'$.
For the particular case when $D=\mathbb D$ or $D=\mathbb C^+,$ we define
the {\it nontangential cluster set}
$C^\Delta(f, x_0)$ of $f$ at $x_0\in\partial D$ in the following way: let
$\Gamma(\alpha)$ be the angle with vertex at $x_0$
bisected by the perpendicular on $\partial D$ at $x_0$, with opening
$\alpha\in(0,\pi).$ Then

$$C^\Delta(f,
x_0)=\overline{\cup_{\alpha\in(0,\pi)}C_{\Gamma(\alpha)}(f,x_0)}.$$
Thus, the existence of nontangential limit of $f$ at $x_0$  means
   that $C^\Delta(f, x_0)$ contains only one point.

The following result (see e.g. Theorem 1.1 in \cite{CL})
   concerns the connectivity
of a cluster set.

\begin{lem}\label{lem1.2}
Let $D$ be a domain in $\mathbb C$ (i.e. an open connected set) and assume
that $D$ is simply connected
at the point $x\in\overline{D}$ 
(i.e. $x$ has a basis of neighbourhoods in
$\mathbb C$ whose intersections with
$D$ are simply
connected). If $f\colon D\to\mathbb C\cup\{\infty\}$ is continuous, then
$C(f,x)$ is connected.
\end{lem}

It is known from Fatou's Theorem (see \cite{CL}) that bounded
analytic functions on $\mathbb D$ 
have good boundary properties,
namely the nontangential limit of such a function
exists at almost all points (in the sense of linear measure) of the
boundary of $\mathbb D$. However, the set of points where the
nontangential limit
does not exist can be also quite rich in some situations, as the following
theorem shows (see e.g Theorem 4.8 in \cite{CL})

Recall  that a subset of a metric space $X$ is said to be {\it
residual}, or of {\it second
Baire category} if it isn't contained in the
union of any sequence of closed subsets of $X$ with empty
interior.

\begin{Th}\label{Th.4.8.CL} If the real or complex function $f(z)$ is
continuous in $|z|<1$, and if for $\theta\in [0,2\pi[$,
   $\{\mathcal G_\theta\}$ is
a  rotation by the angle $\theta$
   of a
continuum $\mathcal G_0$ such that $\mathcal G_0\cap\{|z|=1\}=\{1\}$, then
$
C_{\mathcal G_\theta}(f,e^{i\theta})=C
(f,e^{i\theta})$ on a residual set of points $e^{i
\theta}$ on $\{|z|=1\}$.
\end{Th}

    This theorem says for us that for a function that has no unconditional
limits at the boundary, the
nontangential limit must fail to exist sometimes.

We will use this result in connection with the following theorem of Seidel
(Theorem 5.4 in \cite{CL}).
\begin{Th}\label{Seidel}
Assume that $f\colon\mathbb D\to\mathbb D$ is analytic, and has
   nontangential limits with modulus one
at almost all points $\theta$ in some given interval
$(\theta_1,\theta_2)\subseteq\partial\mathbb D$.
Then for any $\theta_0\in(\theta_1,\theta_2)$ either $f$ extends
analytically through $\theta_0$, or
$C(f,\theta_0)=\overline{\mathbb D}$.
\end{Th}

Another useful auxiliary result is Theorem 5.2.1 from the same \cite{CL}:
\begin{Th}\label{521}
Let $f$ be  meromorphic in the domain $\mathcal D$ bounded by a smooth
   curve $\gamma$. Consider $z_0\in \gamma$ and suppose also that $f$
extends in a meromorphic function in an open set containing
$\gamma\backslash \{z_0\}$.
Then we have $\partial C_\mathcal D(f,z_0)
\subseteq C_\gamma(f,z_0)$, where $\partial A$ denotes the boundary (in
$\mathbb C\cup
\{\infty\}$) of $A\subseteq\mathbb C\cup\{\infty\}$.
\end{Th}

We shall also use the following theorem, which can be seen as a "nontangential limit" version of the analytic continuation principle (see \cite{CL}).
\begin{Th}[Riesz-Privalov]\label{26.03.07.1} Let $f$ 
 be an analytic function on $ \mathbb  D$. Assume that there exists a set $A$ of nonzero 
linear measure in $\partial \mathbb  D$ such that the nontangential limit of f exists at each point of $A$, and equals zero. 
Then $f (z) = 0$ for all $z\in \mathbb  D$.
\end{Th}

Consider now an analytic function
$f\colon\mathbb D\longrightarrow\overline{\mathbb D}.$
The {\it Denjoy-Wolff point} of $f$ is characterized by the fact that
it is the uniform limit on compact subsets of the iterates $f^{\circ n}=
\underbrace{f\circ f\circ\cdots\circ f}_{n\ {\rm times}}$
   of $f$. We state the
following theorem of Denjoy and Wolff  as it appears in Milnor's book
\cite{Milnor}, as Theorem 4.2.
Recall that an hyperbolic rotation  around some point
$z_0\in\mathbb D$ is a function of the form $z\mapsto
e^{i\theta}\frac{z-z_0}{1-\overline{z_0}z},$ $\theta\in\mathbb R$.

\begin{Th}\label{Denj}
Let $f\colon\mathbb D\to\overline{\mathbb D}$ be an analytic
function. Then either $f$ is a hyperbolic rotation around some point
$z_0\in\mathbb D$,
or the sequence of
functions $f^{\circ n}$ converges uniformly on compact subsets of $\mathbb
D$ to
a unique point $w\in\overline{\mathbb D}$, called the Denjoy-Wolff point
of $f$.
\end{Th}

Note that if 
the analytic function 
$f\colon\mathbb D\to\overline{\mathbb D}$ 
is not a hyperbolic rotation and has a fixed point $c\in\mathbb D$, 
then $c$
has to be the
Denjoy-Wolff point of $f$. In fact, the Denjoy-Wolff point
$w\in\overline{\mathbb D}$ can be equivalently characterized (see
Chapter 5 of \cite{shapiro})
   by being the unique point that
satisfies exactly one
of the following two conditions:
\begin{enumerate}
\item[(1)] $|w|<1$, $f(w)=w$ and $|f'(w)|<1$;
\item[(2)] $|w|=1,$ $\sphericalangle\lim_{z\to w}f(z)=w, $ and
$$\lim_{\stackrel{z\longrightarrow
w}{\sphericalangle}}\frac{f(z)-w}{z-w}\leq1.
$$
\end{enumerate}

Since the unit disc is conformally equivalent to the upper half-plane via
the
conformal automorphism of the extended complex plane
$z\mapsto\frac{z-i}{z+i},$
the above theorem and equivalent characterization of the Denjoy-Wolff
point
applies to self-maps of the upper half-plane $\mathbb C^+:=\{z\in\mathbb
C\colon
\Im z>0\}$, with the difference that when infinity is the Denjoy-Wolff of
$f$,
we have $$\sphericalangle\lim_{z\to\infty}f(z)/z\geq1.$$

\subsection{Free convolution, related transforms}
\subsubsection{Cauchy transform and Voiculescu transform}
We now recall some basic notions
about free convolution. Let us remind the reader
that for a probability $\mu$ on $\R$,
we denote $G_\mu$ its  Cauchy-Stieljes transform
$$G_\mu(z)=\int \frac{1}{z-x}\, d\mu(x),\quad z\in\C\backslash \R$$
and $F_\mu(z)=1/G_\mu(z)$. Note that $F_\mu\colon
\C^+\ra\C^+$.

The following theorem characterizes functions which appear as reciprocals
{ (in the sense of multiplication)} of Cauchy-Stieltjes
transforms of probabilities on the real line. For the proof and an
in-depth analysis of the subject,
we refer to \cite{Akhieser}, Chapter 3.

\begin{Th}\label{Nevanlinna}
Let $F\colon\mathbb C^+\to\mathbb C^+$ be an analytic function. Then there exist
$a\in\mathbb R$, $b\ge0$ and a positive finite measure $\rho$ on $\mathbb R$
so that
$$F(z)=a+bz+\int_{\mathbb R}\frac{1+tz}{t-z}d\rho(t),\quad z\in\mathbb
C^+.$$
 Moreover, 
$F$
is the reciprocal
of a Cauchy-Stieltjes transform of a probability measure on the real line
if and only if $b=1$.
The triple $(a,b,\rho)$ satisfies $a=\Re F(i)$, $b=\lim_{y\to+\infty}F(iy)/iy,$ and $b+\rho(\mathbb R)=\Im F(i)$.
\end{Th}

\begin{rmq}\label{ImF>Imz}
{\rm
An immediate consequence of Theorem \ref{Nevanlinna}
is that for
any probability measure $\sigma$ on $\mathbb R$, we have $\Im
F_\sigma(z)\geq
\Im z$
for all $z\in\mathbb C^+,$ with equality for any value of $z$ if and only
if
$\sigma$ is a point mass. In this case, the measure $\rho$ in the
statement of
   Theorem \ref{Nevanlinna} is zero.}
\end{rmq}

The function $F_\mu$ can be seen to be invertible
in a 
set of the
form
$$\Gamma_{\alpha,M}=\{z\in\C:
|z|\ge M, |\Im z|\ge \alpha |\Re z|\}
$$
   for some $M,\alpha>0$.

   The Voiculescu transform of $\mu$
(see paragraph 5 of \cite{BercoviciVoiculescuIUMJ})
is then given on 
$F_\mu(\Gamma_{\alpha,M})$ 
by $\phi_\mu(z)=F_\mu^{-1}(z) -z$.

The free convolution of two probability measures $\mu$ and $\nu$
on the real line is then characterized by
the fact that
\begin{equation}\label{freeco}
\phi_{\mu\boxplus\nu}(z)=\phi_\mu(z)+\phi_\nu(z)
\end{equation}
on the common component of their domain that contains $i[s,+\infty)$
for some large enough $s>0$ (for details, we refer again to \cite{BercoviciVoiculescuIUMJ}.)
Note here that Voiculescu's transform $\phi_\mu$ and
the so-called $R$-transform are related by $R_\mu(z)=\phi_\mu(1/z)$.

Another useful property of Cauchy-Stieltjes transforms of free
convolutions of probability measures
is subordination: for any $\mu,\nu$, there exist unique analytic functions
$\omega_1,\omega_2\colon\mathbb C^+\to\mathbb C^+$
so that $G_\nu(\omega_1(z))=G_\mu(\omega_2(z))=G_{\mu\boxplus\nu}(z)$ for
all $z\in\mathbb C^+$
and $\lim_{y\to+\infty}\omega_j(iy)/iy=1,$ $j=1,2.$ This has been proved
by Biane in \cite{Biane1}.

In the following, we shall also need the following lemmas.

\begin{lem}\label{Fatou} [Fatou's theorem] Let $f\colon\mathbb 
C^+\to\mathbb C$
be an analytic function. If $\mathbb C\setminus f(\mathbb C^+)$ contains a 
half-line,
then $f$ admits finite nontangential limits at Lebesgue-almost all points 
of the real line.
\end{lem}
This lemma follows from Theorem 2.1 of \cite{CL}, and conformal 
transformations.

\begin{lem}\label{19.09.06.8}Let $\nu$ be a \pro measure on the real line.

\begin{itemize}\item[(i)] For almost all (with respect to the Lebesgue
measure) real numbers $x$, the nontangential limit, at $x$, of
$-\ff{\pi}\Im G_\nu$ exists and is equal to the density, at $x$, of the
absolutely continuous part of $\nu$ (with respect to the Lebesgue
measure).
\item[(ii)] Let $I$ be an open interval of the real line. Then we have
equivalence between: \begin{itemize}
\item[(a)] The restriction of $G_\nu$ to $\C^+$ 
extends analyticaly to  an open set
of $\C$  containing $\C^+\cup I$.
\item[(b)] The restriction of $\nu$ to $I$ admits an analytic density.
\end{itemize}
Moreover, in this case, the density of the restriction of $\nu$ to $I$ is
$x\in I \mapsto -\ff{\pi}\Im G_\nu(x),$  
 where $G_\nu(x)$ is the value, at $x$, of the extension  mentioned in (a) of the restriction of $G_\nu$ to $\C^+$. \end{itemize}\end{lem}

\begin{pr} Part (i) is 
Theorem 3.16 from Chapter II of \cite{S&W}.

(ii) Suppose (a) to be true. Let us define, for $t>0$,
$\mc{C}_t=\f{tdx}{\pi(t^2+x^2)}$.
 It is the law of $tX$ when $X$ is a
   $\mc{C}_1$-distributed
random variable (hence a standard Cauchy variable),
and thus  converges weakly to $\delta_0$ when $t$ tends to zero.
Let us also define, for $t\geq 0$, $$\rho_t : x\in \R\mapsto\begin{cases}
   -\ff{\pi}\Im G_\nu(x+it)&\textrm{if $t>0$ or $x\in I$,}\\ 0& \textrm{in
the other case.}\end{cases}$$ Then for all $t>0$,
   $\rho_t$ is  the density of $\nu*\mc{C}_t$,
and converges weakly (i.e. against any continuous bounded function)
   to $\nu$ as $t$ tends to zero. So it suffices to prove that for all
$f$ compactly supported continuous function on $I$,
   $\int f(x)\rho_t(x)\ud x$ tends to $\int f(x)\rho_0(x)\ud x$ when
$t$ goes to zero, which is an easy application of dominated convergence
   theorem.

Suppose (b) to be true. 
It suffices to prove that for all $x\in I$, there is $\eps_x>0$ \st the
restriction of $G_\nu$ to $\C^+$ admits an analytic extension $g_x$ to
$\C^+\cup  B(x, \eps_x)$. (We denote by $B(x, \eps_x)$ the open ball of 
center $x$
and radius $\eps_x$.) Indeed, in this case, since for all $x,x'\in I$,
$g_x, g_{x'}$ coincide on $B(x, \eps_x)\cap B(x', \eps_{x'})\cap \C^+$,
one can define an analytic function on $$\C^+\cup(\ds\cup_{x\in I} B(x,
\eps_x))$$ which coincides with $G_\nu$ on $\C^+$ and with $g_x$ on every
$B(x, \eps_x)$. So let us fix $x\in I$.

Without loss of generality, we may assume that (1) $x=0$, (2) the analytic
density function $f$ of $\nu$ is defined analytic on $[-c,c]$ for some
$c>0$, (3) this function has radius of convergence $>c$ with power series
$f(t)=\sum_{n=0}^\infty a_nt^n,$ and
(4) the support of $\nu$ is contained in  $[-c,c]$ (since for all finite
measures $\nu_1, \nu_2$, $G_{\nu_1+\nu_2}=G_{\nu_1}+G_{\nu_2}$ and
$G_{\nu_1}$ is analytic outside of the support of $\nu_1$). It will be
enough to show that $G_\nu$ extends analytically through $(-c,c)$. Let
$\log$ be the analytic function on $\C\bck\{-it\ste t\in [0,+\infty)\}$
whose derivative is $\ff{z}$. For any $|z|\leq c$,  $z\in\C^+$, we have
the integral
$$G_{\nu}(z)=\int_{\mathbb
R}\frac{f(t)}{z-t}dt=\int_{[-c,c]}\sum_{n=0}^\infty
a_n\frac{t^n}{z-t}dt=\int_{[-c,c]}\sum_{n=0}^\infty
a_n\frac{t^n-z^n+z^n}{z-t}dt$$
$$=-\sum_{n=0}^\infty
a_n\left(\sum_{j=0}^{n-1}z^{n-j}\frac{c^{j+1}-(-c)^{j+1}}{j+1}+z^n[\log(z-c)-\log(z+c)]\right)$$
$$=\left(\sum_{n=0}^\infty
a_nz^n\right)[\log(z+c)-\log(z-c)]-\sum_{n=0}^\infty
a_n\left(\sum_{j=0}^{n-1}z^{n-j}\frac{c^{j+1}-(-c)^{j+1}}{j+1}\right).$$
(We can commute integral with sum because the function of $t\to$
$|f(t)/(z-t)|$ is obviously bounded uniformly on (a neighbourhood of,
even) $[-c,c]$, for any
$z$ in a compact subset of $\mathbb C^+\cap B(0,c)$ and the sum
is absolutely convergent because of the power series condition.)
We claim that in fact this formula defines an extension of $G_\nu$ on
$B(0,c).$ Indeed, the first sum is obviously convergent, while for the
second we have
$|z^{n-j}\frac{c^{j+1}-(-c)^{j+1}}{j+1}|\leq
c^{n-j}\cdot2|c^{j+1}|\frac{1}{j+1}=2\frac{c^{n+1}}{j+1},$ so, since
$$\sum_{n=0}^\infty a_nc^{n+1}\sum_{j=0}^n\frac{1}{j+1}<\infty,$$
(recall the radius of convergence), so must be the second term in the sum
above, for all $|z|\leq c.$ Thus, the function
$G_\nu$ admits an analytic extension on an open set containing  $\C^+\cup
(-c,c)$, and the formula is given above.
\end{pr}


\begin{lem}\label{s}Let $\nu$ be a \pro measure on the real line with 
support equal to $\R$ and which is   concentrated on a set of null 
Lebesgue measure. 
Then the cluster set of the restriction of its Cauchy transform to the 
upper half plane at any real point is the closure, in $\C\cup\{\infty\}$, 
of the lower half plane.
\end{lem}

\begin{pr} Let us apply the upper half-plane version of theorem 
\ref{Seidel} to the opposite of $\nu$'s Cauchy transform. Since by (ii) of 
lemma \ref{19.09.06.8}, for all real number $x$, $-G_\nu$ does not extend 
analytically to $x$, the theorem will imply what we want to prove. 
Indeed, part (i) of lemma \ref{19.09.06.8} implies that {the imaginary 
part} of $G_\nu$ has a null nontangential limit at Lebesgue-almost all 
points $x$ of the real line.
Since by Lemma \ref{Fatou} 
$G_\nu$ admits a finite nontangential limit at Lebesgue-almost all 
$x\in\mathbb R$, it follows that $-G_\nu$ satisfies the upper half-plane 
version of
Theorem \ref{Seidel}. This completes the proof.

\end{pr}

\begin{lem}\label{d}The set of symmetric \pro measures on the real line 
with support $\R$ and which are concentrated on a set of null Lebesgue 
measure is dense in the set of symmetric \pro measures for the topology of 
weak convergence.
\end{lem}

\begin{pr} Let $D$ be the set of symmetric \pro measures on the real line 
with support $\R$ and which are concentrated on a set of null Lebesgue 
measure. Since the set of symmetric \pro measures which are  finite convex 
combination of Dirac masses is dense in the set of symmetric \pro 
measures, it suffices to prove that for all real numbers $a$, 
$\ff{2}(\delta_{-a}+\delta_a)$  is in the closure of $D$.  This is clear, 
since if $\nu\in D$, then for all $\eps\in (0,1)$, 
$\f{1-\eps}{2}(\delta_{-a}+\delta_a)+\eps\nu\in D$.
\end{pr}

\subsubsection{Free infinite divisibility}

One can extend the notion of infinitely divisible law from
classical convolution
   to
free convolution: a probability measure $\mu$ is said to be
$\boxplus$-infinitely divisible if for any $n\in\mathbb N$ there exists
a probability $\mu_n$ so that
$\mu=\underbrace{\mu_n\bxp\mu_n\bxp\cdots\bxp\mu_n}_{n\ \rm
times}$.
It can be shown that any such measure embeds naturally in a semigroup of
measures $\{\mu^{\boxplus t}\colon
t\geq0\}$ so that $t\mapsto\mu^{\boxplus t}$ is continuous in the weak
topology,
$\mu^{\boxplus 1/n}=\mu_n$ for all $n\in\mathbb N$, $\mu_0=\delta_0$, and
$\mu^{\boxplus s+t}=
\mu^{\boxplus s}\bxp\mu^{\boxplus t}$ for all $s,t\geq0.$
It follows easily from (\ref{freeco}) that for any $t\geq0$, we have
$\phi_{\mu^{\boxplus t}}(z)=t\phi_\mu(z)$ for all points $z$ in the
common domain of the two functions.

   In \cite{BercoviciVoiculescuIUMJ}, Bercovici and Voiculescu have
completely
described infinitely divisible probability measures with respect to free
additive convolution in terms of their 
Voiculescu transforms.

\begin{Th}\label{infdiv+}
\
\begin{trivlist}
\item[\ {\rm (i)}] A probability measure $\mu$ on $\mathbb R$ is
$\boxplus$-infinitely divisible if and only if $\phi_\mu$ has an analytic
extension defined on $\mathbb C^+$ with values in $\mathbb C^-\cup\mathbb
R$.
\item[\ {\rm (ii)}] Let $\phi\colon\mathbb C^+\to
\mathbb C^-\cup\mathbb R$
be an analytic function. Then $\phi$ is a continuation of $\phi_\mu$ for
some $\boxplus$-infinitely divisible measure $\mu$ if and only if
$$\lim_{\stackrel{z\longrightarrow\infty}{\sphericalangle}}\frac{\phi(z)}{z}
=0.$$
\end{trivlist}
\end{Th}

   In the following lemma, we describe some properties of free convolutions
with $\boxplus$-infinitely divisible measures.

\begin{lem}\label{lem0}
Let $\mu$ be a 
probability which is not
purely atomic,
and $\nu$ an arbitrary probability
measure. We know that
there exist two subordination functions
$\omega_1$ and $\omega_2$ from $\C^+$ into $\C^+$
   so that
$F_\nu\circ\omega_1=F_\mu\circ\omega_2=F_{\mu\bxp\nu}.$ Moreover:
\begin{enumerate}
\item[{\rm(1)}] $\omega_1(z)+\omega_2(z)=F_{\mu\bxp\nu}(z)+z$, 
$z\in\mathbb
C^+$;
\end{enumerate}
Assume in addition that $\mu$ is $\boxplus$-infinitely divisible. Then
\begin{enumerate}
\item[{\rm(2)}] $\omega_1$ is the inverse function of 
$H(w)=w+\phi_\mu(F_\nu(w))$,
$w\in
\mathbb C^+,$ so in particular it has a continuous extension to $\mathbb
R\cup\{\infty\}
$ and $\omega_1(x)$ is finite for all $x\in\mathbb R$. Moreover,
$\omega_1(z)$ is the Denjoy-Wolff
point of the function $g_z\colon\mathbb C^+\to\mathbb C^+,$
$g_z(w)=w+z-H(w),$  for all $z\in\mathbb C^+\cup\mathbb R$;
\item[{\rm(3)}] $F_\mu$ has a continuous extension to $\mathbb 
R\cup\{\infty\}$,
analytic outside a discrete set in $\R$,
and $F_\mu(x)$ is finite for all $x\in\mathbb R$;
\item[{\rm(4)}] $\omega_2$ and $F_{\mu\bxp\nu}$ extend continuously to 
$\mathbb
R$. Moreover, if $F_{\mu\bxp\nu}(x)\in\mathbb R$, then
so are $\omega_1(x)$ and $\omega_2(x)$;
\item[{\rm(5)}] If all existing nontangential limits of $\phi_\mu$ at 
points of
$\mathbb R$ belong to $
\mathbb C^-$, 
then $F_\mu(x)\in\mathbb C^+$ for all $x\in\mathbb R$.
\end{enumerate}
\end{lem}

\noindent
\begin{rmq}\label{rmqdensity}
(4) already shows that
at $x\in\R$ such that $F_{\mu\bxp\nu}(x)\neq 0$,
 $\mu\bxp\nu$
is absolutely continuous with respect
to Lebesgue measure, with density $\frac{\Im(F_{\mu\bxp\nu})}{\pi
|F_{\mu\bxp\nu}|^2}(x)$.
Note also that at a point $x$
such that $F_{\mu\bxp\nu}(x)= 0$,
Lemma \ref{lem0} (4) implies that  $\omega_2(x)\in\mathbb R$, and by
   the definition of $\omega_2,$
together with Lindel\"{o}f's Theorem \ref{lindel},
$$0=F_{\mu\bxp\nu}(x)=
\lim_{\stackrel{z\longrightarrow x}{\sphericalangle}}F_\mu(\omega_2(z))=
\lim_{\stackrel{\ \
z\longrightarrow\sphericalangle\omega_2(0)}{\sphericalangle\ \ \ \
}}F_\mu(z).$$
Since $\mu$ is infinitely divisible,
   Proposition 5.1 (1) of
   \cite{BBercoviciMult} then
guarantees that $t\omega_2(0)$ is an atom of $\mu^{\boxplus t}$ for all
$t<1$.

\end{rmq}
\begin{pr}
Item (1) has been proved in \cite{BVReg},
and can be easily checked from \eqref{freeco} 
and analytic continuation.

   Item (2) is a direct consequence of the fact that $H$ satisfies the
conditions imposed on the function denoted also $H$ in Theorem 4.6
\cite{BBercoviciMult} (namely it is analytic in $\C^+$, it decreases the 
imaginary part,
and its derivative has strictly positive limit at infinity
  - in this case equal to one). Hence, by Theorem 4.6
\cite{BBercoviciMult} (2), $H$ is
invertible from $\omega_1(\C^+)$ onto
$\C^+$ and, since by \eqref{freeco}, $F_\nu=F_{\nu\bxp\mu}\circ H$,
its inverse is exactly $\omega_1$. Moreover, by Theorem 4.6, part (2)
of \cite{BBercoviciMult},
$\omega_1$ extends continuously to $\R$ with image in $\C^+\cup\R$, while
part (3) of the same theorem guarantees that $\omega_1(z)$ is the 
Denjoy-Wolff
point of $g_z$ for all $z
\in\mathbb C^+\cup\mathbb R$.

Item (3) follows from Proposition 2.8 (a) in \cite{serbanthesis}.

We prove now (4). Assume that there exists $x_0\in\mathbb R$ so that
$C(F_{\mu\bxp\nu},x_0)$ contains more than
one point (and hence, by Lemma
\ref{lem1.2}, is a continuum). Since by (2) $\omega_1(x_0)$ is finite,
Theorem 4.1
of \cite{BA} allows us to conclude that if
$C(F_{\mu\bxp\nu},x_0)\cap\mathbb C^+\neq\varnothing,$
then $F_{\mu\bxp\nu}$ extends analytically to $x_0$, providing a
contradiction. Thus, $C(F_{\mu\bxp\nu},x_0)$
contains either an interval, or the complement of an interval,
in $\mathbb R$. By (2) $\omega_1(x_0)$
exists and then, by (1),
  $\omega_1(x_0)\in\mathbb R$.
By definition, for any $c\in C(F_{\mu\bxp\nu},x_0)$ in such an 
interval there exists a sequence
$\{z_n^c\}_n\subset\mathbb C^+$
converging to $x_0$ so that
$\lim_{n\to\infty}F_{\mu\bxp\nu}(z_n^c)=c.$
Using (1) we obtain

$$\lim_{n\to\infty}\omega_2(z_n^c)=\lim_{n\to\infty}F_{\mu\bxp\nu}(z_n^c)-\omega_1(z_n^c)+z_n^c
=c-\omega_1(x_0)+x_0.$$
Now part (3) of the lemma allows us to conclude that
$$c  =  \lim_{n\to\infty}F_{\mu\bxp\nu}(z_n^c)
 = \lim_{n\to\infty}F_{\mu}(\omega_2(z_n^c))
 = F_\mu(c-\omega_1(x_0)+x_0),$$
for all $c$ in an interval. The analyticity of $F_\mu$ outside a discrete
subset of $\mathbb R$
implies (by analytic continuation) that $F_\mu(z)=z-(x_0-\omega_1(x_0))$
for all $z\in\mathbb C^+$, so
$\mu=\delta_{x_0-\omega_1(x_0)},$ contradicting our hypothesis. The
continuity of $\omega_2$ on the real
line follows now immediately from (1) and (2).

To prove item (5) assume that $F_\mu(x)\in\mathbb R$. We know that
equality $z=
F_\mu(z)+\phi_\mu(F_\mu(z))$ extends to $\mathbb R$. If $\phi_\mu$ has no
nontangential limit at the point $F_\mu(x),$ then, since the continuity of
$F_\mu$
guarantees  the existence of $\lim_{z\to x}\phi_\mu(F_\mu(z)),$
Lindel\"of's Theorem \ref{lindel}
provides a contradiction ($\phi_\mu$ has limit along the path
$F_\mu(x+i\mathbb R^+)$
but not nontangential, at $F_\mu(x)$). If $\phi_\mu$ has nontangential
limit at $F_\mu(x)
$, then, by hypothesis, it must be complex, so obviously $\Im
F_{\mu}(x)=-\Im\sphericalangle\phi_\mu(F_\mu(x))>0$. Contradiction again.
\end{pr}

\subsection{Rectangular free convolution, related transforms}
\subsubsection{Introduction to the rectangular free convolution and to the
related transforms}\label{def-rect}

We recall 
\cite{benaych.rectangular} the construction of the
{\it rectangular $R$-transform with ratio $\la$},
and of the {\it rectangular free convolution $\arc$ with ratio $\la$}, for
$\la\in [0,1]$: one can summarize the different steps of the construction
of the rectangular $R$-transform with ratio $\la$ in the following chain

$$
\begin{array}{l}\ds\underset{\substack{{\textrm{sym. prob.}}\\
{\textrm{measure}}}}{\mu}\,\,\longrightarrow \,\,\underset{
\substack{{\textrm{Cauchy}}\\
{\textrm{transf.}}}}{G_\mu}\,\,\longrightarrow\,\, H_\mu(z)=\la
G_\mu\lf(\ff{\sqrt{z}}\ri)^2+(1-\la)\sqrt{z}G_\mu\lf(\ff{\sqrt{z}}\ri)\,\,\longrightarrow\,\,\\
\ds  \underset{
\textrm{rect. $R$-transf. with ratio $\la$}}{C_\mu(z)=U\lf(
\f{z}{H_\mu^{-1}(z)}-1\ri),}\end{array}$$
where for all $z=\rho e^{i\theta}$, with $\rho\in (0,+\infty),
\theta \in [0, 
2\pi)$, $\sqrt{z}=\rho^{1/2}e^{i\theta/2}$ (note that $\sqrt{\cdot}$  is 
analytic on
$\C\backslash \R^+$), and $U$
is the inverse of $T-1$, where $$T(z)=(\la z+1)(z+1),\textrm{ i.e. } U(z)=
    \f{-\la-1+\lf[(\la+1)^2+4\la z\ri]^{1/2}}{2\la} \textrm{ (when $\la=0$,
$U(z)=z$)},$$
where $z\mapsto z^{1/2}$ is the analytic version of the square root on the
complement of the real non positive half line \st $1^{1/2}=1$ (i.e. for all $z=\rho e^{i\theta}$, with $\rho>0, \theta \in (-\pi, 
\pi)$, $z^{1/2}=\rho^{1/2}e^{i\theta/2}$).

Note that  the rectangular $R$-transform with ratio $1$ (resp. $0$), for
symmetric distributions,
is linked to  the \trv  by the relation
$C_\mu(z)=\sqrt{z}\vfi_\mu(1/\sqrt{z})$ (resp. $C_\mu(z)=z\vfi_\rho(z)$,
where $\rho$ is the push-forward of $\mu$ by the function $t\to t^2$).

The rectangular free convolution of two symmetric \pro measures $\mu,\nu$ 
on the
real line is the unique symmetric \pro measure 
whose 
rectangular
$R$-transform is the sum of the rectangular $R$-transforms 
of $\mu$ and $\nu$, 
 and it is denoted 
by $\mu\arc\nu$. Then, we have 
\begin{equation}\label{10.01.06.1}C_{\mu\arc\nu}=C_\mu+C_\nu.
\end{equation}
If $\la=0$, $\mu\arc\nu$ is
the symmetric law which push-forward by  $t\to t^2$
is the free convolution of the push-forwards by  $t\to t^2$ of $\mu$ and 
$\nu$, and if
$\la=1$, it is $\mu\bxp\nu$.

\begin{rmq}\label{rmqhow}[How to compute $\mu$ when we know $C_\mu$ ?] First, we have
$z/H_\mu^{-1}(z)=T(C_\mu(z))$, for $z\in \C\backslash \R^+$
small enough. From this, we can compute $H_\mu(z)$ for $z\in \C\backslash
\R^+$ small enough.
Then we can use the equation, for  $z\in \C\backslash \R^+$,
\begin{equation}\label{28.12.06.hangthedj}\ff{z}H_\mu(z)=
\la\lf(\ff{\sqrt{z}}G_\mu(\ff{\sqrt{z}})\ri)^2+(1-\la)\ff{\sqrt{z}}G_\mu(\ff{\sqrt{z}}).\end{equation}
Moreover, when $z\in \C\backslash \R^+$ is small enough, $1/\sqrt{z}$ is
large and in $\C^-$, so
$\ff{\sqrt{z}}G_\mu(\ff{\sqrt{z}})$ is close to $1$. $\ff{z}H_\mu(z)$ is
also close to $1$, and for $h,g$
complex numbers close to $1$, $$h=\la g^2+(1-\la)g\Leftrightarrow
g=V(h),\textrm{ with }V(z)
=\f{\la-1+((\la-1)^2+4\la z)^\ff{2}}{2\la}=U(z-1)+1.$$
So one has, for $z\in  \C\backslash \R^+$ small enough,
\begin{equation}\label{H->G}
\ff{\sqrt{z}}G_\mu(\ff{\sqrt{z}})=V(\f{H_\mu(z)}{z}).
\end{equation}
\end{rmq}

\subsubsection{A few remarks about $H_\mu$}\label{Remarks.about.H}

(a) One has, for $z\in  \C\backslash \R^+$ small enough,
$$\ff{\sqrt{z}}G_\mu(\ff{\sqrt{z}})=V(\f{H_\mu(z)}{z}).$$
Note that the function
$\ff{\sqrt{z}}G_\mu(\ff{\sqrt{z}})$
is analytic on $\C\bck \R^+$, hence
$$U\lf(\f{H_\mu(z)}{z}-1\ri)=V\lf(\f{H_\mu(z)}{z}\ri)-1=\frac{1}{\sqrt{z}}G_\mu\left(\frac{1}{\sqrt{z}}\right)-1$$
admits an analytic  extension to $\C\bck \R^+$. Note that one cannot
assert that this extension is given by
the same formula on the whole  $\C\bck \R^+$, but we know, by analytic
continuation,
that if one denotes this extension by $M_\mu$, one has, for all $z\in
\C\bck \R^+$,
$$\lf[2\la M_\mu(z)+1+\la\ri]^2=(\la-1)^2+4\la \f{H_\mu(z)}{z},\quad{\rm or, \ equivalently,}\quad H_\mu(z)=zT(M_\mu(z)).$$
Let us observe that $M_\mu(z)=\psi_\mu(\sqrt{z})=\psi_{\mu^2}(z)$, where
$\psi_\mu(z)=\int\frac{zt}{1-zt}d\mu(t)$
is the so-called moment generating function of $\mu$, and $\mu^2$ is the probability 
on $[0,+\infty)$ given by $\int f(t)\,d\mu^2(t)=\int f(t^2)\,d\mu(t)$ for all Borel bounded functions $f$. Hence, as noticed in 
Proposition 6.2 of \cite{BercoviciVoiculescuIUMJ}, $M_\mu$ maps the upper half-plane into itself and 
the left half-plane $i\mathbb C^+$ into the disc with  diameter
the interval $(\mu(\{0\})-1,0)$.

(b) $H_{\mu}$ maps $\C\setminus \R^+$ into itself and maps $i\mathbb C^+\cap
\mathbb C^+$
into $i\mathbb C^+$. Indeed, for $z\in \mathbb C^+$, 
$$H_\mu(z)=\underbrace{
G_\mu\lf(\ff{\sqrt{z}} \ri)
}_{
\in \C^+
}\lf(\underbrace{
\la G_\mu\lf(\ff{\sqrt{z}}\ri)+(1-\la)\sqrt{z}
}_{ \in \C^+}
\ri),$$ and the product of two elements of $\C^+$ cannot belong to $\R^+$.
(As a consequence, Lemma \ref{Fatou} guarantees that the restriction  of 
$H$ to
the upper half-plane has
nontangential limits at almost all points of the positive half-line.) The second statement
follows from part (a) above and the definition of $T$: let $z=x+iy$ be so that 
$x<0$ and $y>0$. Then
\begin{eqnarray*}
z(M_\mu(z)+1) & = & z\left(\int_{[0,+\infty)}\frac{tz}{1-tz}\,d\mu^2(t)+1\right) = \int_{[0,+\infty)}\frac{z}{1-tz}\,d\mu^2(t)\\ 
& = & \int_{[0,+\infty)}\frac{x(1-tx)-ty^2}{(1-tx)^2+(ty)^2}\,d
\mu^2(t)+iy\int_{[0,+\infty)}\frac{1}{(1-tx)^2+(ty)^2}\,d\mu^2(t).
\end{eqnarray*}
Since $x<0$, the real part of the above expression is negative, as, since
$y>0$, its imaginary part is  positive.
From (a) above,
$$
H_\mu(z)=zT(M_\mu(z)) = z(M_\mu(z)+1)(\lambda M_\mu(z)+1).
$$
Since $M_\mu(z)$ belongs to the upper half of the disc of diameter $(\mu(\{0\})-1,0)$,
$\lambda M_\mu(z)+1$ belongs to $\mathbb C^+\cap(-i\mathbb C^+)$, so that the 
product of $\lambda M_\mu(z)+1$ and $z(M_\mu(z)+1)$ must belong to $i\mathbb C^+
.$

(c)
Using the two previous remarks, we observe that if there exist
$r,c\in(0,+\infty)$
and a sequence $\{z_n\}_n\subset\mathbb C^+$ so that
$\lim_{n\to\infty}z_n=r$ and $\lim_{
n\to\infty}H_\mu(z_n)=c$, then the set $\{M_\mu(z_n)\colon n\in\mathbb
N\}$ has at most two
limit points, either both negative (if $r>c$) or one negative and one
non-negative (if $r\le c$).
Indeed, the formula above guarantees that

$$\lim_{n\to\infty}M_\mu(z_n)\in\left\{\frac{-(1+\lambda)\sqrt{r}\pm\sqrt{r(1-\lambda)^2+4\lambda
c}}{2\lambda\sqrt{r}}\right\}.$$

(d)
Let us define two properties, for functions defined on $\mathbb
C\setminus\mathbb R^+$. $$\mathrm{(P1)}\quad \quad\forall z\in\mathbb
C\setminus\mathbb R^+ , f(\overline{z})=-\overline{f(z)}.$$
$$\mathrm{(P2)}\quad \quad\forall z\in\mathbb C\setminus\mathbb R^+ ,
f(\overline{z})=\overline{f(z)}.$$It is easy to see that $\sqrt{\cdot}$
has the property (P1) and that for  $\mu$ symmetric \pro measure 
$G_\mu(1/\sqrt{z})$ has also
property (P1), hence $H_\mu$ has property (P2).
As a consequence, in view of (b), 
 $H_\mu(\mathbb R^-)\subset \mathbb R^-$ and $H_\mu(i\mathbb C^+)\subseteq
i\mathbb C^+$. 
Similarly, $H_\mu^{-1}$ also satisfies property (P2)
and hence also $C_\mu$ satisfies property (P2), so, in particular, $C_\mu((-a,0))\subseteq
\mathbb R$ for any $a>0$ so that $(-a,0)$ is included in the domain of $C_\mu$.

(e)
Let us denote by $x_0$ the largest number in 
$(-\infty,0)$ so that $H_\mu'(x_0)=0$ (we do not
exclude the case $x_0=-\infty$). Since $H_\mu(0)=0$ 
and $H_\mu'(0)=1$, $H_\mu^{-1}$ and 
 $C_\mu$ are  defined, and analytic, on the interval
$(H_\mu(x_0),0)$, and moreover, $C_\mu((H_\mu(x_0),0))\subseteq\mathbb R^-.$
Indeed, $H_\mu^{-1}$ is obviously defined and analytic on $(H_\mu(x_0),0)$, and
as $H_\mu(\mathbb R^-)\subseteq\mathbb R^-$ (by (b) and (d) above), we have
$\frac{x}{H_\mu^{-1}(x)}>0$ for all $x\in(H_\mu(x_0),0).$ Thus, 
$C_\mu(x)=U\left(\frac{x}{H_\mu^{-1}(x)}-1\right)$
is defined and analytic on $(H_\mu(x_0),0).$

To show that $C_\mu((H_\mu(x_0),0))\subseteq\mathbb R^-$ it is enough to prove that
$U\left(\frac{H_\mu(x)}{x}-1\right)<0$ for any $x\in(x_0,0)$. (As observed in Remark 
\ref{rmqhow}, the derivative of $H_\mu$ in zero is one, so that $H_\mu$ is increasing
on $(x_0,0)$.) This statement is due to
 the inequality $\frac{H_\mu(x)}{x}<1,$ $x\in(x_0,0)$. Now, $1/\sqrt{x}\in i\mathbb R^-$ and $\mu$ is symmetric,
so $G_\mu(1/\sqrt{x})\in i \mathbb R^+$, which implies that
$|G_\mu(1/\sqrt{x})|=\Im G_\mu(1/\sqrt{x})$ for all $x<0$. Remark \ref{ImF>Imz} implies that

$$\left|\frac{1}{\sqrt{x}}\right|=-\Im \frac{1}{\sqrt{x}}<-\Im F_\mu\left(\frac{1}{\sqrt{x}}\right)
=\frac{1}{\Im G_\mu\left(\frac{1}{\sqrt{x}}\right)},$$
so $\left|\frac{1}{\sqrt{x}}G_\mu\left(\frac{1}{\sqrt{x}}\right)\right|<1$, for any $\mu\neq\delta_0,$
$x<0$. Thus,
$0<\frac{1}{\sqrt{x}}G_\mu\left(\frac{1}{\sqrt{x}}\right)<1,\quad x<0.$
The definition of $H_\mu$ and the fact that $0<\lambda<1$ imply now the desired result.

(f) 
We have 
 \begin{eqnarray}
\lim_{x\to-\infty}H_\mu(x)&=&-(1-\lambda)\int t^{-2}\,d\mu(t)\quad(=-\infty \textrm{ if }\mu(\{0\})>0),\nonumber\\
\lim_{x\to-\infty}\frac{H_\mu(x)}{x}&=&
\lambda\mu(\{0\})^2+(1-\lambda)\mu(\{0\}).\label{latome}\\
\nonumber
\end{eqnarray}
This follows from the definition of $H_\mu$ together with the monotone convergence theorem:
recall that 

\begin{eqnarray*}
H_\mu(x) & = & \lambda G_\mu(1/\sqrt{x})^2+(1-\lambda)\sqrt{x}G_\mu(1/\sqrt{x})\\
& = & \lambda\left(\int\frac{1}{\frac{-i}{\sqrt{|x|}}-t}\,d\mu(t)\right)^2+(1-\lambda)
i\sqrt{|x|}\int\frac{1}{\frac{-i}{\sqrt{|x|}}-t}\,d\mu(t)\\
& = & -\lambda\left(\int\frac{\sqrt{|x|}}{1+t^2|x|}\,d\mu(t)\right)^2-(1-\lambda)
\int\frac{|x|}{1+t^2|x|}\,d\mu(t).
\end{eqnarray*}
(We have used the fact that $\mu$ is symmetric in the last equality.) Since $\lim_{x\to
-\infty}\frac{1}{1+t^2|x|}=\chi_{\{0\}}(t)$ and the convergence is dominated by $1$,
 \eqref{latome} follows.

Now observe that the functions $f_x(t)=\frac{|x|}{1+t^2|x|},$ $x<-1,t\in\mathbb R,$
satisfy $f_{x_1}(t)>f_{x_2}(t)$ iff $|x_1|>|x_2|$, $f_x(t)<t^{-2}$ and $\lim_{x\to-\infty}f_x
(t)=t^{-2}$, $t\in\mathbb R$, with the convention $1/0=+\infty$. So by the monotone 
convergence theorem,

$$\lim_{x\to-\infty}\int\frac{|x|}{1+t^2|x|}\,d\mu(t)=\int t^{-2}\,d\mu(t)\in(0,+\infty].$$
If $\int t^{-2}\,d\mu(t)<+\infty$, we deduce that
$$\lim_{x\to-\infty}\int\frac{\sqrt{|x|}}{1+t^2|x|}\,d\mu(t)=\lim_{x\to-\infty}
|x|^{-\frac{1}{2}}\lim_{x\to-\infty}\int\frac{|x|}{1+t^2|x|}\,d\mu(t)
=0$$
so that indeed
$\lim_{x\to-\infty}H_\mu(x)=-(1-\lambda)\int t^{-2}\,d\mu(t)$. This is also true when  $\int t^{-2}\,d\mu(t)=+\infty$.

\subsubsection{Free rectangular infinite divisibility}\label{10.01.06.12}
As for the free convolution, for any $\la \in [0,1]$, one can extend the 
notion of infinitely divisible law
   to rectangular
free convolution with ratio $\la$: a symmetric probability measure $\mu$ 
is said to be
$\arc$-infinitely divisible if for any $n\in\mathbb N$ there exists
a symmetric probability $\mu_n$ so that
$\mu=\underbrace{\mu_n\arc\mu_n\arc\cdots\arc\mu_n}_{n\ \rm
times}$.
It can be shown that any such measure embeds naturally in a semigroup of
measures $\{\mu^{\arc t}\colon
t\geq0\}$ so that $t\mapsto\mu^{\arc t}$ is continuous in the weak
topology,
$\mu^{\arc 1/n}=\mu_n$ for all $n\in\mathbb N$, $\mu_0=\delta_0$, and
$\mu^{\arc (s+t)}=
\mu^{\arc s}\arc\mu^{\arc t}$ for all $s,t\geq0.$
It follows easily from (\ref{10.01.06.1}) that for any $t\geq0$, we have
$C_{\mu^{\arc t}}(z)=tC_\mu(z)$ for all points $z$ in the
common domain of the two functions.

   In \cite{fbg05.inf.div}, the  infinitely divisible probability measures 
with respect to $\arc$ are completely described in terms of their 
rectangular $R-$transforms.

   \begin{Th}\label{10.01.06.2}
   A symmetric \pro measure on the real line is $\arc$-infinitely divisible 
\ssi there is a symmetric positive finite measure $G$ on the real line \st 
$C_\mu$ extends to $\C\bck\R^+$ and is given by the following formula:
   \begin{equation}\label{10.01.06.3}
  \forall z\in \C\backslash\R^+, \quad\quad C_\mu(z)=z\int_\R\f{1+t^2}{1-zt^2}\ud 
G(t).
   \end{equation}
   In this case, $G$ is unique and is called the {\it L\'evy measure} of 
$\mu$.
   \end{Th}

\begin{rmq}It would be useful, in order to know if the measures to which 
Lemma \ref{auxm
} can be applied are all $\arc$-infinitely 
divisible, to know if, as for the Voiculescu transform and $\bxp$-infinite 
divisibility, any symmetric \pro measure whose rectangular $R$-transform 
extends analytically to   $\C\bck\R^+$ is actually $\arc$-infinitely 
divisible. Unfortunately, the proof of the analogous result in the square 
case involves the fact that the Voiculescu transform of any \pro measure 
takes its values in the closure of the lower half-plane, and we still did 
not find the analogue of that 
fact 
 in the rectangular context.
\end{rmq}

In the following we shall describe some more or less obvious consequences of Theorem 
\ref{10.01.06.2}.
First we record for future reference
the geometry of the preimage of the complex plane via $T(z)=(\lambda z+1)
(z+1)$:
\begin{rmq}\label{T}
{\rm\begin{enumerate}
\item[(i)] $T^{-1}(\{0\})=\{-1/\lambda,-1\}$ and $T'(-(\lambda+1)/2\lambda)=0$;
\item[(ii)] $T((-\infty,-1/\lambda])=T([-1,+\infty))=\mathbb R^+,$ and $T$ is injective
on each of these two intervals;
\item[(iii)] $T((-1/\lambda,-(1+\lambda)/2\lambda])=T([-(1+\lambda)/2\lambda,-1))
=[-(1-\lambda)^2/4\lambda,0)$ and $T$ is injective
on each of these two intervals;
\item[(iv)] $T(-(1+\lambda)/2\lambda+i\mathbb R^+)=T(-(1+\lambda)/2\lambda-i\mathbb R^+)
=(-\infty,-(1-\lambda)^2/4\lambda]$ and $T$ is injective on each of these two sets;
\item[(v)] $\Re T(x+iy)=0$ iff $\lambda y^2=(\lambda x+1)(x+1).$ In particular, the
pre-image of the imaginary axis is an equilateral hyperbola whose branches go through
$-1$ and $-1/\lambda$ and the tangents at these points to the hyperbola are vertical.
More general, $\Re T(x+iy)/\Im T(x+iy)=c\ge0$ if and only if $x^2-y^2-2cxy
+\left(1+\frac1\lambda\right)(x-cy)+\frac1\lambda=0$. That is, the pre-image via $T$ of
any non-horizontal line going through the origin is an equilateral hyperbola going 
through $-1$ and $-1/\lambda$ and whose tangents at these two points are parallel
to the line $cy=x$.

\end{enumerate}

Let us denote $K_1=\{z\in\mathbb C^+\colon \Re z>-(1+\lambda)/2\lambda\}$,
$K_2=\{z\in\mathbb C^+\colon \Re z<-(1+\lambda)/2\lambda\}$,
$K_3=\{z\in\mathbb C^-\colon \Re z>-(1+\lambda)/2\lambda\}$,
$K_4=\{z\in\mathbb C^-\colon \Re z<-(1+\lambda)/2\lambda\}$.

Note that using the formula $T(z)=\la\lf[\lf(z+\f{\la+1}{2\la}\ri)^2-\f{(1-\la)^2}{4\la^2}\ri],$
one easily sees that $K_1=T^{-1}(\C^+)\cap \C^+, K_2=T^{-1}(\C^-)\cap \C^+,K_3=T^{-1}(\C^-)\cap \C^-,K_4=T^{-1}(\C^+)\cap \C^-.$

}
\end{rmq}

\begin{lem}\label{rect-id}
Let $\mu$ be a $\arc$-infinitely divisible probability measure. Then
\begin{enumerate}
\item $H_\mu$ is the right inverse of the analytic function $\mathbb C\setminus\mathbb R^+
\ni w\mapsto \frac{w}{T(C_\mu(w))}$, hence injective;
\item $\mu(\{0\})>0$ if and only if $\lim_{w\to-\infty}C_\mu(w)\in(-1,0]$. In that case,
\begin{equation}\label{23.02.07.1}\mu(\{0\})=\frac{-(1-\lambda)+\sqrt{(1-\lambda)^2+4\lambda 
T\left(\displaystyle\lim_{w\to-\infty}C_\mu(w)\right)}}{2\lambda}=1+\displaystyle\lim_{w\to-\infty}C_\mu(w),\end{equation}
or, equivalently, $\displaystyle\lim_{w\to-\infty}C_\mu(w)=\mu(\{0\})-1.$
\item $\pi>\arg H_\mu(z)\ge\arg z$ for all $z\in \mathbb C^+$, with equality if and
only if $\mu=\delta_0$. In particular, $C_\mu(H_\mu(\mathbb C^+))\subset K_1$;
\item $H_\mu$ is analytic around infinity whenever $\lim_{x\to-\infty}C_\mu(x)<-1.$
\end{enumerate}
\end{lem}

\begin{proof}
By the definition of $C_\mu$, Theorem \ref{10.01.06.2}, parts (d) and (e) 
of subsection \ref{Remarks.about.H}, and Remark \ref{rmqhow}
we obtain that $T(C_\mu(H_\mu(z)))=H_\mu(z)/z$ for all $z\in(-\infty,0)$, and,
by part (b) of subsection \ref{Remarks.about.H} and analytic continuation, 
for all $z\in\mathbb C\setminus\mathbb R^+$. This proves item 1.

We prove now item 2. Since the case $\mu=\delta_0$ is trivial, we exclude it from our
analysis. This allows us to assert that the L\'evy measure of $\mu$
 has a positive
mass, hence, by (\ref{10.01.06.3}),  that  $C_\mu((-\infty,0))\subset (-\infty, 0)$. Note that (\ref{10.01.06.3}) implies also that $\lim_{w\to-\infty}C_\mu(w)$ exists in $[-\infty, 0)$.
As $H_\mu((-\infty,0))\subseteq(-\infty,0),$ (by (b), (d) of subsection \ref{Remarks.about.H})
we have $H_\mu(x)/x>0$ for all $x\in(-\infty,0)$, 
and hence the relation $T(C_\mu(H_\mu(z)))=H_\mu(z)/z$
implies that $T(C_\mu(H_\mu((-\infty,0))))\subset \R^{+*}$ 
and therefore 
 $C_\mu(H_\mu((-\infty,0))\subset (-\infty,-\frac{1}{\lambda}]
\cup [-1,0]$. 
Since
$\lim_{x\uparrow0}C_\mu(H_\mu(x))=0,$ the continuity of $x\ra C_\mu(H_\mu(x))$
on $\R^-$ 
implies  that $C_\mu(H_\mu((-\infty,0)))\subseteq(-1,0).$

Using  part (f) of subsection \ref{Remarks.about.H} and part (1) of this lemma,
$$\lambda\mu(\{0\})^2+(1-\lambda)\mu(\{0\})=\lim_{x\to-\infty}H_\mu(x)/x=
\lim_{x\to-\infty}T(C_\mu(H_\mu(x))).$$
If $\mu(\{0\})>0$, then by (f) of subsection \ref{Remarks.about.H},  $\lim_{x\to-\infty}H_\mu(x)=-\infty$, hence  $$T(\lim_{w\to-\infty}C_\mu(w))=\lim_{w\to-\infty}
T(C_\mu(w))=\lim_{x\to-\infty}T(C_\mu(H_\mu(x)))=\lambda\mu(\{0\})^2+(1-\lambda)\mu(\{0\})\in(0,1),$$ so, since $C_\mu(H_\mu((-\infty,0)))\subseteq(-1,0),$ we have
$\lim_{w\to-\infty}C_\mu(w)=\mu(\{0\})-1\in(-1,0).$

Conversely, assume that $\lim_{w\to-\infty}C_\mu(w)\in(-1,0)$. We claim that
$\lim_{x\to-\infty}H_\mu(x)=-\infty.$ Indeed, assume to the contrary that
this limit is finite, and denote it by $c\in(-\infty,0).$ Then we have 
$0=\lim_{x\to-\infty}H_\mu(x)/x=\lim_{x\to-\infty}T(C_\mu(H_\mu(x)))=T(C_\mu(c)),$
so that $C_\mu(c)=-1$ or $-1/\la$. But $C_\mu$ is increasing on $(-\infty,0)$ (it follows easily from the differentiation of (\ref{10.01.06.3})), and
we have assumed that $\lim_{x\to-\infty}C_\mu(x)>-1.$ This is a contradiction.
The statement concerning the mass at  the origin follows since
by \eqref{latome}
\begin{eqnarray*}
0<\lambda\mu(\{0\})^2+(1-\lambda)\mu(\{0\})=\lim_{x\to-\infty}\frac{H_\mu(x)}{x} & = &
\lim_{x\to-\infty}T(C_\mu(H_\mu(z))) =  \lim_{w\to-\infty}T(C_\mu(w)).
\end{eqnarray*}
To conclude, one can easily deduce \eqref{23.02.07.1} from the previous equation.

To prove item 3, we claim first that $H_\mu(\mathbb C^+)\subseteq
\mathbb C^+.$  Assume this is not the case:  there exists a point $z_1$
in the upper half-plane so that $H_\mu(z_1)\in\mathbb C^-\cup \R$. Observe that by the relations \begin{equation}\label{23.02.07.2}\forall \alpha\in (0,\pi),  \quad \lim_{\substack{z\to 0\\ |\arg z-\pi|<\alpha}}H_\mu(z)/z=1,\end{equation}
there is  a point $z_0\in\mathbb C^+$ so that $H_\mu(z_0)
\in\mathbb C^+.$   Consider a segment $\gamma$
uniting $z_0$ and $z_1$. Then there must be a point $z_r\in\gamma$ so that $H_\mu(z_r)
\in\mathbb R$ and $H_\mu([z_0,z_r))\subset\mathbb C^+.$ Since $H_\mu(\mathbb C\setminus
\mathbb R^+)\subseteq\mathbb C\setminus
\mathbb R^+,$ we must have
$H_\mu(z_r)<0.$ But  then $C_\mu(H_\mu(z_r))<0$, so $T(C_\mu(H_\mu(z_r)))\in\mathbb R$.
Thus, we contradict the relation $H_\mu(z_r)=z_rT(C_\mu(H_\mu(z_r))).$ 
This assures us that $H_\mu(\mathbb C^+)\subseteq\mathbb C^+.$
  
To conclude the proof of item 3, we have to prove that $H_\mu(z)/z\in\mathbb C^+$ whenever $z\in\mathbb C^+.$ 
This is equivalent to $T(C_\mu(H_\mu(z)))\in\mathbb C^+$ whenever $z\in\mathbb C^+,$
i.e., since by \eqref{10.01.06.3}, $C_\mu(\C^+)\subseteq \C^+$, to $C_\mu(H_\mu(\mathbb C^+))\subseteq K_1$. Note first that by \eqref{23.02.07.2}, there are some points $z\in\C^+$ for which $H_\mu(z)/z\in\mathbb C^+$, i.e. $C_\mu(H_\mu(z))\in K_1$. Hence the inclusion $C_\mu(H_\mu(\mathbb C^+))\subseteq K_1$ can fail only if $C_\mu(H_\mu(\mathbb C^+))$ intersects
the line $-(1+\lambda)/2\lambda+i\mathbb R^+$, so that there exists a point
$w_0\in\mathbb C^+$ with the property that $T(C_\mu(H_\mu(w_0)))<0$. But then we
obtain that $H_\mu(w_0)=w_0T(C_\mu(H_\mu(w_0)))\in\mathbb C^-$, a contradiction.
This proves item 3.


We proceed now with proving item 4. 
As observed in the beginning of the proof of item 2, we must have $C_\mu(H_\mu((-\infty,0)))\subseteq(-1,0),$
so that $\lim_{x\to-\infty}C_\mu(H_\mu(x))\ge-1$, hence the hypothesis $$\lim_{x\to-\infty}C_\mu(x)<-1$$ implies $\lim_{x\to-\infty}H_\mu(x)=c
\in(-\infty,0)$
(note here that the limit exists by \eqref{latome}). Thus, $H_\mu(z)$ is analytic around infinity if and only if $$W_\mu : z\in \C\bck(0,+\infty) \mapsto \begin{cases}H_\mu(1/z)&\textrm{ if $z\neq 0$},\\ c&\textrm{ if $z= 0$},\end{cases}$$
extends analytically around  zero. The relation $H_\mu(z)=zT(C_\mu(H_\mu(z)))$ allows us to write
$$W_\mu(z)=\frac1z(\lambda C_\mu(W_\mu(z))+1)(C_\mu(
W_\mu(z))+1),$$
hence for $z\in \C\bck\R^+$ small enough so that $(1-\lambda)^2+4\lambda zW_\mu(z)\notin \R^-$,
$$C_\mu(W_\mu(z))-\frac{-\lambda-1+[(1-\lambda)^2+4\lambda zW_\mu(z)]^{1/2}}{2\lambda}=0.$$
This relation holds  for $z$ in $I\cap (-\infty,0)$, where $I\subset\mathbb R$ is a small enough interval centered at zero, as $\mu\neq\delta_0$. Thus, let us define $f\colon I\times (I+c)\to\mathbb R,$ by 
$$f(z,w)=C_\mu(w)-\frac{-\lambda-1+[(1-\lambda)^2+4\lambda zw]^{1/2}}{2\lambda}$$
(recall that $c=\lim_{x\to-\infty}H_\mu(x)<0$, hence if $I$ is small enough, $f$ is well defined). This function satisfies $f(z,W_\mu(z))=
0$ for all $z\in I\cap (-\infty,0)$. Hence $f(0,c)=\lim_{z\uparrow 0}f(z,W_\mu(z))=0.$
We observe that 
$$\partial_w f(z,w)=C_\mu'(w)-\frac{z}{[(1-\lambda)^2+4\lambda zw]^{1/2}},$$
so that $\partial_w f(0,w)=C_\mu'(w)>0$ for all $w\in (I+c)\cap(-\infty,0)$.
Thus, 
the conditions of the implicit function theorem are satisfied, so we conclude that
there exists a unique real map $g$, analytic on some subinterval $J$ of $I$, 
centered at zero, so that $g(0)=c$ and $f(x,g(x))=0$ for all $x\in J$.
The uniqueness guarantees that $g(x)=W_\mu(x)$ on their common
domain, and hence it is an analytic extension to the interval $J$ of $W_\mu(z)=H_\mu
(1/z)$. 
This concludes the proof.
\end{proof}

\section{The square case}
Below, we prove that
free convolution is regularizing,
namely that we can find a set of probability
measures (roughly $\boxplus$-infinitely divisible distribution whose
Voiculescu transform is sufficiently nice) such that  any probability
measure, once convoluted by one of
these measures, has a density with respect to Lebesgue measure which
is analytic and positive everywhere.
The fact that we require the density
to be analytic everywhere or positive everywhere
will be seen
in Proposition \ref{negative} to impose that
these regularizing measures have no finite second
moment. We shall give also some examples of such measures
after the proof of the Theorem.

\subsection{A result of analyticity}\label{sectionanalysq}
\begin{Th}\label{Thinfty}
Let $\mu$ be a $\boxplus$-infinitely divisible distribution. Assume that
the Voiculescu transform
$\phi_\mu$ satisfies the following conditions:
\begin{enumerate}
\item For any $x\in\mathbb R$, either $\sphericalangle\lim_{z\to
x}\phi_\mu(z)\in\mathbb C^-,$ or
$\phi_\mu$ has no nontangential limit at $x$;
\item Either {\rm (i)}
$\sphericalangle\lim_{z\to\infty}\phi_\mu(z)=\infty$, or {\rm(ii)}
$\sphericalangle\lim_{z\to\infty}\phi_\mu(z)\in\mathbb C^-,$ or {\rm
(iii)} $C^\Delta(\phi_\mu,\infty)$ contains more than one point.
\end{enumerate}
   Then $\mu\bxp\nu$ has a
positive, everywhere analytic density for all  probability measure $\nu$.
\end{Th}

\begin{pr}
We have, on a
neighbourhood of infinity,
$$F_{\nu}^{-1}(z)+\phi_\mu(z)=z+\phi_\nu(z)+\phi_\mu(z)=z+\phi_{\mu\bxp\nu}(z)=
F_{\mu\bxp\nu}^{-1}(z).$$
We replace $z$ by $F_{\mu\bxp\nu}(z)$, for $z$ in an appropriate
truncated cone $\Gamma_{\alpha,M}$ - see 2.2.1 above 
(possible since $F_{\mu\bxp\nu}(z)$ is defined on the
whole upper half plane
and is equivalent to $z$ as $z$ goes to infinity in a nontangential way),
and get
$$z-\phi_\mu(F_{\mu\bxp\nu}(z))=F_{\nu}^{-1}(F_{\mu\bxp\nu}(z)).$$
Applying $F_\nu$ in both sides and using analytic continuation
(recall that $\phi_\mu$ extends to $\C^+$
by Theorem \ref{19.09.06.8}), we obtain
\begin{equation}\label{cup.of.tea}
F_{\mu\bxp\nu}(z)=F_\nu(z-\phi_\mu(F_{\mu\bxp\nu}(z))),\quad
z\in\mathbb C^+.
\end{equation}
We will show that $F_{\mu\bxp\nu}$ extends analytically to $\mathbb
R$, and
$F_{\mu\bxp\nu}(x)\in\mathbb C^+$ for all $x\in\mathbb R$.
   This will imply the
theorem
according to Lemma \ref{19.09.06.8}. Let us fix a real number $x$.
Observe first that the existence of a continuous extension
with values in $\C^+\cup\R\cup\{\infty\}$
of $F_{\mu\bxp\nu}$ at $x$ is guaranteed by Lemma \ref{lem0}, (4).

We are first going to prove that we do not have
$\lim_{z\to x}F_{\mu\bxp\nu}(z)=\infty$
(i.e $\lim_{z\to x}|F_{\mu\bxp\nu}(z)|=\infty$).
Suppose that  this happens. Then we have

\begin{eqnarray}
x-\omega_1(x) & = & \lim_{z\to x}z-\omega_1(z) =  \lim_{z\to x}\omega_2(z)-F_\mu(\omega_2(z))\nonumber\\
& = & \lim_{z\to x}\phi_\mu(F_\mu(\omega_2(z))) =  \lim_{z\to x}\phi_\mu(F_{\mu\bxp\nu}(z))\nonumber\\
& = &
\lim_{\stackrel{w\longrightarrow\infty}{\sphericalangle}}\phi_\mu(w)\label{eqi}\\
\nonumber
\end{eqnarray}
We have used part (2) of Lemma \ref{lem0} in the first equality,
(1) of Lemma \ref{lem0} in the second equality,
definition of $\phi_\mu$ and Theorem \ref{infdiv+} in the third,
and Lindel\"{o}f's
Theorem \ref{lindel}
   in the 
last
 equality.
We next show that  any of the three hypotheses
of Theorem \ref{Thinfty}.(2)
are in contradiction with \eqref{eqi}.

Indeed, if (i) holds, then \eqref{eqi} implies that $x-\omega_1(x)=\infty$
which
   contradicts
part (2) of Lemma \ref{lem0}.
(iii) clearly can not hold since \eqref{eqi} implies that
$C^\Delta(\phi_\mu,\infty)$
contains only one point.
Finally,
assume  that (ii) of item 2 of our theorem holds.
Then, \eqref{eqi}
implies that $\omega_1(x)\in\mathbb C^+.$ Thus,  we get the
contradiction $$
\infty=\lim_{z\to x}F_{\mu\bxp\nu}(z)=\lim_{z\to
x}F_{\mu}(\omega_1(z))
=F_{\mu}(\omega_1(x))\in\mathbb C^+$$


We next prove that  $c=\lim_{z\to x}F_{\mu\bxp\nu}(z)$
cannot be real. So, we assume $c\in\R$ and obtain a contradiction
based on the fact that we then have 
\begin{equation}\label{answer.alice.19.09.06.1}
c=\lim_{z\to x}
F_{\mu\bxp\nu}(z)=\lim_{z\to x}F_\nu(z-\phi_\mu(
F_{\mu\bxp\nu}(z)))\in\mathbb R.
\end{equation}
We observe first that $\lim_{z\to x}\phi_\mu(
F_{\mu\bxp\nu}(z))$ exists.
Indeed, we see as  in \eqref{eqi},
$$\lim_{z\to x}\phi_\mu(
F_{\mu\bxp\nu}(z))=\lim_{z\to x}\omega_2(z)-F_\mu(\omega_2(z))=
\lim_{z\to x}z-\omega_1(z)=x-\omega_1(x)$$
If $\phi_\mu$ has nontangential limit at $c$, then by our assumption from
item 1, it must belong to the lower half-plane. Lindel\"{o}f's Theorem
\ref{lindel}
guarantees  that $\lim_{z\to x}\phi_\mu(
F_{\mu\bxp\nu}(z))$ equals the nontangential limit of $\phi_\mu$ at
$c$,
so
$$\lim_{z\to x}F_\nu(z-\phi_\mu(
F_{\mu\bxp\nu}(z)))=F_\nu(x-\lim_{\stackrel{w\longrightarrow
c}{\sphericalangle}}\phi_\mu(w))\in\mathbb
C^+,$$
contradicting equation (\ref{answer.alice.19.09.06.1}).

If $\phi_\mu$ has no nontangential limit at $c$, then it is obvious from
the existence of $\lim_{z\to x}\phi_\mu(
F_{\mu\bxp\nu}(z))$, of $c=\lim_{z\to x}F_{\mu\bxp\nu}(z)$,
finiteness of $c$, and from Lindel\"of's Theorem \ref{lindel}
that $c$ must belong to the upper half-plane.

Hence,
we have proved that, for
any $x\in \R$,
$c= F_{\mu\bxp \nu}(x)=\lim_{z\ra x} F_{\mu\bxp \nu}(z)\in
\C^+.$
We finally prove that $F_{\mu\bxp\nu}$
extends analytically in the neighbourhood of $x\in\R$
by using the implicit function theorem.
Note that the hypothesis that $\sphericalangle\lim_{z\to t}\phi_\mu(z)\in
\mathbb C^-$ for all $t\in\mathbb R$ for which this limit exists
implies that $\mu$ is not a Dirac measure, hence that
$\phi_\mu(\C^+)\subset
\mathbb C^-$. Let us introduce the function $f(v,w)=F_\nu(v-\phi_\mu(w)),$
defined on
$$\{(v,w)\in\C\times \C^+\ste \Im v >\Im \phi_\mu(w)\},$$which contains
$(\C^+\cup \R)\times \C^+$.  One has
$$f(x,c)=\ds F_\nu(x-\phi_\mu(c))=\lim_{z\to
x}F_\nu(z-\phi_\mu(F_{\mu\bxp\nu}(z)))=\lim_{z\to
x}F_{\mu\bxp\nu}(z)=c.$$
In other words, $c$ is the Denjoy-Wolff point
of the function $f(x,\cdot)$.
Since $f(x,\cdot)$ 
is not a conformal automorphism of $\mathbb C^+,$ we have
$$\lf |\f{\partial }{\partial w}f(x,c)\ri |<1
$$ (and in particular $\neq1$). So with $g(v,w)=w-f(v,w)$, we have
$$\f{\partial }{\partial w}g(x,c)\neq 0,$$
hence,  by the implicit function
theorem,  there exists an analytic function $L$, defined in a connected
neighborhood $V$ of $x$ and a neighborhood $W$ of $c$ such that  for all
$(v,w)\in V\times W$,
$$g(v,w)=0\Leftrightarrow w=L(v).$$
By (\ref{cup.of.tea}) the function $L$ coincides with the function
$F_{\mu\bxp \nu}$ on $V\cap \C^+$, so the function $F_{\mu\bxp
\nu}$
admits an analytic extension to $V$, with value $c\in \C^+$ at  $x$. Lemma
\ref{19.09.06.8} allows us to conclude.
   \end{pr}

\noindent
{\bf Examples.} In this series of examples,
we provide explicit examples of measures satisfying the hypotheses
of Theorem \ref{Thinfty}.

\begin{enumerate}
\item  We give here an example of a Voiculescu transform that satisfies
condition (i). Let
$$\phi_\mu(z)=\frac{1}{z+i}-\sqrt{z},\quad z\in\overline{\mathbb C^+},$$
where $\sqrt{\cdot}$ is the natural continuous extension of the square
root defined
on $\mathbb C\setminus [0,+\infty)$ so that $\sqrt{-1}=i$ to $\mathbb
R\cup\{\infty\}$.
Theorem \ref{infdiv+} guarantees that $\phi_\mu$ is the Voiculescu
transform of a $\boxplus$-infinitely
divisible probability. Clearly $\Im\phi_\mu(z)<0$ for all $z\in\mathbb
C^+\cup\mathbb R$, $\phi_\mu(\infty)
=\infty,$ and $\phi_\mu$ is obviously continuous on $\overline{\mathbb
C^+},$ so $\phi_\mu$
satisfies the condition (i) in the previous theorem. 
Moreover, we have $$\inf_{x\in\mathbb R}|\Im\phi_\mu(x)|=\lim_{x\to+\infty}\frac{1}{x^2+1}=0.$$
This also shows that $\mu$ is not a convolution with a Cauchy law.
It is an easy exercise to observe, based on (\ref{freeco}), that
$\mu=\lambda_1\bxp\lambda_2$, where $\lambda_j$ are both infinitely
divisible, $\lambda_2$ is a $\boxplus$-stable distribution (see
\cite{BercoviciVoiculescuIUMJ}) whose density is given by $x\mapsto
\frac{\sqrt{4x-1}}{2\pi x},$ $x\in[1/4,\infty)$, and
$F_{\lambda_1}(z)=(z-i+\sqrt{(z+i)^2-4})/2$, $z\in\mathbb C^+.$

\item  We observe that, if $\phi_\mu$ extends continuously to $\mathbb R\cup\{\infty\}$,
 condition (ii) in the above theorem
can be 
 reduced to requiring that $\mu$ is the free convolution
of some probability measure by a Cauchy law. Indeed, since $\phi_\mu$
is continuous on $\mathbb R\cup\{\infty\}$, if $\Im \phi_\mu(x)<0$ for all
$ x\in \mathbb R$ and $\phi_\mu(\infty)\in-\mathbb C^+$, then $\Im
\phi_\mu(x)$ must actually be
   bounded away from zero (by continuity). By Theorem \ref{infdiv+}, there exists $c>0$ so that
   $\phi_\mu(z)+ci=\phi_\mu(z)-(-ci)$ is still a Voiculescu transform of an
infinitely divisible distribution, say $\eta$. Then
$\phi_\mu(z)=\phi_\eta(z)+(-ci)$, and $\mathfrak C(z):=-ci$ is
the Voiculescu transform of a Cauchy distribution.

\item
The example of probability measure that satisfies condition (iii) will be
constructed in terms of the Voiculescu transform, as an explicit limit of
compactly supported probabilities, each whose density is an algebraic 
function.
Specifically, we will construct two sequences $\{a_n\}_n$ and $\{t_n\}_n$
of real numbers and functions $f_n(z)=a_n\frac{1+t_nz}{t_n-z},$ 
$z\in\mathbb C
\setminus \{t_n\}$, so that $g_n=\sum_{j=1}^nf_j$ converges on the upper
half-plane to the nonconstant analytic function $g$ and
$C_{i\mathbb R_+}(g,\infty)=iC_{i\mathbb R_+}(\Im g,\infty)\supset 
i[7,+\infty]
.$

Let us recall that the $R$-transform of the free Poisson law with 
parameter $k$ is
$R_p(z)=\frac{k}{1-z},$ so that its dilation by $t$ has an $R$-transform 
given by
$\frac1tR_p(\frac{z}{t})=\frac{k}{t-z}.$ Since $\phi_\mu(z)=R_\mu(1/z),$
we can write
$$a\frac{1+\frac{t}{z}}{\frac1z-t}=\frac{a}{t}\cdot\frac{z+t}{\frac1t-z}=\frac{a}{t}\left[
\frac{\frac1t+t}{\frac1t-z}-1\right]=\frac{\frac{a}{t^2}+a}{\frac1t-z}-\frac{a}{t},$$
so we conclude that $f_n$ is just minus the Voiculescu transform of the 
translation
with $a_n/t_n$ of the dilation with $1/t_n$ of the free Poisson law of 
parameter
$\frac{a_n}{t_n^2}+a_n$. 

First let us enumerate some properties of the functions $f_n$.
\begin{trivlist}
\item{(j)} $y\mapsto\Im f_n(iy)$ is a smooth function from $[0,+\infty)$ 
into itself
and $\Im f_n(iy)=a_n\frac{y(1+t_n^2)}{t_n^2+y^2}$;
\item{(jj)} $\max_{y\in[0,+\infty)}\Im f_n(iy)=\Im 
f_n(it_n)=a_n(1+t_n^2)/2t_n
$. Moreover, the function  $y\mapsto\Im f_n(iy)$ increases from zero to
$a_n(1+t_n^2)/2t_n$ on the interval $[0,t_n]$, after which it decreases 
back to
zero;
\item{(jjj)} There are exactly two points $y_n^+$ and $y_n^-$, right and 
left,
respectively, from $t_n$, so that
$\Im f_n(iy_n^+)=\Im f_n(iy_n^-)=1.$ We have
$$y_n^-=\frac{a_n(1+t_n^2)-\sqrt{a_n^2(1+t_n^2)^2-4t_n^2}}{2},$$
so that $\lim_{|t_n|\to\infty}y_n^-=1/a_n$. Moreover, for any $a_n>0$, 
$y_n^-<2/a_n,$
and if $0<a_n<1,$ then we also have $1/a_n<y_n^-.$
\end{trivlist}

Let us observe also that if we replace $f_n$ by the sum
$$f_n(z)=\frac{a_n}{2}\frac{1+zt_n}{t_n-z}+\frac{a_n}{2}\frac{1-zt_n}{-t_n-z},$$
then we do not change the imaginary part of $f_n(iy),$ while we insure 
that
$\Re f_n(iy)=0$ for all $y\geq0$, so from now on we will replace $f_n$ 
with
this new function. (This will correspond to the free additive convolution 
of two free
Poisson laws as above.)

Let $a_1=1,$ $t_1=2$. Choose $0<a_2<a_1/2$ so that $\Im 
f_1\left(i\frac{1}{2a_2}
\right)<1/10$
and $\frac{1}{2a_2}>t_1$.
Item (jj) guarantees that we can make such a choice. Choose
$t_2>t_1$ so that 
$\Im f_2(it_2)=a_2(1+t_2^2)/2t_2>2$. 
This condition can be fulfilled because of item (jj)
.
Observe that the monotonicity of $\Im f_n(i
\cdot)$
on $[t_n+\infty)$ implies that $\Im f_1(iy_2^-)<\Im 
f_1\left(i\frac{1}{2a_2}
\right)<1/10$.

Assume now that we have constructed $a_j$, $t_j$, $1\le j\le n-1$ so that 
$0<
a_j<a_{j-1}/2$, $f_{j-1}\left(i\frac{1}{2a_j}\right)<1/10^{j-1},$
$\frac{1}{2a_j}>t_{j-1}$,
$t_j>t_{j-1}$, $\Im f_j(it_j)=a_j(1+t_j^2)/2t_j>j$, and 
$y_j^->\frac{1}{2a_j},$
for all $j$ between $1$ and $n-1$. We choose $0<a_n<a_{n-1}/2$
small enough so that
$\Im f_{n-1}\left(i\frac{1}{2a_n}\right)<1/10^{n-1},$ 
$\frac{1}{2a_n}>t_{n-1},$
  (using item (jj) above),
and $t_n>t_{n-1}$ large enough so that 
$\Im f_n(it_n)=a_n(1+t_n^2)/2t_n>n$.
As before, construction is permitted by
using item (jj).

Observe now that the sequence $\{a_n\}_n$ constructed this way is 
positive,
decreasing, and satisfies $1<\sum_{n=1}^\infty a_n<2.$ Moreover,
$\{y_n^-\}_n$ is, by item (jjj),
in its own turn an increasing sequence, and $y_n^->t_j$ for
all $j<n$
so that $\Im f_j(iy_n^-)\le\Im f_j(iy_{j+1}^-)<\Im f_j
\left(i\frac{1}{2a_{j+1}}
\right)<1/10^{j}$, by monotonicity of $\Im f_j$ on $[t_n,+\infty)$. Thus,
$$\Im(f_1+f_2+\cdots+f_{n-1}+f_n)(iy_n^-)<\frac{1}{10}+\frac{1}{10^2}+\cdots
\frac{1}{10^{n-1}}+1<\frac{10}{9}.$$
On the other hand, for any $m>n$, we have
\begin{eqnarray*}
\Im f_m(iy_n^-) = a_m\frac{y_n^-(1+t_m^2)}{t_m^2+(y_n^-)^2} & = &
\underbrace{a_n\frac{y_n^-(1+t_n^2)}{t_n^2+(y_n^-)^2}}_{\Im 
f_n(iy_n^-)=1}\cdot\frac{a_m}{a_n}\cdot
\frac{(1+t_m^2)(t_n^2+(y_n^-)^2)}{(1+t_n^2)(t_m^2+(y_n^-)^2)}\\
& \leq & 
\frac{a_m}{a_n}\cdot\frac{(1+t_m^2)(t_n^2+\left(\frac{2}{a_n}\right)^2)}{(1+t_n^2)t_m^2}\\
& < & 2\frac{a_m}{a_n}\left(1+\frac{4}{a_n^2t_n^2}\right)\\
& < & \frac{1}{2^{m-n-1}}+\frac{1}{2^{m-n-3}n^2},
\end{eqnarray*}
so that $\Im(f_{n+1}+f_{n+2}+\cdots+f_m)(iy_n^-)<4$ for all $m>n,$ when 
$n>1$ is large enough.
(We have used in the 
last inequality the fact that the choice of 
$t_n$ so that
$\Im f_n(it_n)>n$ implies that $a_n>(2nt_n)/(1+t_n^2)>n/t_n,$ so that 
$1/n>1/(a_nt_n),$ and 
our choice that $a_n<a_{n-1}/2$, 
$n\ge1$).
Also,
$\Im(f_1+f_2+\cdots+f_{n-1}+f_n)(t_n)>\Im f_n(it_n)>n,$
for all $n\in\mathbb N$. As seen before, $\Re f_n(iy)=0$ for all $y\ge0.$

Now, it is easy to verify that $(1+t_nz)/(t_n-z)$
are uniformly bounded on, say, $i[0,1]$,
so that $\sum_{n=1}^\infty f_n$ is convergent and the limit $g$ is an
analytic self-map of the upper half-plane. We observe that $\Im
g(iy_n^-)\leq4+10/9<7$, $\Im g(it_n)\geq n$, and $\Re 
g(i[0,+\infty))=\{0\}$
for all $n\in\mathbb N$, while $\lim_{n\to\infty}
t_n=\lim_{n\to\infty}y_n^-=\infty.$ Thus, $g$ has no radial, hence
no nontangential, limit at infinity.

At the same time, since all $-f_n$s are Voiculecu transforms, so is $-g$.
\end{enumerate}

We next consider the case where
we do not restrict ourselves
to regularizing measures which are infinitely divisible.

\begin{Th}
Let $\mu$ be a Borel probability measure on $\R$
so that the function $h_\mu(z)=F_\mu
(z)-z$ satisfies the following property:
\begin{enumerate}
\item[{\rm (H)}]
For any $x\in\mathbb R\cup\{\infty\}$,
either $\sphericalangle h_\mu(x)$ does not
exist, or $\sphericalangle h_\mu(x)\in\mathbb C^+;$
\end{enumerate}
Then for any Borel probability measure $\nu$ on $\R$ which is not a point 
mass,
$\mu\bxp\nu$ is absolutely continuous
with respect to the Lebesgue measure  and has a positive
 analytic density with
respect to the Lebesgue measure.
\end{Th}
\begin{rmq} For $x=\infty,$ {\rm (H)} is in fact equivalent 
to the fact that either $\sphericalangle \phi_\mu(\infty)$
does not exist or $\sphericalangle \phi_\mu(\infty)$ belongs to $\C^+$.
Indeed, as observed in \cite{BercoviciVoiculescuIUMJ}, $z+\phi_\mu(z)$ belongs to a neighbourhood of 
infinity  for sufficiently large $z$ in some 
 truncated cone $\Gamma_{\alpha, M}$. 
Thus, by the definition of $\phi$ and $h$,
$h_\mu(z+\phi_\mu(z))=-\phi_\mu(z)$ for sufficiently large $z$ in such a 
cone. Since 
$z+\phi_\mu(z)$ tends nontangentially to infinity when $z$ tends nontangentially to infinity, 
$\sphericalangle \phi_\mu(\infty)$ exists iff $\sphericalangle h_\mu(\infty)$ exists, and they are equal.
The significance of this fact for our problem will be seen in the next subsection.
\end{rmq}

\begin{pr}
We claim first that $\sphericalangle F_{\mu\bxp\nu}(x)$ exists and 
belongs
to the upper half-plane for all $x\in\mathbb R$. Indeed, with the 
notations
from Lemma \ref{lem0},
by Theorem 3.3 of \cite{BA}, the nontangential limits of $\omega_1$ and
$\omega_2$ at $x$ exist. As observed in part (1) of Theorem 3.3 in 
\cite{BA},
if $\sphericalangle\omega_2(x)\in\mathbb C^+$, then the result is true.
Assume first that $\sphericalangle\omega_2(x)\in\mathbb R$.
Theorem \ref{lindel} guarantees that if $\sphericalangle 
F_{\mu\bxp\nu}(x)
=\sphericalangle(F_{\mu}\circ\omega_2)(x)$ exists, it must equal the
nontangential limit of $F_\mu$ in $\sphericalangle\omega_2(x)$, so
  the
nontangential limit of $h_\mu$ in $\sphericalangle\omega_2(x)$
exists, and hence, by (H), belongs to the upper half-plane.

Suppose  the
nontangential limit of $F_\mu$, and hence of $h_\mu$,
in $\sphericalangle\omega_2(x)$ does not exist.
Then
by Theorem \ref{lindel}, $h_\mu\circ\omega_2$ has no nontangential limit 
at $x
.$ But
by part (1) of Lemma \ref{lem0} and definition of $h_\mu$, we have
$$\sphericalangle\omega_1(x)=x+
\lim_{\stackrel{z\longrightarrow x}{\sphericalangle}}h_\mu(\omega_2(z)),$$
which implies that $\sphericalangle\omega_1(x)$ does not exist, 
contradicting
Theorem 3.3 of \cite{BA}.

The last possible case is when $\sphericalangle\omega_2(x)=\infty$.
As before, if $\sphericalangle (h_\mu\circ\omega_2)(x)$ exists, then
by Theorem \ref{lindel} it
must coincide with $\sphericalangle h_\mu(\sphericalangle\omega_2(x))=
\sphericalangle h_\mu(\infty)$, so,
by our hypothesis (H), it must belong to the upper half-plane. Thus,
$\sphericalangle\omega_1(x)\in\mathbb C^+,$ so that, by Theorem 3.3 in
\cite{BA}, $\sphericalangle  F_{\mu\bxp\nu}(x)\in\mathbb C^+.$
Assume now that $\sphericalangle (h_\mu\circ\omega_2)(x)$ does not exist,
so that
there exists an infinite set $W$ of points $c\in \mathbb C^+\cup\mathbb R$
for which there is a sequence $\{z_n^c\}_n$ converging to $x$ 
nontangentially
so that $\lim_{n\to\infty}(h_\mu\circ\omega_2)(z_n^c)=c$. But then
$$\sphericalangle\omega_1(x)=
\lim_{\stackrel{z\longrightarrow x}{\sphericalangle}}\omega_1(z)=
\lim_{n\to\infty}\omega_1(z_n^c)=\lim_{n\to\infty}z_n^c+h_\mu(\omega_2(z_n^c))
=x+c$$
for any $c\in W.$ This contradicts the existence of the nontangential 
limit
of $\omega_1$ at $x$.

This establishes the existence of nontangential limits of 
$F_{\mu\bxp\nu}$
at all points $x\in\mathbb R$ and the fact that $\sphericalangle 
F_{\mu\bxp\nu}(x)
\in\mathbb C^+$ for all $x\in\mathbb R$.
We claim that $\sphericalangle\omega_2(x)\in\mathbb C^+$. Indeed,
it is easy to see that  $\sphericalangle\omega_2(x)$ is finite, since 
otherwise
we would have, by Lemma \ref{lem0}, Theorem \ref{lindel} and Theorem 
\ref{Nevanlinna}, that
$\sphericalangle F_{\mu\bxp\nu}(x)=\sphericalangle 
(F_{\mu}\circ\omega_2)(x)
=\sphericalangle F_{\mu}(\sphericalangle\omega_2(x))=\sphericalangle 
F_{\mu}(\infty)=\infty,$
which is a contradiction. Thus, by Lemma \ref{lem0}, part (1), we have 
that
$\sphericalangle\omega_1(x)$ is also finite. Moreover, at least one of 
$\sphericalangle\omega_1(x)$, $\sphericalangle\omega_2(x)$
must then belong to the upper half-plane. Remark \ref{ImF>Imz} guarantees 
that in fact
both must be in the upper half-plane. Theorem 3.3 of \cite{BA} and Lemma 
\ref{19.09.06.8} concludes the proof.
\end{pr}

\subsection{A result of non existence of analytic densities}\label{negan}

\begin{propo}\label{negative}
Assume that $\mu$, a $\boxplus$-infinitely divisible probability so that
$\mu^{\boxplus t}$ has no atoms for some $t<1$
(the existence of $\mu^{\boxplus t}$ for $t<1$ is guaranteed by the
infinite divisibility of $\mu$), has finite second
moment. Then there exists a  probability measure $\nu$  on $\mathbb R$ so
that the density of
$\mu\bxp\nu$ is not analytic everywhere.
Moreover, the density of $\mu\bxp\nu$
vanishes at a point.
\end{propo}

\begin{rmq}
\begin{itemize}
\item
Note that the fact that the density
of $\mu\bxp\nu$ may easily  vanish 
inside the support of the measure 
 was already foreseen by
P. Biane in \cite{Biane2}, Proposition 6,
where he proved that if $\nu$ is a probability measure with
continuous strictly positive density on $]-\e,0[\cup ]0,\e[$
for some $\e>0$ such that
   $$\int x^{-2} d\nu(x)<\infty$$
and $\mu=\sigma_t$ is the semicircular variable
with covariance $t$, then $\frac{d\mu\bxp \sigma_t}{dx}(0)=0$
for $t>0$ small enough. Our proof extends this
phenomenon to any $\boxplus$-infinitely divisible probability measure
$\mu$
with finite second moment, under the 
(technical) hypothesis that
   $\mu^{\boxplus t}$ has no atoms for some $t<1$. 
\item
In the case where the probability measure $\mu$ has a finite second moment,
by theorem 1.3  and remark 1.1
of \cite{benaychIUMJ}, if
$m$, $v$ denote respectively the mean and the variance of $\mu$, one has, as
$z$ goes to infinity non tangentially, $\phi_\mu(z)=m+v/z+o(1/z)$, hence
the second hypothesis of theorem \ref{Thinfty} cannot be satisfied. However,
the existence of a finite second moment for $\mu$ has no incidence on the
first hypothesis of theorem \ref{Thinfty}. Indeed, as it will be explained in the proof of proposition
\ref{negative}, up to a translation, $\phi_\mu$ can be expressed as the
Cauchy transform of the finite positive measure $d\sigma(t)=(1+t^2)dG(t)$,
where $G$ is the L\'evy measure of $\mu$. Hence up to the addition of a real number,  the non tangential limit of
$\phi_\mu$ at any real number $x$ only depends on the restriction of $G$
to a neighborhood of $x$, which, by proposition 2.3 of \cite{benaychIUMJ},
is independent of 
 the existence of a finite second moment for $\mu$.
\end{itemize}
\end{rmq}

\begin{pr}
First note that by theorem 3.1 of \cite{BBercoviciAtoms}, the hypothesis
that $\mu^{\boxplus t}$ has no atoms for some $t<1$
implies that $\mu$ has no atoms. Hence by Remark \ref{rmqdensity}
or  Theorem 7.4 of
    \cite{BVReg}, for any $\nu$ which is not a point mass,
   $\mu\bxp\nu$ has a density. 

Observe  that if the density of $\mu\bxp\nu$ has a hole in the support
(meaning a nontrivial
interval on which it is zero), it cannot be analytic on $\mathbb R$ by the
identity principle.
Similarly, the set of zeros of the density must be discrete in 
$\mathbb
R$. Thus,
we may assume that $\mu$ satisfies these two conditions, since otherwise
we would readily obtain the probability
measure $\nu$ of the proposition
by taking $\nu=\delta_a$ for some $a\in\mathbb R.$

The strategy of the proof is as follows; we first show that
   $g:z\in\C^+\ra F_\nu(-\phi_\mu(z))$ has infinity as Denjoy-Wolff point
under a certain condition. Regarding $F_{\mu\bxp\nu}(0)$
as a fixed point of this map will guarantee that either
$F_{\mu\bxp\nu}(0)$ is infinite or belongs to $\R$.
We will show that under some hypothesis on $\nu$, it has to be infinite,
which will prove that $\frac{d\mu\bxp\nu}{dx}(0)=0$
and also that $\frac{d\mu\bxp\nu}{dx}$ is not analytic at
the origin. 

To study the
Denjoy-Wolff point of $g$, we first shall
write both $\phi_\mu$ and $F_\nu$
as (roughly speaking) Cauchy-Stieljes transforms of some measures
   $\sigma$ and $\rho$  on $\R$. 

   Recall (\cite{BercoviciVoiculescuIUMJ}) that there exists a real number
$\gamma$
   and a positive finite measure $G$ on the real line, called the {\it
L\'evy measure} of $\mu$ such that for all $z\in \C^+$,
$$\phi_\mu(z)=\gamma+\int_\mathbb R\frac{1+tz}{z-t}dG(t). $$
By proposition 2.3 of \cite{benaychIUMJ}, the finiteness of the second
moment of $\mu$ is equivalent
to the finiteness of the second moment of its L\'{e}vy measure $G$. Thus,
we can represent
the Voiculescu transform of $\mu$ as
$$\phi_\mu(z)=\gamma+\int_\mathbb R\frac{1+tz}{z-t}dG(t)=
\gamma+\int_\mathbb R(t+\frac{1+t^2}{z-t})dG(t)=\gamma'+\int_\mathbb
R\frac{d\sigma(t)}{z-t},\quad z\in\mathbb C^+,
$$
where $\gamma'\in\mathbb R$, $d\sigma(t)=(1+t^2)\times dG(t)$.
   Because $G$ has finite second
moment, $\sigma$ has finite mass.
   By a translation of $\mu$,
we may assume that $\gamma'=0$, so that $\phi_\mu$ is simply the Cauchy
transform of
the positive finite measure $\sigma$.

Let $\nu$ be a Borel probability on $\mathbb R$.  By Theorem
\ref{Nevanlinna}, we can write
   the reciprocal of
its Cauchy transform as
$$F_\nu(z)=a+z+\int_\mathbb R\frac{1+tz}{t-z}d\rho(t)=
a+z+\int_\mathbb R\frac{d\rho(t)}{t-z}+\int_\mathbb
R\frac{zt}{t-z}\,d\rho(t),$$
for all $z\in\mathbb C^+,$ where $a\in\mathbb R,$ and $\rho$ is a positive
finite
measure. We will show that if $\nu$ is so that
$\rho(\{0\})\geq\sigma(\mathbb R),$
then the density of $\mu\bxp\nu$
vanishes at the origin. So let $\nu$ satisfy this condition. We have
by Theorem \ref{infdiv+}(ii) 
\begin{equation}\label{aaa}
\lim_{y\to+\infty}\frac{F_\nu(-\phi_\mu(iy))}{iy}=
\lim_{y\to+\infty}\int_\mathbb R\frac{d\rho(t)}{iy(t+\phi_\mu(iy))}+
\lim_{y\to+\infty}\frac{-\phi_\mu(iy)}{iy}\int_\mathbb
R\frac{t\,d\rho(t)}{t+\phi_\mu(iy)}.
\end{equation}

Observe that $\phi_\mu(iy)$ approaches zero nontangentially
when $y$ tends to infinity. Indeed,
since $\lim_{y\to+\infty}iy\phi(iy)=\sigma(\mathbb R)>0,$ we have
$\lim_{y\to+\infty}y\Re\phi_\mu(iy)=0,\quad\lim_{y\to+\infty}y\Im\phi_\mu(iy)=-\sigma(
\mathbb R),$
so, given $0<\varepsilon<\sigma(\mathbb R)/2,$ there exists
$y_\varepsilon>1$ so that
for all $y\ge y_\varepsilon,$ we have $|y\Re\phi_\mu(iy)|<\varepsilon,$ $
|\sigma(\mathbb R)+y\Im\phi_\mu(iy)|<\varepsilon.$ Thus,
$$\frac{|\Im\phi_\mu(iy)|}{|\Re\phi_\mu(iy)|}=
\frac{|y\Im\phi_\mu(iy)|}{|y\Re\phi_\mu(iy)|}>\frac{\sigma(\mathbb
R)-\varepsilon}{
\varepsilon}>1,$$ for all $y\ge y_\varepsilon.$
Now,
$$\lim_{y\to+\infty}\int_\mathbb R\frac{t\,d\rho(t)}{t+\phi_\mu(iy)}  = 
\lim_{y\to+\infty}\int_\mathbb
R\left(1-\frac{\phi_\mu(iy)}{t+\phi_\mu(iy)}\right)d\rho(t)
 = \rho(\mathbb R)-\rho(\{0\}).$$
Since, by Theorem \ref{infdiv+} above,
$\lim_{y\to+\infty}\phi_\mu(iy)/iy=0,$ we conclude that the
second limit in the equation (\ref{aaa}) vanishes. 

On the other hand, if we denote $f_y(t)=\frac{1}{iy(t+\phi_\mu(iy))},$
$t\in\mathbb R,$
$y>1,$ then $\lim_{y\to+\infty}f_y(t)=\frac{1}{\sigma(\mathbb
R)}\chi_{\{0\}}(t)$
pointwise, where
$\chi_A$ is the characteristic function of $A$. Also,
$$|f_y(t)|^2=\frac{1}{y^2(t+\Re\phi_\mu(iy))^2+y^2(\Im\phi_\mu(iy))^2}\leq
\frac{1}{y^2(\Im\phi_\mu(iy))^2}<\frac{4}{\sigma(\mathbb R)^2},$$
for all $y>y_\varepsilon.$ So by the dominated convergence theorem,
$$\lim_{y\to+\infty}\int_\mathbb
R\frac{d\rho(t)}{iy(t+\phi_\mu(iy))}=\frac{\rho(\{0\})}{
\sigma(\mathbb R)}.$$


By \eqref{aaa}, we conclude that

$$\lim_{y\to+\infty}F_\nu(-\phi_\mu(iy))/iy=\frac{\rho(\{0\})}{
\sigma(\mathbb R)}\geq1,$$
which insures that the analytic function
$g\colon\mathbb C^+\to\mathbb C^+$ defined by
$g(z)=F_\nu(-\phi_\mu(z)),$ has infinity as its Denjoy-Wolff point.

We next show that this implies that
$F_{\mu\bxp\nu}(0)$ belongs to $\R\cup\{\infty\}$.
So, we suppose that
$F_{\mu\bxp\nu}(0)\in\mathbb C^+$ to get a contradiction
   ($F_{\mu\bxp\nu}$ extends continuously to $\R$ by Lemma \ref{lem0}
(4)).
   Note that by (\ref{cup.of.tea}), the relation
$F_\nu(z-\phi_\mu(F_{\mu\bxp\nu}(z)))=F_{\mu\bxp\nu}(z),$
   gives, by letting $z$ going to zero nontangentially,
\begin{equation}\label{13.09.06.paris.1}
F_\nu(-\phi_\mu(F_{\mu\bxp\nu}(0)))= F_{\mu\bxp\nu}(0).
\end{equation}
Thus, $F_{\mu\bxp\nu}(0)$ should be a fixed point
of $g$ in $\C^+$, and thus
its Denjoy-Wolff point; this is in contradiction
with the previous statment that
the Denjoy-Wolff point of $g$ is infinity.
Hence one has $F_{\mu\bxp\nu}(0)\in\mathbb \R\cup\{\infty\}.$

Observe that
$ F_{\mu\bxp\nu}(0)\neq 0$. Indeed, by Remark \ref{rmqdensity},
this 
equality 
would imply that $t\omega_2(0)$ is an atom of $\mu^{\boxplus t}$ for 
all
   $t<1$, contradicting the hypothesis.

Thus, $\Im F_{\mu\bxp\nu}(0)=0$ and $ F_{\mu\bxp\nu}(0)\neq 0$,
so that $\Im  G_{\mu\bxp\nu}(0)=0$.
Part (4) of Lemma \ref{lem0} tells us that $F_{\mu\bxp\nu}$  is
continuous on $\mathbb R$. In particular,
   since $F_{\mu\bxp\nu}(0)\neq 0$, $G_{\mu\bxp\nu}(x)$
will be continuous and finite for $x$ in some open interval
$I$ around zero. 
Lemma \ref{19.09.06.8} (i)
guarantees that ${\mu\bxp\nu}$ will have a continuous density on $I$ 
which vanishes  at zero.

We show below a more precise statment to
prove the breaking of analyticity
at the origin
when $\rho(\{0\})=\sigma(\R)$, namely that $F_{\mu\bxp
\nu}(0)=\infty$.

So,
assume now  that $\rho(\{0\})=\sigma(\mathbb R).$ Then we claim that in
fact
$
F_{\mu\bxp\nu}(0)=\infty$. Indeed, by Lemma
\ref{lem0} (2),
$\omega_1(0)$ 
is the Denjoy-Wolff point of the function $g_0(w)=-\phi_\mu(F_\nu(w))$.
We next show that this point must be the origin. 

Indeed,
observe that $F_\nu(iy)/iy$ goes to one as $y\in\R^+$ goes
to infinity. Hence, $F_\nu(iy)$ approaches infinity nontangentially
when $y\to+\infty$. Thus, since $\phi_\mu$
approaches zero nontangentially
at infinity,  $g_0(y)=-\phi_\mu(F_\nu(iy))$
converges to zero as $y$ goes to zero yielding 
$g_0(0)=0$. Also,
$$\lim_{y\to0}\frac{g_0(y)}{iy}=
\lim_{y\to0}\frac{-\phi_\mu(F_\nu(iy))}{iy}=\frac{\lim_{y\to0}\phi_\mu(F_\nu(iy))F_\nu(iy)}{\lim_{y\to0}-iyF_\nu(iy)}=
\frac{\sigma(\mathbb R)}{\rho(\{0\})}=1.$$
The Denjoy-Wolff theorem and the remarks following it 
imply that zero is the Denjoy-Wolff point for $g_0$, so 
by uniqueness of the Denjoy-Wolff point, 
$\omega_1(0)=0.$

We know that $F_\nu$ has infinite nontangential limit at zero  (because we
supposed that $\rho$ has an atom at zero), so this, coupled with the
equation $F_{\mu\bxp\nu} (z)=F_\nu(\omega_1(z))$, with the existence of
the limit of $F_{\mu\bxp\nu}$ at zero and Lindel\"{o}f's Theorem
\ref{lindel}, implies that
$F_{\mu\bxp\nu}(0)=\infty$.


Observe that $\omega_1$ is not analytic in zero.
Indeed, if  it were analytic, it would
have a finite derivative in zero. 
   However, with $H$ the function given in  Lemma \ref{lem0}(2), the
previous estimates show that
$H(0)=0$
and
$$H'(0)=\lim_{y\ra 0}\frac{H(iy)}{iy}=1-\lim_{y\to0}\frac{-\phi_\mu(F_\nu(iy))}{iy}=0$$
   implies, by
Proposition 4.7 (5)
in \cite{BBercoviciMult}, that
   $\lim_{y\to0}\omega_1(iy)/iy=\omega'_1(0)=1/H'(0)=\infty$.
As a consequence, 
 $G_{\mu\bxp \nu}$ is not differentiable
at the origin, since if it were
\begin{eqnarray*}
\lim_{y\to0}-\frac{G_{\mu\bxp\nu}(iy)-G_{\mu\bxp\nu}(0)}{iy-0} & =
&
-\lim_{y\to0}\frac{1}{iyF_\nu(\omega_1(iy))}\\
& = &
\lim_{y\to0}-\frac{1}{\omega_1(iy)F_\nu(\omega_1(iy))}\cdot\frac{\omega_1(iy)}{iy}.
\end{eqnarray*}
The second factor above has just been shown to converge to infinity. For
the first factor,
observe that $\sphericalangle\lim_{z\to0}zF_\nu(z)={\rho(\{0\})}.$ Since
$y\mapsto\omega_1(iy)$
is a smooth path in the upper half-plane ending at zero, Theorem \ref{521}
guarantees that there
exists a subsequence $y_n\to0$ so that
$\lim_{n\to\infty}\omega_1(iy_n)F_\nu(\omega_1(iy_n))=\rho(\{0\}).$
So the limit above either does not exist, or is infinite. In both cases,
we conclude that $G_{\mu\bxp\nu}$
is not differentiable at  zero. Thus by lemma \ref{19.09.06.8}, the
density of $\mu\bxp\nu$ is not analytic in zero.
\end{pr}

%


   { The} above proposition provides a large class of examples of free convolutions whose densities 
have cusps
in their support (points where
the density vanishes and is not analytic), and relates this phenomenon to the finiteness of
second moments. 
We show below that it is possible that
the density of $\mu\bxp\nu$ vanishes at a
point, but is still analytic.

\begin{propo}
Let $\mu$ be the semicircular distribution.
Then, there exists a probability measure
$\nu$ on $\R$ so that the density
   $\rho(x)=\frac{d\mu\bxp\nu}{dx}(x)$ vanishes at the origin, is strictly
positive on $]-\e,0[\cup]0,\e[$ for some $\e>0$
but is analytic at the origin.
\end{propo}
\begin{pr}
Let $\mu$ be the semicircular distribution,
so that $\phi_\mu(z)=1/2z,$ $z\in\overline{\mathbb C^+}.$ We claim that
$\nu$  given by its reciprocal Cauchy-Stieljes
$$F_\nu(z)=z+i-1+\frac{z-i}{z+i}-\frac{1}{z},\quad z\in\mathbb C^+,$$
will satisfy the properties of the proposition.

Indeed, with the notations from the proof of the previous proposition,
$$g_0(w)  = -\phi_\mu(F_\nu(w))
= -\frac{1}{2\left(w+i-1+\frac{w-i}{w+i}-\frac1w\right)}
 = -\frac{w(w+i)}{2w^3+4iw^2-4(1+i)w-2i},$$
for all $w\in\mathbb C^+,$ so in fact $g_0$ extends analytically around
zero, and moreover
$$g_0'(0)=\lim_{w\to0}g_0(w)/w=
1/2<1
\Rightarrow |g'_0(0)|<1\,
\mbox{ and }\, g_0(0)=0.$$ Thus, zero is the Denjoy-Wolff point of $g_0$,
and, by Lemma \ref{lem0} (2),
we conclude that $\omega_1(0)=0.$

We next show that $\omega_1$
extends analytically around the origin.
In fact,
the function $H(w)=w+\phi_\mu(F_\nu(w))$ has, by   Lemma \ref{lem0}.(2),
$\omega_1$ as right inverse. At the same time,
$H$ extends analytically around zero,
and $H'(0)=1-\frac{1}{2}=\frac{1}{2}\neq0$,
so $H$ is locally invertible around zero. The analyticity of $\omega_1$
around the origin  follows from the
implicit function theorem.

We now conclude that $G_{\mu\bxp \nu}$
vanishes at the origin,
is analytic in a neighborhood
of the origin and has negative
imaginary part in a neighborhood of the
origin (except at the origin itself); this
will prove the lemma according to Lemma \ref{19.09.06.8}.
Now, $F_\nu$ is meromorphic around zero, with a
simple pole at zero, so $G_\nu(0)=0$ and $G_\nu$ is
analytic around zero.
By Lemma \ref{lem0}, we have
$0=G_\nu(0)=G_\nu(\omega_1(0))=G_{\mu\bxp\nu}(0),$
and $G_{\mu\bxp\nu}(z)=
G_\nu(\omega_1(z))$ is analytic on a neighbourhood
of zero in $\mathbb C$. We finally show that
$G_{\nu\bxp\mu}$ has positive imaginary
part in $]-\e,0[\cup]0,\e[$ for some $\e>0$.
For that, since $G_{\nu\bxp\mu}(z)=G_\nu(\omega_1(z))$,
it is enough to show that $G_\nu(z)$ has  negative imaginary part
for $z$ so that $0<|z|<\e'$
(since $\omega_1$ is analytic and null at the origin).

But, a straightforward computation gives
$$G_\nu(r)  =  \frac{1}{F_\nu(r)}
 = \frac{r(r+i)}{r^3+2ir^2-2(i+1)r-i}
 = \frac{r(r+i)}{(r^3-2r)+i(2r^2-2r-1)},$$
so
$$
\Im G_\nu(r) = -\frac{r^2(r-1)^2}{r^2(r^2-2)^2+(2r^2-2r-1)^2}<0,$$
for all $r\in\mathbb R\setminus\{0,1\}.$

We notice that in fact we have only used, with the notations from the
previous proposition, the facts that
$\rho(\{0\})>\sigma(\mathbb R),$ that $G_\nu$ is analytic around zero, and
$\phi_\mu$ around infinity. Thus, a much larger class of
such pairs of measures $\mu,\nu$ provide an analytic density around zero
which is zero in zero.
\end{pr}

\begin{rmq}
Note that our 
construction of the example in the proposition above were based on the fact that $\rho(\{0\})>0$.
In that case,
$F_\nu(z)\approx_{z\ra 0} -\frac{1}{z}\rho(\{0\})$.
This is equivalent to the fact that $\int t^{-1} d\nu(t)=0$
and
$$\int \frac{1}{t^2} d\nu(t)=\frac{1}{\rho(\{0\})}.$$
Therefore,
if $\mu$ is an infinitely divisible measure with
finite second moment and positive density, and if
we denote again by $\sigma$ the finite
measure on $\R$
given by
$$\phi_\mu(z)=\gamma+\int_\mathbb R\frac{1+tz}{z-t}dG(t)=
\gamma+\int_\mathbb R(t+\frac{1+t^2}{z-t})dG(t)=\gamma'
+\int_\mathbb R\frac{d\sigma(t)}{z-t},\quad z\in\mathbb C^+,
$$
where $\gamma'\in\mathbb R$, $d\sigma(t)=(1+t^2)\times dG(t)$, we have
shown that we have three possibilities, which follow the intuition

\begin{enumerate}
\item If $\int \frac{1}{t^2} d\nu(t)<1/\sigma(\mathbb R)$,
so $\nu$ does not put much mass in the neighborhhood of the
origin,
$\mu\bxp\nu$ has a density which vanishes at the
origin but has a finite derivative at 
the origin. The proof was detailed above in a specific
example in the last proposition  but could be generalized.

\item If $\int \frac{1}{t^2} d\nu(t)=1/\sigma(\mathbb R)$,
which corresponds to a critical amount
of mass around the origin, the density has a cusp
at the origin (at least under the asumptions that $\mu^{\boxplus t}$
has no atoms for $t<1$).

\item If $\nu(\{0\})>0$,
   we have an analytic strictly positive density in
zero whenever $\mu^{\boxplus t}$ lacks atoms for all $t>0$.
   Indeed, if one assumes
   $\nu(\{0\})>0$, then $ \ds\lim_{z\to 0}\lf|\f{F_\nu(z)}{z}\ri|= \infty $
would imply $ \ds\lim_{z\to 0}zG_\nu(z)= 0 $, which is obviously false, by
(2) of lemma 2.17 of \cite{BA}.  Thus if $F_{\mu\boxplus\nu}(0)\in
\mathbb R$, 
then  $\omega_1(0)\in\mathbb C^+$ (by lemma \ref{lem0}), which is impossible, and if  $ F_{\mu\boxplus\nu}(0)=\infty$, then 
$F_{\mu\boxplus\nu}(0)$ is the    Denjoy-Wolff point of $g$ (equation (\ref{answer.alice.19.09.06.1})),  which is also impossible because
  $-\phi_\mu(\infty)=0$ and $F_\nu(0)=0$.
Thus,    $F_{\mu\bxp\nu}(0)$ must
belong to $\C^+$.
\end{enumerate}
\end{rmq}

\section{The rectangular case}
\subsection{Main result}
We shall fix $\lambda\in(0,1)$, and assume that all probability measures
are symmetric. We prove here an analogue to
Theorem \ref{Thinfty} for the rectangular convolution,
which says that the restriction to the upper half-plane of the function 
$H$ extends continuously to $\R^+$, analytically outside a closed
set of Lebesgue measure zero.
We shall see that this implies that $\mu\arc\nu$ admits an analytic density on the 
complement of that set.
  Unlike for the square case, 
we did not succeed to get rid of this closed negligible set
where the density
could stop being analytic. We however can give sufficient 
conditions so that the density is continuous  everywhere
(Corollary \ref{final.corollary.rectangular}).
Examples which satisfy our conditions are provided in section
\ref{10.01.06.80}. 
Moreover, as we shall discuss later,
the density often vanishes around the origin, which is, in the rectangular
setting,  a very specific point.
A consequence of this fact is that the full strength of Theorem \ref{Thinfty} cannot be
achieved in the rectangular case: given a $\arc$-infinitely divisible probability $\mu$,
there exists a symmetric probability measure $\nu\neq\delta_0$ so that the density of
$\mu\arc\nu$ is not everywhere analytic. 
We shall  study this phenomenon 
  in the last  paragraph.

Our main tool will be an ad-hoc subordination result for the functions $H$:

\begin{lem}\label{auxm}
Let $\mu,\nu$ be two symmetric probability
measures on $\mathbb R$. 
Assume that the rectangular $R$-transform $C_\mu$ of $\mu$ extends
analytically to $
\mathbb C\setminus\mathbb R^+$ (this happens for example if $\mu$ is $\arc
$-infinitely divisible - see Theorem \ref{10.01.06.2}). Then there exist two 
unique meromorphic functions $\omega_1,\omega_2$ on $\mathbb C\setminus
\mathbb R^+$ so that $H_\mu(\omega_1(z))=H_\nu(\omega_2(z))=H_{\mu\arc\nu}(z),
$ $\omega_j(\bar z)=\overline{\omega_j(z)},$ $z\in\mathbb C\setminus\mathbb R^+$, 
and $\lim_{x\uparrow0}\omega_j(x)=0,$ $j\in\{1,2\}$. Moreover, 
\begin{enumerate}
\item[{\rm(i)}] $\omega_2$ is injective and analytic on $\mathbb C\setminus\mathbb
R^+$; it is the right inverse of the meromorphic function $k(w)=\frac{H_\nu(w)}{T[C_\mu(
H_\nu(w))+M_\nu(w)]},$ $w
\in \mathbb C\setminus\mathbb R^+$;
\item[{\rm(ii)}] $\arg z\le\arg\omega_2(z)<\pi$, $z\in\mathbb C^+.$
\end{enumerate}

\end{lem}

\begin{pr}
There
exists an $\varepsilon>0$ so that for $z\in(-
\varepsilon,0)$, by taking $w=H_{\mu\arc\nu}(z)$ in relation
$C_{\mu\arc\nu}(w)=
C_\mu(w)+C_\nu(w),$ we have
$$U\left(\frac{H_{\mu\arc\nu}(z)}{z}-1\right)-C_\mu(H_{\mu\arc\nu}(z))=
U\left(\frac{H_{\mu\arc\nu}(z)}{H_\nu^{-1}(H_{\mu\arc\nu}(z))}-1\right).$$
Applying $T$ in both sides gives
$$T\left[U\left(\frac{H_{\mu\arc\nu}(z)}{z}-1\right)-C_\mu(H_{\mu\arc\nu}(z))
\right]
=\frac{H_{\mu\arc\nu}(z)}{H_\nu^{-1}(H_{\mu\arc\nu}(z))}.$$
By part (b) of section \ref{Remarks.about.H},
  $H_{\mu\arc\nu}(z)$ doesn't vanish on
$\C\bck \R^+$,
hence in the interval $(-\varepsilon, 0)$ 
where the previous equation is valid,
its left hand
term doesn't vanish. So on $(-\varepsilon, 0)$ we have
\begin{equation}\label{H}
H_{\mu\arc\nu}(z)=H_\nu\left(\frac{H_{\nu\arc\mu}(z)}{T\left[U\left(\frac{H_{\mu\arc\nu}(z)}{z}-1\right)-C_\mu(H_{\mu\arc\nu}(z))
\right]}\right).
\end{equation}
This equation holds for $z\in(-\varepsilon,0)$, and, by analytic
continuity, in all
points of the connected component of the domain
of analyticity  of the right hand term
which contains
$(-\varepsilon,0)$. Thus, if we denote $f(z,w)=\frac{w}{T( U(\frac{w}{z}-1)-C_\mu(w))},
$ and let $\omega_2(z)=f(z,H_{\mu\arc\nu}(z)),$ we have proved that $H_{\mu\arc\nu}
(z)=H_\nu(\omega_2(z))$ for $z$ in some domain containing the interval
$(-\varepsilon,0)$. We shall argue in the following that this equation can be extended
to all points of $\mathbb C\setminus\mathbb R^+$.

Note that since $C_\mu$ is analytic on $\C\bck \R^+$, by (b) of
the section \ref{Remarks.about.H}, $C_\mu(H_{\mu\arc\nu}(z))$ is defined 
on $\C\bck
\R^+$.
Moreover, by (a) of the section \ref{Remarks.about.H},
$z\mapsto U\left(\frac{H_{\mu\arc\nu}(z)}{z}-1\right)$
admits an analytic extension to $\C\bck \R^+$, denoted by
$M_{\mu\arc\nu}$. So any point
$z$ of $\C\bck\R^+$ which is in the boundary of the domain of the right
hand term of \eqref{H} satisfies either
$$\omega_2(z):=  \frac{H_{\nu\arc\mu}(z)}
{T\left[M_{\mu\arc\nu}(z)-C_\mu(H_{\mu\arc\nu}(z))
\right]}\in\R^+\, 
\mbox{ or }\, T\lf[M_{\mu\arc\nu}(z)-C_\mu(H_{\mu\arc\nu}(z))\ri]=0
.$$

$\bullet$ We first discuss the case when $\nu$ has the property that 
for any $x$
in $\mathbb R$, the Cauchy transform of
$\nu$ does not extend continuously to $x$.
This happens for example if $\nu$ is concentrated on a set of Lebesgue 
measure zero
and has support equal to $\mathbb R$, according to Lemma \ref{s}.

Consider the connected component of the domain of the right hand side of
(\ref{H}) that contains $(-\varepsilon,0).$ Assume first that
$z_0\in\mathbb C\setminus
\mathbb R^+,$ and yet $z_0$ is in the boundary of this component, which
implies that
either $\omega_2(z_0)\in[0,+\infty)$ (
because of part (b) of section \ref{Remarks.about.H}), or
$T\lf[M_{\mu\arc\nu}(z_0)-C_\mu(H_{\mu\arc\nu}(z_0))\ri]=0.$ The functions
$H_{\mu\arc\nu}$ and $\omega_2$ are
analytic and, respectively, meromorphic, in $z_0$
 and
$H_{\mu\arc\nu}(z_0)\in\mathbb C\setminus\mathbb R^+$.

Observe that if the first situation occurs, 
there must exist a whole (1-dimensional) analytic connected variety $V$
given by
the relation $\omega_2(.)\in(0,+\infty)$ to which $z_0$
belongs because $\omega_2$ is open as a meromorphic function.
But now for any point $\zeta\in V$, we will have that
$$\lim_{z\to\zeta}H_{\mu\arc\nu}(z)=\lim_{z\to\zeta}H_\nu\left(\frac{H_{\nu\arc\mu}(z)}{T\left[M_{\mu\arc\nu}(z)
-C_\mu(H_{\mu\arc\nu}(z))
\right]}\right)=\lim_{z\to\zeta}H_\nu(\omega_2(z)).$$
The left hand-side exists always and equals $H_{\nu\arc\mu}(\zeta)$, while 
the right hand side cannot
exist at least for a set of second Baire
category - Theorem \ref{Th.4.8.CL}. Specifically, $\omega_2(V)$ must be (by the identity theorem for
analytic functions)
a nontrivial interval in $(0,+\infty)$; for any $r\in \omega_2(V)$, we have a
$\zeta\in V\subset\mathbb C\setminus\mathbb R_+$ so that $\omega_2(\zeta)=r$, and
so, by Lindel\"{o}f's
theorem \ref{lindel}, since
$H_{\mu\arc\nu}(\zeta)=\lim_{z\to\zeta}H_{\nu}(\omega_2(z)),$ we have
$$H_{\mu\arc\nu}(\zeta)=\lim_{z\to\zeta}H_{\nu}(\omega_2(z))=\lim_{\stackrel{w\longrightarrow
r}{\sphericalangle}}H_\nu(w).$$
Theorem \ref{Th.4.8.CL}, together with the equation above, implies that
there is a point of $\omega_2(V)$ where the cluster set of $H_\nu$
is a single point. Hence by (\ref{H->G}), we have a contradiction with the
fact that  for any $x$
in $\mathbb R$, the Cauchy transform of $\nu$ does not extend continuously
to $x$.

Assume now that $z_0$ is so that
$T\lf[M_{\mu\arc\nu}(z_0)-C_\mu(H_{\mu\arc\nu}(z_0))\ri]=0$. 
Observe that since $z_0$ is assumed to be in $\mathbb C\setminus\mathbb
R^+,$
$H_{\mu\boxplus_\lambda\nu}(z_0)\in\mathbb C\setminus\mathbb R^+$ by
the section \ref{Remarks.about.H}, (b).
Hence the function
$\omega_2$ is meromorphic on a neighbourhood of
$z_0$, with a pole
at $z_0$. We conclude that $\lim_{z\to
z_0}\frac{H_{\nu\arc\mu}(z)}{T\left[M_{\mu\boxplus_\lambda\nu}(z)-C_\mu(H_{\mu\arc\nu}(z))
\right]}=\infty.$
Consider now a small enough ball $W\subset\mathbb C\setminus\mathbb R^+$
around $z_0$ so that $z_0$ is the only pole of
$\omega_2$ in $W$, and consider a connected component of the intersection of
this ball with the domain of the
function in the right hand side of (\ref{H}). Clearly $\omega_2(W)$ is a
neighbourhood of infinity and $W$
will contain $p\geq1$ analytic varieties that are mapped by $\omega_2$ onto
$(-\infty,-M)\cup(N,+\infty),$ for some
large enough $M,N>0$ ($p$ is the order of the pole at $z_0$). The
preimages of $(N,+\infty)$ divide $W$ into $p$
sectors. By (\ref{H}), if a point in one of these sectors belongs to the
connected component of the domain of the right hand term
of (\ref{H}) which contains $(-\varepsilon,0)$ then all that sector will
belong to it, and the two (distinct or not! - it might be a slit circle, if
$p=1$)
boundaries of the sector are mapped inside $\mathbb R^+.$ Thus, we are
reduced to the previous case.

We conclude that $\omega_2(z)=f(z,H_{\mu\arc\nu}(z))$ maps $\mathbb C\setminus
\mathbb R^+$ into itself.

$\bullet$ We now generalize the previous result to
any probability measure $\nu$. To this end, 
 we can approximate in the weak topology arbitrary symmetric
probabilities $\nu$
with probabilities concentrated on a set of zero Lebesgue measure and 
which have
total support, according to Lemma \ref{d}.
Since $\arc$ is continuous,
$H_{\mu\arc\nu_n}$ converges  to $H_{\mu\arc \nu}$,
and hence  $f(z,H_{\mu\arc\nu_n}(z)):\C\backslash \R^+\ra \C\backslash
\R^+$
converges to $f(z,H_{\mu\arc\nu}(z))$.
This implies that either $f(z,H_{\mu\arc\nu}(z))$
takes values also in $\mathbb C\setminus\mathbb R^+,$ or
it is constant. { But it cannot be constant, by equation (\ref{H}) and by
the fact that $H_{\mu\arc \nu}(z)$ is equivalent to $ z$ as $z$ tends to
zero in $\C\bck\R^+$ (see \cite{benaych.rectangular}, proposition 4.1).
Hence $f(z,H_{\mu\arc\nu}(z))$
takes values  in $\mathbb C\setminus\mathbb R^+$.} This proves that $\omega_2$
maps $\mathbb C\setminus\mathbb R^+$ into itself, and thus the equation
$H_\nu\circ\omega_2=H_{\mu\arc\nu}$ holds on $\mathbb C\setminus\mathbb R^+$.

$\bullet$ We finally show  that  $\omega_2$
satisfies the announced properties. 

First, 
it follows immediately from the definition of $\omega_2$
 and part (d) of subsection 
\ref{Remarks.about.H} that $\omega_2(\bar z)=\overline{\omega_2(z)}$ for all $z\in
\mathbb C\setminus\mathbb R^+.$ The uniqueness 
of $\omega_2$ on $(-\varepsilon,0)$ 
and the analyticity of $\omega_2$ 
in $\mathbb C\setminus\mathbb R^+$ proved above 
shows that $\omega_2$ is uniquely defined on all  $\mathbb C\setminus\mathbb R^+$.

We prove next properties (i) and (ii). As shown above, $H_\nu\circ\omega_2=H_{\mu
\arc\nu},$ so that for $\varepsilon>
0$ small enough, by subsection \ref{Remarks.about.H} (a), $M_\nu(\omega_2(x))=
C_\nu(H_\nu(\omega_2(x)))=C_\nu(H_{\mu\arc\nu}(x)),$ $x\in(-\varepsilon,0).$
Thus, $T\left[C_\mu(H_\nu(\omega_2(x)))+M_\nu(\omega_2(x))\right]=
T\left[C_\mu(H_{\mu\arc\nu}(x))+C_\nu(H_{\mu\arc\nu}(x))\right]=T\left[C_{\mu\arc\nu}
(H_{\mu\arc\nu}(x))\right]=T\left[M_{\mu\arc\nu}(x)\right].$ Now it follows immediately 
from the definition of 
$k$  in Lemma \ref{auxm} and subsection \ref{Remarks.about.H} (a) that $k(\omega_2
(z))=z$ for $z\in(-\varepsilon,0)$, $\varepsilon>0$ small enough, and  by analytic 
continuation, for $z\in\mathbb C\setminus\mathbb R^+$. This proves (i).

Let us recall that $\displaystyle\lim_{x
\uparrow0}\omega_2(x)=0$ and, by the definition of $T$ and properties of the function 
$H$, $\displaystyle
\lim_{x\uparrow0}\frac{\omega_2(x)}{x}=[T(0)]^{-1}\lim_{x\uparrow0}\frac{H_{\mu\arc\nu}(x)}{x}
=1$. Also, for $\varepsilon>0$ small enough, $\omega_2((-\varepsilon,0))\subseteq(
-\infty,0).$ Thus, the derivative of the
analytic function $\omega_2$ on the interval $(-\varepsilon,0)$ is positive for $
\varepsilon>0$
small enough, and so there is a small enough cone $\mathfrak K$
with vertex at zero and bisected by the negative half-line so that $\omega_2(\mathfrak
K\cap\mathbb C^+)\subseteq\mathbb C^+$ and $\omega_2(\mathfrak K\cap\mathbb
C^-)\subseteq\mathbb C^-$. Clearly, since $\omega_2(\mathbb C\setminus\mathbb
R^+)\subseteq\mathbb C\setminus\mathbb R^+,$ $\omega_2(\mathbb  C^+)\not
\subseteq\mathbb C^+$  implies that there exists a point $z_0\in\mathbb C^+$ with the 
property that $\omega_2(z_0)\in(-\infty,0).$ Assume such a point exists. Then from the 
equation (\ref{H}) and subsection \ref{Remarks.about.H} (d) we obtain that $H_{\mu
\arc\nu}(z_0)\in(-\infty,0)$, so that $T(M_{\mu\arc\nu}(z_0)-C_\mu(H_{\mu\arc\nu}
(z_0)))>0$. As observed in Remark \ref{T}, this requires that
$M_{\mu\arc\nu}(z_0)-C_\mu(H_{\mu\arc\nu}(z_0))\in\mathbb R.$ Since, by the
same subsection \ref{Remarks.about.H} (d), we have $C_\mu(\mathbb R^-)
\subseteq\mathbb R$, it follows that $M_{\mu\arc\nu}(z_0)\in\mathbb R$. 
But then, according to subsection \ref{Remarks.about.H} (a), $H_{\mu\arc\nu}(z_0)
=z_0T(M_{\mu\arc\nu}(z_0))\not\in\mathbb R$, a contradiction. We have now proved 
that $\omega_2$ preserves half-planes, and thus $\arg\omega_2(z)<\pi$ for $z\in
\mathbb C^+$.

Next we show that $\omega_2$ increases the argument. 
It is known from 
Theorem \ref{Nevanlinna} that 

$$\omega_2(z)=a+bz+\int_\R\frac{1+tz}{t-z}d\rho(t), z\in\mathbb C^+$$
for some $a\in\mathbb R,$ $b\ge0$ and positive finite measure $\rho$ on the real 
line.  
Since $\omega_2((-\infty,0))\subseteq(-\infty,0) $ and $\omega_2$ is analytic on 
the negative half-line, $\rho$ must be supported on $\mathbb R^+$. Moreover,
$0=\lim_{x\uparrow0}\omega_2(x)=a+\lim_{x\uparrow0}\int_\R\frac{1+tx}{t-x}d\rho(t)
=\int_\mathbb R\frac1td\rho(t)+a.$ Thus,
$a=-\int_\mathbb R\frac1td\rho(t).$ We conclude that

\begin{eqnarray*}
\omega_2(z) & = & a+bz+\int_\R\frac{1+tz}{t-z}d\rho(t)=  bz+\int_{\mathbb R^+}\left(\frac{1+tz}{t-z}-\frac1t\right)d\rho(t)\\ 
& = & z\left(b+\int_{\mathbb R^+}\frac{t^2+1}{t(t-z)}d\rho(t)\right).
\end{eqnarray*}
It is trivial to see that the factor in the parenthesis above maps $\mathbb C^+$ into
itself. Thus, $\arg\omega_2(z)\ge\arg z.$ This proves (ii).
Let us define 

$$\omega_1(z)=\frac{H_{\mu\arc\nu}(z)}{T[M_{\mu\arc\nu}(z)-M_\nu(\omega_2(z))]},
\quad z\in\mathbb C\setminus\mathbb R^+.$$

This function is obviously defined and meromorphic on $\mathbb C\setminus\mathbb
R^+$ and analytic continuation shows immediately that for $z\in(-\varepsilon,0)$

$$\omega_1(z)=\frac{H_{\mu\arc\nu}(z)}{T[M_{\mu\arc\nu}(z)-C_\nu(
H_{\mu\arc\nu}(z))]},$$
so that, as for $\omega_2$, $H_\mu(\omega_1(z))=H_{\mu\arc\nu}(z)$. 
This equality obviously extends by analytic continuation to $\mathbb C\setminus
\mathbb R^+$. However, we do not exclude the possibility that $H_\mu$ has 
an analytic continuation through the positive half-line that does not coincide with
the one provided by the formula in the subsection \ref{def-rect}.

It follows easily from the definition of $\omega_1$ and subsection
\ref{Remarks.about.H} that $\overline{\omega_1(z)}=\omega_1(\bar z)$ for all
$z\in\mathbb C\setminus\mathbb R^+$ and $\lim_{x\uparrow0}\omega_1(x)=0.$
The uniqueness of $\omega_1$ is determined by the same argument as in the case of 
$\omega_2$.
\end{pr}

Next, we study the boundary behaviour of the restriction of the subordination function 
$\omega_2$ to the upper half-plane.
\begin{lem}\label{bdry-subord-rect}
Let $\mu,\nu$ and $\omega_2$ be as in the Lemma \ref{auxm}. Then $\omega_2
|_{\mathbb C^+}$ extends continuously to $(0,+\infty).$
\end{lem}
\begin{pr}
Throughout the proof we will consider only $\omega_2
|_{\mathbb C^+}$ and we will denote it as $\omega_2$.
Assume that $r\in(0,+\infty)$ is so that the cluster set
$C(\omega_2,r)$  of $\omega_2$ at $r$
 is nontrivial, and hence, by Lemma \ref{lem1.2} an uncountably 
infinite closed connected subset of $\mathbb C^+\cup\mathbb R\cup\{\infty\}$.
Consider first the case when there exists at least one element 
$c\in C(\omega_2,r)\cap(\mathbb C^+
\cup(-\infty,0))$, and thus, by connectivity of $C(\omega_2,r)$, infinitely many. Fix such 
a point $c$, and let $\{z_n^{(c)}\}_{n\in\mathbb N}\subseteq\mathbb C^+$ be a 
sequence with the property that $\lim_{n\to\infty}z_n^{(c)}=r$ and $\lim_{n\to\infty}
\omega_2(z_n^{(c)})=c$. Passing to the limit in the equation 
$k(\omega_2(z_n^{(c)}))=z_n^{(c)}$, where $k$ is the function from Lemma 
\ref{auxm}, provides $k(c)=r$ for all $c\in C(\omega_2,r)\cap(\mathbb C^+
\cup(-\infty,0))$, and hence, by analytic continuation, for all $c\in\mathbb C\setminus
\mathbb R^+.$ This implies that $k$ is the constant function $r$, an obvious 
contradiction to Lemma \ref{auxm}. 

If $C(\omega_2,r)\subseteq\mathbb R^+\cup\{\infty\},$  then $C(\omega_2,r)\cap
\mathbb R^+$ must be a nontrivial closed interval, by Lemma \ref{lem1.2}. As $\omega_2$ maps $\mathbb C^+$ into itself, for all $c\in
C(\omega_2,r),$ with the possible exception of two points, there exists a sequence
$\{z_n^{(c)}\}_{n\in\mathbb N}\subseteq\mathbb C^+$ so that 
$\lim_{n\to\infty}z_n^{(c)}=r$, $\lim_{n\to\infty}
\omega_2(z_n^{(c)})=c$ and $\Re\omega_2(z_n^{(c)})=c$.  
As shown in subsection \ref{Remarks.about.H} (b), $H_\nu$ has nontangential limits
at almost all points of $\mathbb R^+$, and by Lemma \ref{Fatou}, so do $M_\nu$ and
$C_\mu$. Thus, $k$ must have nontangential limits at almost all points of $\mathbb 
R^+$. We have obtained that for Lebesgue-almost all points $c\in C(\omega_2,r),$

$$r=\lim_{n\to\infty}z_n^{(c)}=\lim_{n\to\infty}k(\omega_2(z_n^{(c)}))=
\lim_{\stackrel{w\longrightarrow c}{\sphericalangle}}k(w),$$
so that $k$ has constant nontangential limit $r$ on a set of nonzero Lebesgue 
measure. This, according to Theorem \ref{26.03.07.1}, implies that $k$ is the
constant function $r$, providing the same contradiction as before.

Thus, $\omega_2|_{\mathbb C^+}$ extends continuously to $(0,+\infty).$
\end{pr}

Now we are ready to prove our first continuity result.

\begin{propo}\label{10.01.06.5}


Let $\mu,\nu$ be two symmetric probability measures on $\mathbb R$, $\mu\neq\delta_0$. Assume that $\mu$ is $\arc$-infinitely divisible. Then for any $x\in(0,+\infty)$, the limits 

$$\lim_{z\to x,z\in\mathbb C^+}M_{\mu\arc\nu}(x)\quad{\rm and}\quad\lim_{z\to x,z\in\mathbb C^+}H_{\mu\arc\nu}(x)$$
exist in $\mathbb C\cup\{\infty\}$. The first limit belongs to ${\mathbb C^+\cup\mathbb
R\cup\{\infty\}}.$

\end{propo}

\begin{pr}
We will
follow idea from the proof of the theorem \ref{Thinfty}. 
Assume that $r\in(0,+\infty)$ is so that $C(M_{\mu\arc\nu},r)$ is nontrivial. Consider first the case when $ C(M_{\mu\arc\nu},r)\cap\mathbb C^+\neq\varnothing$, and thus, by Lemma \ref{lem1.2}, is uncountably infinite. Let $c\in C(M_{\mu\arc\nu},r)\cap\mathbb C^+$ and $\{z_n^{(c)}\}_{n\in\mathbb N}\subset
\mathbb C^+$ be so that $\lim_{n\to\infty} z_n^{(c)}=r$ and $\lim_{n\to\infty}
M_{\mu\arc\nu}( z_n^{(c)})=c$. We know from Lemma \ref{bdry-subord-rect} that $\omega_2(r):=\lim_{z\to r}\omega_2|_{\mathbb C^+}(z)$ exists in $\overline{\mathbb C^+}.$ Using the definition of $\omega_2$ and subsection \ref{Remarks.about.H} (a), we have:

\begin{eqnarray*}
\omega_2(r)=\lim_{n\to\infty}\omega_2(z_n^{(c)}) & = & 
\lim_{n\to\infty}\frac{H_{\mu\arc\nu}( z_n^{(c)})}{T\left[M_{\mu\arc\nu}( z_n^{(c)})-C_\mu(H_{\mu\arc\nu}( z_n^{(c)}))\right]}\\
& = & \lim_{n\to\infty}\frac{ z_n^{(c)}T(M_{\mu\arc\nu}( z_n^{(c)}))}{T\left[M_{\mu\arc\nu}( z_n^{(c)})-C_\mu(z_n^{(c)}T(M_{\mu\arc\nu}( z_n^{(c)})))\right]}\\
& = & \frac{rT(c)}{T[c-C_\mu(rT(c))]},\quad c\in C(M_{\mu\arc\nu},r)\cap\mathbb C^+.
\end{eqnarray*}
Thus, the meromorphic function $g_r\colon\mathbb C^+\cup(-1/\lambda,-1)\cup\mathbb C^-\to\mathbb C\cup\{\infty\},$ given by $g_r(c)=\frac{T(c)}{T[c-C_\mu(rT(c))]}$ is, by analytic continuation, constant, equal to $\omega_2(r)/r.$ It is trivial to observe that this implies $\omega_2(r)\not\in\{0,\infty\}.$ 

We shall express $C_\mu$ as a function of $s=T(c)$ to obtain a contradiction. Indeed, consider $c\in\left(-\frac{1}{2\lambda}-\frac12,-1\right).$ Then $s=T(c)$ if and only if 
$c=\frac{-1-\lambda+\left[(1-\lambda)^2+4\lambda s\right]^{1/2}}{2\lambda},$ $s\in\left(\frac{1}{4}(2-\lambda-\lambda^{-1}),0\right)$ (recall the notations from section 2.3.1.) 
Thus,
\begin{equation}\label{star}
\frac{rs}{\omega_2(r)}=T\left[\frac{-1-\lambda+\left[(1-\lambda)^2+4\lambda s\right]^{1/2}}{2\lambda}-C_\mu(rs)\right].
\end{equation}
As it is known that $C_\mu((-\infty,0))\subseteq(-\infty,0)$ and $\lim_{x\uparrow0}
C_\mu(x)=0$, we conclude that for $s\in\left(-\frac{(1-\lambda)^2}{4\lambda},0\right)$ 
close enough to zero,

$$T\left[\frac{-1-\lambda+\left[(1-\lambda)^2+4\lambda s\right]^{1/2}}{2\lambda}-C_\mu
(rs)\right]\in \R,$$ so that $\omega_2(r)\in \R\bck\{0\}$. Thus, since $\lim_{x\uparrow0}C_\mu(x)=0$, 
(\ref{star}) is equivalent to

$$C_{\mu}(rs)=\frac{-1-\lambda+\left[(1-\lambda)^2+4\lambda s\right]^{1/2}}{2\lambda}- 
\frac{-1-\lambda+\left[(1-\lambda)^2+4\lambda \frac{rs}{\omega_2(r)}\right]^{1/2}}{2
\lambda}.$$ But this implies either that $\omega_2(r)=r,$ so that $C_\mu(s)=0$ and 
thus $\mu=\delta_0$, or that $C_\mu$ is not analytic in the point $-\frac{r(1-\lambda
)^2}{4\lambda}\in(-\infty,0)$, an obvious contradiction.

Now consider the case when $C(M_{\mu\arc\nu},r)\subseteq \mathbb R\cup\{\infty\}.$ 
By subsection \ref{Remarks.about.H} (a) and Remark \ref{T}, in this case $C(H_{\mu
\arc\nu}, r)\subseteq\left[-\frac{r(1-\lambda)^2}{4\lambda},+\infty\right]$ is a nontrivial 
interval. As in the proof of Lemma \ref{bdry-subord-rect}, for any $d\in C(H_{\mu\arc\nu}, 
r)\setminus\{-\frac{r(1-\lambda)^2}{4\lambda},\infty\}$, with the possible exception of two points, there exists a 
sequence $\{z_n^{(d)}\}_{n\in\mathbb N}\subset\mathbb C^+$ so that $\lim_{n\to\infty} 
z_n^{(d)}=r,$ $\lim_{n\to\infty}H_{\mu\arc\nu}(z_n^{(d)})=d$ and $\Re 
H_{\mu\arc\nu}(z_n^{(d)})=d,$ $n\in\mathbb N$.

Let us observe that, by subsection \ref{Remarks.about.H} (a) and (c), we have

$$\lim_{n\to\infty}M_{\mu\arc\nu}(z_n^{(d)})\in\left\{\frac{-1-\lambda\pm\left[(1-\lambda
)^2+4\lambda\frac{d}{r}\right]^{1/2}}{2\lambda}\right\},$$
where we have the sign plus when $H_{\mu\arc\nu}(z_n^{(d)})$ tends to $d$ from
$\mathbb C^+$, and the sign minus when $H_{\mu\arc\nu}(z_n^{(d)})$ 
tends to $d$ from
$\mathbb C^-$. By dropping if necessary to a subsequence, we may assume that
$\lim_{n\to\infty}M_{\mu\arc\nu}(z_n^{(d)})$ exists. It is clear from the definition of 
the cluster set and the above considerations that $C(H_{\mu\arc\nu}, r)\setminus
\{-\frac{r(1-\lambda)^2}{4\lambda},\infty\}\subseteq A_+\cup A_-$, where 
$$A_+=\left\{d\in\mathbb R\setminus
\{\frac{r(1-\lambda)^2}{-4\lambda},\infty\}\colon
\begin{array}[t]{l}
\displaystyle \exists\{z_n^{(d)}\}_{n\in\mathbb N}\subseteq\mathbb C^+{\textrm{ so that}}\lim_{n\to\infty} z_n^{(d)}=r,\Re 
H_{\mu\arc\nu}(z_n^{(d)})=d, \\
\displaystyle\lim_{n\to\infty}H_{\mu\arc\nu}(z_n^{(d)})=d,H_{\mu\arc\nu}(z_n^{(d)})\in\mathbb C^+ 
\end{array}\right\},$$
$$A_-=\left\{d\in\mathbb R\setminus
\{\frac{r(1-\lambda)^2}{-4\lambda},\infty\}\colon
\begin{array}[t]{l}
\displaystyle \exists\{z_n^{(d)}\}_{n\in\mathbb N}\subseteq\mathbb C^+{\textrm{ so that}}\lim_{n\to\infty} z_n^{(d)}=r,\Re 
H_{\mu\arc\nu}(z_n^{(d)})=d, \\
\displaystyle\lim_{n\to\infty}H_{\mu\arc\nu}(z_n^{(d)})=d,H_{\mu\arc\nu}(z_n^{(d)})\in\mathbb C^- 
\end{array}\right\}.$$
are two (not necessarily disjoint) sets. Thus, at least one of $A_+,A_-$ has nonzero 
Lebesgue measure. Denote $C_\mu^+$ the restriction of $C_\mu$ to the upper 
half-plane and $C_\mu^-$ the restriction of $C_\mu$ to the lower half-plane.

Assume first that $A_+$ has nonzero Lebesgue measure. Then again
\begin{eqnarray*}
\omega_2(r)=\lim_{n\to\infty}\omega_2(z_n^{(d)}) & = & 
\lim_{n\to\infty}\frac{H_{\mu\arc\nu}( z_n^{(d)})}{T\left[M_{\mu\arc\nu}( z_n^{(d)})-C_\mu(H_{\mu\arc\nu}( z_n^{(d)}))\right]}\\
& = & \frac{d}{T\left[\frac{-1-\lambda+\left[(1-\lambda
)^2+4\lambda\frac{d}{r}\right]^{1/2}}{2\lambda}-\displaystyle\lim_{\stackrel{w\longrightarrow d}{\sphericalangle}}C^+_\mu(w)\right]},\quad d\in A_+.
\end{eqnarray*}
By the Riesz-Privalov Theorem we obtain again that 

$$\omega_2(r)T\left[\frac{-1-\lambda+\left[(1-\lambda
)^2+4\lambda\frac{d}{r}\right]^{1/2}}{2\lambda}-C_\mu^+(d)\right]=d,\quad d\in\mathbb
C^+.$$ 
Recalling that $C_\mu$ extends analytically to the negative half-line and considering 
values of $d\in(-\infty,0)$ close enough to zero, we observe as before that 
$\omega_2
(r)\in\mathbb R\backslash\{0\}$ and by analytic continuation
$$C_\mu(d)=\frac{\left[(1-\lambda
)^2+4\lambda\displaystyle\frac{d}{r}\right]^{1/2}}{2\lambda}-\frac{\left[(1-\lambda
)^2+4\lambda\displaystyle\frac{d}{\omega_2(r)}\right]^{1/2}}{2\lambda}, $$
providing the same contradiction as in the previous case.

Assume next that $A_-$ has nonzero Lebesgue measure, so that 
\begin{eqnarray*}
\omega_2(r)=\lim_{n\to\infty}\omega_2(z_n^{(d)}) & = & 
\lim_{n\to\infty}\frac{H_{\mu\arc\nu}( z_n^{(d)})}{T\left[M_{\mu\arc\nu}( z_n^{(d)})-C_\mu(H_{\mu\arc\nu}( z_n^{(d)}))\right]}\\
& = & \frac{d}{T\left[\frac{-1-\lambda-\left[(1-\lambda
)^2+4\lambda\frac{d}{r}\right]^{1/2}}{2\lambda}-\displaystyle\lim_{\stackrel{w\longrightarrow d}{\sphericalangle}}C^-_\mu(w)\right]},\quad d\in A_-.
\end{eqnarray*}
Exactly as for $A_+$, we obtain that $\omega_2(r)\in\mathbb R\backslash\{0\}$, and, from the Riesz-Privalov
Theorem, the formula
$C_\mu(d)=\frac{-\left[(1-\lambda)^2+4\frac{\lambda}{r}\right]^{1/2}+
\left[(1-\lambda)^2+4\frac{\lambda}{\omega_2(r)}\right]^{1/2}}{2\lambda}$, which
provides again the same contradiction.

Thus, we have established that the limit
$$M_{\mu\arc\nu}(x)=\lim_{z\to x,z\in\mathbb C^+}M_{\mu\arc\nu}(z)$$
exists for any $x\in(0,+\infty).$ The existence of the similar limit for $H_{\mu\arc\nu}$
follows immediately from subsection \ref{Remarks.about.H} (a), and since
$M_{\mu\arc\nu}(\mathbb C^+)\subseteq\mathbb C^+\cup\{0\},$ it follows that
$M_{\mu\arc\nu}(x)\in\mathbb C^+\cup\mathbb R\cup\{\infty\}.$
\end{pr}

\begin{cor}
Under the assumptions of Proposition \ref{10.01.06.5}, the absolutely continuous part
(with respect to the Lebesgue measure) of $\mu\arc\nu$ is continuous 
outside a closed set of zero Lebesgue measure, and its singular
part, if it exists, is 
supported 
on a closed subset of $\mathbb R$ of zero Lebesgue 
measure. 
\end{cor}

\begin{pr}
Recall that, by part (1) of Lemma 2.17 in \cite{BA}, the nontangential limit of the Cauchy 
transform of $\mu\arc\nu$ is infinite for almost all points in the support of the singular
part of $\mu\arc\nu$. Thus, using Proposition \ref{10.01.06.5} and the equality
$M_{\mu\arc\nu}(z)=\frac{1}{\sqrt{z}}G_{\mu\arc\nu}\left(\frac{1}{\sqrt{z}}\right)-1$
from subsection \ref{Remarks.about.H} (a), we can state that the support of the singular part of $\mu\arc\nu$ is concentrated on 
$S=S^+\cup S^-\cup\{0\}$, where $S^+=\{x\in(0,+\infty)\colon M_{\mu\arc\nu}(1/x^2)=
\infty\}$ and $S^-=-S^+.$

By the Riesz-Privalov Theorem (Theorem \ref{26.03.07.1}) it follows that the set
$S^+$ must be of zero Lebesgue measure, and by Proposition \ref{10.01.06.5}, it 
follows that $S^+$, being the preimage of a point via a continuous map, must be
closed in $(0,+\infty)$. This proves the second statement of the corollary. 

The first statement follows from Lemma \ref{19.09.06.8} (i): the Cauchy transform
$G_{\mu\arc\nu}$ extends continuously and with finite values to $\mathbb R\setminus
S$, so that the density of ${\mu\arc\nu}$ with respect to the Lebesgue measure is
continuous on this set.
\end{pr}

Next, we show that, under some stricter conditions imposed on the rectangular
$R$-transform of $\mu$, we can guarantee that ${\mu\arc\nu}$ has much better
regularity properties. This will follow as a corollary of the proposition below.

\begin{propo}\label{rect-last}
We assume, in addition to
the hypotheses of Proposition \ref{10.01.06.5}, that 
 $$\lim_{x\to-\infty} [C_{\mu}(x)]^2/x\neq 0.$$ Then for any $x\in(0,+\infty)$,
$$M_{\mu\arc\nu}(x)=\lim_{z\to x,z\in\mathbb C^+}M_{\mu\arc\nu}(x)\quad{\rm and}
\quad H_{\mu\arc\nu}(x)=\lim_{z\to x,z\in\mathbb C^+}H_{\mu\boxplus_\lambda\nu}(z)
$$
are finite.
\end{propo}

\begin{pr}
Fix $x\in(0,+\infty)$. 
The existence of the limits has been established
 in Proposition \ref{10.01.06.5}. We shall prove the statement for $H_{\mu\arc\nu}$,
and the statement for $M_{\mu\arc\nu}$ will follow from subsection \ref{Remarks.about.H}, (a).
 We shall prove that this limit is finite   by exploiting 
the asymptotic behaviour of $C_\mu\circ H_{\mu\arc\nu}$ and $M_{\mu\arc\nu}$ as $H_{\mu\arc\nu}$
tends to infinity in order to obtain a contradiction.


Let $c=H_{\mu\arc\nu}(x)$.
Assume towards contradiction that $c=\infty.$ Let $\ell=\lim_{z\to x}\omega_2(z)$, 
where the limit is considered from the upper half-plane (the limit exists by Lemma 
\ref{bdry-subord-rect}). 
By Theorem \ref{lindel}, together
with the above, this implies that $$\lim_{\stackrel{z\longrightarrow \ell
}{\sphericalangle}}
H_\nu(z)=\infty
.$$
 Subsection \ref{Remarks.about.H} 
(a) guarantees that if $H_{\mu\arc\nu}(z)$ tends to infinity as $z$ tends to $x$, then so 
does $M_{\mu\arc\nu}(z)$ and moreover $H_{\mu\arc\nu}(z)/M_{\mu\arc\nu}(z)^2$
tends to $\lambda x$ as $z\to x.$ Also, since $T(M_{\mu\arc\nu}(z))$ and $H_{\mu\arc
\nu}(z)=zT(M_{\mu\arc\nu}(z))$ belong to $\mathbb C\setminus\mathbb R^+$ for $z\in
\mathbb C^+$, we have
$\lim_{z\to x}\frac{\sqrt{H_{\mu\arc\nu}(z)}}{M_{\mu\arc\nu}(z)}=\sqrt{\lambda x}.$
We will use this fact to determine the possible values of $\ell$.



Let us observe that the existence of $\ell$ guarantees the existence of
 $k:=\lim_{t\to-\infty}C_\mu(t)/\sqrt{t}$. Indeed, as the limit of $H_{\mu\arc\nu}$
 at $t$ is infinite, and
 \begin{eqnarray}
 \ell & = & \lim_{z\to x,z\in\mathbb C^+}
 \frac{H_{\mu\arc\nu}(z)}{\lambda(M_{\mu\arc\nu}(z)-C_\mu(H_{\mu\arc\nu}(z)))^2
 +(1+\lambda)(M_{\mu\arc\nu}(z)-C_\mu(H_{\mu\arc\nu}(z)))+1}\nonumber\\
 & = & \frac{1}{\lambda} 
\left( \frac{1}{\sqrt{\lambda x}} - \lim_{z\to x,z\in\mathbb C^+}
\frac{ C_\mu(H_{\mu\arc\nu}(z))}{\sqrt{H_{\mu\arc\nu}(z)}}\right)^{-2}.
\label{elllimit}\\
\nonumber
\end{eqnarray}
On the other hand, as $C_\mu$
satisfies $\arg C_\mu(z)\in(\arg z,\pi),$
(by \eqref{10.01.06.3})
 $C_\mu(\bar{z})=\overline{C_\mu(z)}$ for
$z\in\mathbb C^+$, and $C_\mu(\mathbb R^-)\subseteq\mathbb R^-$, it follows that
Theorem \ref{lindel} applies to the map $w\mapsto\frac{C_\mu(w)}{\sqrt{w}}.$ Thus,
$k:=\lim_{x\to-\infty}\frac{C_\mu(x)}{\sqrt{x}}$ exists 
and is purely imaginary ($k=i|k|$) since  
 $C_\mu(x)$ is negative for $x$ negative.
Thus, \eqref{elllimit}
gives 
$$
\lambda\ell  = \frac{1}{((\lambda x)^{-\frac{1}{2}}-i|k|)^2}.
$$

It follows immediately from equation (\ref{H}) and analyticity of $H_\nu$ on
$\mathbb C^+$ that when $|k|\in(0,+\infty)$, $\ell\in\mathbb C^+$ and we obtain the  contradiction with the fact that  $\infty=H_\nu(\ell)$.

Assume that $k$ is infinite. Then $\ell=0$. But this contradicts Theorem \ref{lindel}
and Remark \ref{rmqhow}: indeed, we obtain that the limit at zero of $H_\nu$ 
along
$\omega_2(z)=\frac{H_{\mu\boxplus_\lambda\nu}(z)}{
T[M_{\mu\arc\nu}(z)-C_\mu(H_{\mu\arc\nu}(z))]}$ (as $z\to x$ from the upper 
half-plane) is infinite (by (\ref{H})), while 
the limit at zero of $H_\nu$ along the negative half-line is zero.

This completes the proof of the proposition since $k$ is not zero if
and only if we have $\lim_{x\to-\infty} [C_{\mu}(x)]^2/x\neq 0$.
\end{pr}

\begin{cor}\label{final.corollary.rectangular}
Under the assumptions of Proposition \ref{rect-last}, $\mu\arc\nu$ is absolutely 
continuous with respect to the Lebesgue measure and its density is continuous.
\end{cor}

\begin{pr}
Since by Proposition \ref{rect-last} $M_{\mu\arc\nu}(x)$ exists and is finite
for all $x\in(0,+\infty),$ the corollary is a consequence of
a variant of  Lemma \ref{19.09.06.8} (which states
the existence of a continuous density of a measure with 
Cauchy transform which extends continuously to the real line)
and the following propositions \ref{atom} and \ref{hole} (which, with lemma \ref{rect-id} (b), allow us to claim that $\mu\arc\nu$ has no atom at the origin). 
\end{pr}

In the following we discuss the issue of analyticity for the density of $\mu\arc\nu$.
\begin{lem}\label{omega2inC+}
Under the hypotheses of Proposition \ref{10.01.06.5}, if $\omega_2(x)\in
\mathbb C^+$, then there exists an $\varepsilon>0$ so that $M_{\mu\arc\nu}$
extends analytically to $(x-\varepsilon,x)\cup(x,x+\varepsilon).$ 
\end{lem}

\begin{pr}
By continuity of $\omega_2$, guaranteed in Lemma \ref{bdry-subord-rect}, there exists $\eta>0$ so that $\omega_2([x-\eta,x+\eta])
\subseteq\mathbb C^+$ is a nontrivial curve in the upper half-plane.  We claim that in 
fact $\omega_2$ is injective on $[x-\eta,x+\eta]$. Indeed, if we assume that $v_1,
v_2\in[x-\eta,x+\eta]$ satisfy $\omega_2(v_1)=\omega_2(v_2),$ then, since $k$ is meromorphic on $\mathbb C\setminus\mathbb R^+$, we obtain
$v_1=k(\omega_2(v_1))=k(\omega_2(v_2))=v_2.$

Let us observe that again since $k(\omega_2(z))=z$ and $k$ is meromorphic on 
$\mathbb C\setminus\mathbb R^+,$ the set $\{w\in\omega_2([x-\eta,x+\eta])\colon
k(w)=\infty{\textrm{ or }} k'(w)=0\}$ is discrete in $\omega_2([x-\eta,x+\eta])$.
If $\omega_2(x)$ belongs to this set, then there exists an $0<\varepsilon\le\eta$ so
that $\omega_2((x-\varepsilon,x))\cup\omega_2((x,x+\varepsilon))$ does not intersect
this set. Thus by the inverse function theorem $\omega_2$ extends analytically
through $(x-\varepsilon,x)\cup(x,x+\varepsilon).$
Otherwise, we apply the inverse function theorem to a neighbourhood of $\omega_2(x)
$ to obtain the same result.

Since $H_{\mu\arc\nu}=H_\nu\circ\omega_2$ and $H_\nu$ is analytic on 
$\mathbb C\setminus\mathbb R^+,$ the statement of the Lemma follows directly from 
subsection \ref{Remarks.about.H} (a).

\end{pr}

\begin{lem}\label{lem48}
Under the hypotheses of Proposition \ref{10.01.06.5}, assume that 
$$H_{\mu\arc\nu}(x)=\lim_{z\to x,z\in\mathbb C^+}H_{\mu\boxplus_\lambda\nu}(z)\in\mathbb C\setminus\mathbb R$$ for some $x\in(0,+\infty).$ Then $\omega_2(x):=
\lim_{z\to x}\omega_2|_{\mathbb C^+}(z)
\in\mathbb C^+.$ 

\end{lem}

\begin{pr}
The equality 
$\omega_2(x)=\frac{H_{\mu\arc\nu}(x)}{T\left[M_{\mu\arc\nu}(x)-C_\mu(H_{\mu\arc\nu}
(x))\right]}$ assures us that $\omega_2(x)$ cannot be infinite. Indeed, assume to the 
contrary that $\omega_2(x)=\infty.$
Then, by Theorem \ref{lindel} and subsection \ref{Remarks.about.H} (f)

$$H_{\mu\arc\nu}(x)=\lim_{z\to x}H_{\mu\arc\nu}(z)=\lim_{z\to x}H_{\nu}(\omega_2(z))
=\lim_{x\to-\infty}H_\nu(x)\in[-\infty,0),$$
a contradiction to our assumption on $H_{\mu\arc\nu}(x).$

Assume first that $H_{\mu\arc\nu}(x)
\in\mathbb C^+.$ 
We show next that $\omega_2(x)\not\in\mathbb R.$ Clearly by subsection \ref{Remarks.about.H} (d), $\omega_2(x)\not\in(-\infty,0].$
Assume again towards contradiction that $\omega_2(x)\in\mathbb R^+$. 

Then  $\frac{H_{\mu\arc\nu}(x)}{x}=T(M_{\mu\arc\nu}(x))$ and $T\left[M_{\mu\arc\nu}(x
)-C_\mu(H_{\mu\arc\nu}(x))\right]$ belong to the same half-line $\chi$ originating at 
zero and
passing through the point $H_{\mu\arc\nu}(x)\in\mathbb C^+$. Thus both points $
M_{\mu\arc\nu}(x)$ and $M_\nu(\omega_2(x))=M_{\mu\arc\nu}(x)-C_\mu(H_{\mu\arc\nu}(x))$ belong to the 
same hyperbola $\mathcal H$ described in Remark \ref{T} (v) whose tangents at the
intersection with $-1/\lambda$ and $-1$ are parallel to $\chi.$ Call these tangents
$\mathcal T_{1/\lambda}$ and $\mathcal T_1$. In particular, since $
H_{\mu\arc\nu}(x)\in\mathbb C^+$, we have that $M_{\mu\arc\nu}(x)\in
\mathcal H\cap K_1$ (recall the notations from Remark \ref{T}), and since 
$\omega_2(\mathbb C^+)\subseteq\mathbb C^+$, we have $M_\nu(\omega_2(x))\in
\mathbb C^+\cup\mathbb R$, and so 
$M_\nu(\omega_2(x))=M_{\mu\arc\nu}(x)-C_\mu(H_{\mu\arc\nu}(x))\in\mathcal H
\cap K_1.$

Now let us recall that $\pi>\arg C_\mu(H_{\mu\arc\nu}(x))\ge\arg H_{\mu\arc\nu}(x)>0$,
so that $ M_{\mu\arc\nu}(x)-C_\mu(H_{\mu\arc\nu}(x))$ has imaginary part strictly
less that the imaginary part of $ M_{\mu\arc\nu}(x)$, so on $\mathcal H\cap K_1$ it
must be below $ M_{\mu\arc\nu}(x)$. But at the same time
$-C_\mu(H_{\mu\arc\nu}(x))$ is in $\mathbb C^-$ and to the right of the line
$\chi\cup-\chi$. Since the tangent $\mathcal T_1$ is parallel to $\chi\cup-\chi$, adding 
this number to $ M_{\mu\arc\nu}(x)$ will give a point in the upper half-plane that is 
necessarily at a greater distance from $\mathcal T_1$ than $ M_{\mu\arc\nu}(x)$, 
and thus it cannot be on $\mathcal H\cap K_1$ in between $-1$ and $M_{\mu\arc\nu}(x
)$ (as this 
part 
of the hyperbola is closer to $\mathcal T_1$ than $M_{\mu\arc\nu}(x)$
is), which provides a contradiction. Thus, if $H_{\mu\arc\nu}(x)\in\mathbb C^+,$ then 
$\omega_2(x)\in\mathbb C^+.$

The case when $H_{\mu\arc\nu}(x)
\in\mathbb C^-$ is similar, and we will only sketch the proof. Indeed, then it is clear that, 
since $M_{\mu\arc\nu}(x)\in\mathbb C^+$ (it cannot be in $\mathbb R$ because of
section \ref{Remarks.about.H} (a)), we must have, with the notations from the previous
case, $M_{\mu\arc\nu}(x)\in\mathcal H\cap K_2$. Now, $-\pi<\arg C_\mu(H_{\mu\arc
\nu}(x))\le\arg H_{\mu\arc
\nu}(x)<0$, and so, as above, $\Im M_{\mu\arc
\nu}(x)<\Im [M_{\mu\arc\nu}(x)-C_\mu(H_{\mu\arc
\nu}(x))].$ This time, however, we obtain that $-C_\mu(H_{\mu\arc
\nu}(x))\in\mathbb C^+,$ and to the {\it right} of $\chi\cup-\chi$. Thus, since $\mathcal
T_{1/\lambda}$ is parallel to $\chi\cup-\chi$, the point $M_{\mu\arc\nu}(x)-C_\mu(
H_{\mu\arc\nu}(x))$ will either be closer to $\mathcal T_{1/\lambda}$ and on its
left side, or it will be on its right side. But the 
part 
 of $\mathcal H\cap K_2$ which 
has an imaginary part greater than the imaginary part of $M_{\mu\arc\nu}(x)$
is on the left side of $\mathcal T_{1/\lambda}$ {\it and} farther away from 
$\mathcal T_{1/\lambda}$ than $M_{\mu\arc\nu}(x)$ is. Contradiction again. Thus, if
$H_{\mu\arc\nu}(x)\in\mathbb C^-,$ then 
$\omega_2(x)\in\mathbb C^+.$
\end{pr}

\begin{lem}\label{lem49}
Under the hypotheses of Proposition \ref{10.01.06.5}, assume that $x\in(0,+\infty)$ is 
so that
$$M_{\mu\arc\nu}(x):=\lim_{z\to x,z\in\mathbb C^+}M_{\mu\arc\nu}(z)\in i(\mathbb R^+
\setminus\{0\})-\frac{1+\lambda}{2\lambda}.$$
Then $\omega_2(x)\in\mathbb C^+.$
\end{lem}
\begin{pr}
The proof of this lemma is immediate. Indeed, by Remark \ref{T} and subsection
\ref{Remarks.about.H} (a), we have $H_{\mu\arc\nu}(x)\in(-\infty, -(1-\lambda)^2/4
\lambda]$, so that $C_\mu(H_{\mu\arc\nu}(x))<0$. But

$$\omega_2(x)=\frac{H_{\mu\arc\nu}(x)}{T\left[M_{\mu\arc\nu}(x)-C_\mu(
H_{\mu\arc\nu}(x))\right]},$$
so by Remark \ref{T}, $\omega_2(x)\in\mathbb C^+.$
\end{pr}

We can now prove the analogue of Theorem \ref{Thinfty} for the rectangular case.

\begin{propo}\label{really-last}
Let $\mu$ and $\nu$ be as in Proposition \ref{rect-last}. Assume in addition that the 
restriction of $C_\mu$ to the upper half-plane extends continuously to $(0,+\infty)$ and
${C_{\mu}}|_{\mathbb C^+}(x)\in\mathbb C^+$ for all $x\in(0,+\infty)$. Then 
$\mu\arc\nu$ is absolutely continuous with respect to the Lebesgue measure, and 
there exists an open set $U\subset\mathbb R$ so that $(\mu\arc\nu)(U)=1$ and
the density $h(x)=\frac{d(\mu\arc\nu)(x)}{dx}$ is analytic on $U$.
\end{propo}

\begin{pr}
We shall use the notations from Proposition \ref{rect-last}.
We know from Proposition \ref{rect-last} that $H_{\mu\arc\nu}(x)$ is finite for any
$x\in(0,+\infty).$ Fix such an $x$.
We show first that $H_{\mu\boxplus_\lambda\nu}(x)\not\in(0,+\infty)$. 
Assume towards contradiction that $H_{\mu\boxplus_\lambda\nu}(x)>0$. 
By Proposition \ref{rect-last}
 $m:=M_{\mu\boxplus_\lambda\nu}(x)$ exists and by subsection 
 \ref{Remarks.about.H} (a), is
real. 
We know that $C_\mu(H_{\mu\arc\nu}(x))\in\mathbb C^+$ by hypothesis. Thus, using (\ref{H}), we get that
\begin{eqnarray*}
H_{\mu\arc\nu}(x)=\lim_{z\to x}H_\nu(\omega_2(z)) & = & \lim_{z\to x}
H_\nu\left(\frac{H_{\nu\arc\mu}(z)}{T\left[M_{\mu\arc\nu}(z)-C_\mu(H_{\mu\arc\nu}(z))
\right]}\right)\\
& = & H_\nu\left(\frac{H_{\mu\arc\nu}(x)}{T\left[m-C_\mu(H_{\mu\arc\nu}(x))\right]}\right).
\end{eqnarray*}
Now, by our hypothesis on $C_\mu,$ we have
$m-C_\mu(H_{\mu\arc\nu}(x))\in\mathbb C^-$. Thus, from the definition of $T$,
$T(m-C_\mu(H_{\mu\arc\nu}(x)))\not\in
[-(1-\lambda)^2/4\lambda,+\infty)\supset[0,+\infty)$. We have reached a
contradiction since $H_\nu(\C\backslash\R^+)\subset \C\backslash\R^+$ by
section
\ref{Remarks.about.H}, (b). 

Thus,  $H_{\mu\arc\nu}((0,+\infty))\subset\mathbb C\setminus(0,+\infty)$. In 
particular, $M_{\mu\arc\nu}(x)\in\mathbb C^+\cup[-1/\lambda,-1].$

Since $[-1/\lambda,-1]$ is a closed set and $M_{\mu\arc\nu}$ is continuous on
$(0,+\infty)$, the set $S=\{x\in(0,+\infty)\colon M_{\mu\arc\nu}(1/x^2)\in[-1/\lambda,-1]\}
$ is closed in $(0,+\infty).$
We claim that the set $S$ 
satisfies $(\mu\arc\nu)(S)=0$. Indeed, the equality $M_{\mu\arc\nu}(1/x^2)
=xG_{\mu\arc\nu}(x)-1$, Lemma \ref{19.09.06.8} (i) and the closeness of $S$ make the 
claim obvious.

We claim next that for any $x\not\in S$, $\omega_2(1/x^2)\in\mathbb C^+$. Indeed,
$x\not\in S$ implies that either $M_{\mu\arc\nu}(1/x^2)\in i(\mathbb R^+
\setminus\{0\})-\frac{1+\lambda}{2\lambda}$, and then the statement follows from
Lemma \ref{lem49}, or $M_{\mu\arc\nu}(1/x^2)\in K_1\cup K_2$ and then
the statement follows from Remark \ref{T} and Lemma \ref{lem48}. 

Now by Lemma \ref{19.09.06.8} (ii), Proposition \ref{rect-last} and Lemma 
\ref{omega2inC+} the statement of the proposition follows. The possibility of the 
existence of an atom at zero will be discarded in the next section.

%
%
%
%
%

\end{pr}

\subsection{Examples}\label{10.01.06.80}
We have a whole family of measures satisfying the previous proposition and 
corollary.  In particular, we are going to see that all $\arc$-stable 
distributions with index strictly
smaller than 2  work. Recall that $\arc$-infinitely 
divisible measures and their L\'evy measures where introduced in section 
\ref{10.01.06.12}.

\begin{propo}\label{10.01.06.8}
Let $G\neq \delta_0$ be a symmetric  positive finite measure on the real line, 
whose 
restriction  to $(0,+\infty)$ admits an analytic positive density. Let $\mu$ be the $\arc$-infinitely divisible measure $\mu$ with L\'evy measure $G$. 

Then  the restriction of 
$C_\mu$ to the upper half-plane extends analytically at any point $x$ of 
$(0,+\infty)$ and satisfies $\Im(C_\mu(x))>0.$
\end{propo}

\begin{pr}Let $\rho$ be the density of the restriction of the positive 
measure  $ G$  to $(0,+\infty)$. By theorem \ref{10.01.06.2}, $C_\mu$ 
extends analytically to $ \C\bck\R^+$ by the formula 
$$C_\mu(z)=z\left(G(\{0\})+2\int_0^{+\infty}\f{(1+t^2)\rho(t)}{1-zt^2}\ud 
t\right)=G_\tau(1/z),
$$ where $\tau$ is push-forward of the  measure $(1+t^2)\ud G(t)$ by the 
function $t\to t^2$, and $G_\tau$ denotes the Cauchy transform of $\tau$. 
Note that $\tau$ is a positive Radon measure, and that  its restriction 
to $(0,+\infty)$ admits the   density $u\to 
\f{(1+u)\rho(u^\ff{2})}{u^\ff{2}}$ on $(0,+\infty)$.
This density is analytic, hence by (ii) of lemma \ref{19.09.06.8} (which 
extends easily to positive measures on the real line which integrate 
$\ff{1+|u|}$, as $\tau$ does) and by the fact that the behavior of 
$G_\tau$ on the lower half-plane can be deduced from its behaviour on the 
upper half-plane by the formula 
$G_\tau(\overline{\cdot})=\overline{G_\tau(\cdot)}$,  the restriction of 
$C_\mu$ to the upper half-plane extends analytically at any point $x$ of 
$(0,+\infty)$ and satisfies, by  (i) of lemma \ref{19.09.06.8}, 
$$\Im(C_\mu(x))= \pi\f{(1+x)\rho(x^\ff{2})}{x^\ff{2}}>0.$$
\end{pr}

It is proved in \cite{fbg05.inf.div} that there is a bijection between the 
set of symmetric $*$-infinitely divisible distributions and the set of
$\arc$-infinitely divisible distributions, which preserves many 
properties, as limit theorems and the fact of being stable. Hence, for all 
$\alpha\in (0,2)$, the set of $\arc$-stable laws $\mu_\alpha$
 with index $\alpha$ is 
the set of $\arc$-infinitely divisible laws which L\'evy measure is of the 
type $t\f{|x|^{1-\alpha}}{1+x^2}\ud x$, where $t$ can be any positive constant, so proposition 
\ref{10.01.06.8} can be applied to them. In  fact, an application of the 
residue formula gives the rectangular $R$-transform with ratio $\la$ of the 
$\arc$-infinitely divisible law  $\mu_\alpha$ with  L\'evy measure 
$\f{|x|^{1-\alpha}}{1+x^2}\ud x$: for all $z\in \C\bck\R^+$, 
$$C_{\mu_\alpha}(z)=-\f{\pi }{2\sin (\pi\alpha/2)}(-z)^{\alpha/2},$$where  the 
power is defined on $\C\bck \R^-$ in relation with the argument with value 
$0$ on the positive half line. 
For  $\alpha\in [1,2)$, $\mu_\alpha$ satisfies the hypothesis
of Proposition \ref{really-last}. As a consequence, for any positive number $t$, the same holds for the $t$-th power, with respect to $\arc$, of $\mu_\alpha$. In deed, if one denotes this measure by  $\mu_\alpha^{ t}$ (it should be denoted by $\mu_\alpha^{\arc t}$, but this notation is a bit hard to swallow), one has $C_{\mu_\alpha^t}(z)=tC_{\mu_\alpha}(z)$.

A matricial model for the measures $\mu_\alpha^{ t}$ was given
in \cite{fbg05.inf.div}, section 5.

 Moreover, for any positive $t$, the density of $\mu_1^t$ 
has been computed in  section 4.2 of \cite{fbg05.inf.div}:
$$\frac{d\mu_1^t}{dx}(x)=
\f{t}{
\pi(\la t^2+x^2)}\lf(1-\f{t^2(\la-1)^2}{4x^2}\ri)^\ff{2}$$on its support $\operatorname{Supp}(\mu_1^t)= { \R\backslash 
\lf(-\f{t(1-\la)}{2}, \f{t(1-\la)}{2}\ri)}$.

Remark however  that the "rectangular Gaussian laws", i.e. the $\arc$-stable laws  with index $2$, which are symmetric square roots of dilations of  Pastur-Marchenko laws, which are the laws $\mu_2^t$, $t>0$, satisfying $C_{\mu_2^t}(z)= tz$,  do not satisfy the  hypotheses 
of Proposition \ref{10.01.06.5}  since $C_{\mu_2^t}((0,\infty))\subset
(0,\infty)$.

\subsection{Study of the density around  the origin}\label{studyzero}

In this section, we study the existence of a hole around the
origin in the support of the free convolution 
$\mu\arc\nu$. Since in our approach the origin itself is 
a very specific point,
we shall study separately the existence of an atom
at the origin and then 
existence of a set $[-\e,\e]$ 
which does not intersect 
the support of $\mu\arc\nu$.

Some of the considerations of this section do not require
the assumptions of Proposition \ref{10.01.06.5}.

\begin{propo}\label{atom}
\begin{enumerate}\item For all symmetric \pro measures $\mu,\nu$,  
$(\mu\arc\nu)(\{0\})\geq \mu(\{0\})+\nu(\{0\})-1$. 
\item Assume that $\mu$ is $\arc$-infinitely divisible ($\nu$ is still an arbitrary symmetric
\pro measure). If $(\mu\arc\nu)(\{0\})>0$, then $\mu(\{0\})+\nu(\{0\})>1$ and
$(\mu\arc\nu)(\{0\})=\mu(\{0\})+\nu(\{0\})-1$.
\end{enumerate}
\end{propo}


\begin{pr}
We prove item 1. Consider a sequence $p_n\geq n$ of positive integers such that $$n/p_n\ninf \la$$ and, on a \pro space $\Omega$,  an independent set  of random variables $$(X_i)_{i\geq 1}, (Y_i)_{i\geq 1}, (U_n)_{n\geq 1}, (V_n)_{n\geq 1}$$ such that \begin{itemize}\item[-] each $X_i$ is distributed according to $\mu$, \item[-] each $Y_i$ is distributed according to $\nu$, \item[-] for all $n$, $U_n$ is an $n$ by $n$ Haar-distributed unitary random matrix, \item[-] for all $n$, $V_n$ is a $p_n$ by $p_n$ Haar-distributed unitary random matrix.\end{itemize} 
Let us define, for all $n$, \begin{itemize}\item[a)] $M_n$ to be the $n$ by $p_n$ random matrix with $|X_1|,\ldots, |X_n| $ on the diagonal, and zeros everywhere else, \item[b)] $N_n$ to be $U_n$ times the $n$ by $p_n$ random matrix with $|Y_1|,\ldots, |Y_n| $ on the diagonal, and zeros everywhere else times $V_n$.\end{itemize}

Let, for all $n$, $d_n$ (resp. $d'_n, d_n''$) be the random variable equal to the number of null singular values of $M_n$ (resp. $N_n, M_n+N_n$). 
 Note that $d_n+(p_n-n)=\operatorname{dim}\ker M_n,$ $d'_n+(p_n-n)=\operatorname{dim}\ker N_n, $ $d''_n+(p_n-n)=\operatorname{dim}\ker M_n+N_n.$ Note also that since $\ker M_n \cap \ker N_n \subset \ker (M_n+N_n)$, one has $$\operatorname{dim}\ker M_n+N_n\geq \operatorname{dim} \ker M_n \cap \ker N_n\geq \operatorname{dim}\ker M_n+\operatorname{dim}\ker N_n
-p_n,$$ hence $$d''_n+(p_n-n)\geq d_n+(p_n-n)+d'_n+(p_n-n)-p_n,$$i.e.\begin{equation}\label{2.3.07.1}d''_n\geq d_n+d'_n-n.\end{equation}

Note that the singular values of $M_n$ (resp. of $N_n$) are $|X_1|,\ldots, |X_n| $ (resp. $|Y_1|,\ldots, |Y_n| $), hence by the law of large numbers, the symmetrization of the singular law of $M_n$ (resp. of  $N_n$) converges almost surely weakly to $\mu$ (resp.$\nu$). Thus by theorem 4.8 of \cite{benaych.rectangular}, the singular law $\operatorname{SL}(M_n+N_n)$ 
of $M_n+N_n$ converges in \pro to $\mu\arc \nu$ in the metric space of the set of \pro measures on the real line endowed with a distance which defines the weak convergence. So for almost all $\omega\in \Omega$, there is a subsequence ${\varphi(n)}$ of the sequence $\operatorname{SL}(M_n+N_n)(\omega)$ which converges weakly to $\mu\arc \nu$. For such an $\omega$, one has \begin{equation}\label{2.3.07.2}\ds(\mu\arc\nu)(\{0\})\geq \limsup_{n\to \infty} \operatorname{SL}(M_{\varphi(n)}+N_{\varphi(n)})(\omega)(\{0\})= \limsup_{n\to \infty} \f{d_{\varphi(n)}''}{{\varphi(n)}}.\end{equation}

Note moreover that for all $n$,  $d_n$ (resp. $d_n'$) is the number of $i$'s in $\{1,\ldots,n\}$ \st $X_i=0$ (resp. $Y_i=0$), hence  the law of large numbers implies also that for almost all $\omega\in \Omega$
\begin{equation}\label{2.3.07.3} \f{d_n(\omega)}{n}\ninf \mu(\{0\}), \quad \f{d_n'(\omega)}{n}\ninf \nu(\{0\}).\end{equation}

Putting together \eqref{2.3.07.1}, \eqref{2.3.07.2}, \eqref{2.3.07.3}, one gets $(\mu\arc\nu)(\{0\})\geq \mu(\{0\})+\nu(\{0\})-1$.

Let us now prove item 2. First of all, we exclude the case $\nu=\delta_0$, which is trivial. The strategy will be to use Lemma \ref{rect-id} and the 
description of atoms given in \ref{Remarks.about.H} (f) together with the formula
(\ref{H}) to prove the equality $(\mu\arc\nu)(\{0\})=\mu(\{0\})+\nu(\{0\})-1$.
We shall first prove that $\lim_{x\to-\infty}C_\mu(x)>-1$. Note that by \ref{Remarks.about.H} (f),
our hypothesis implies that $\lim_{x\to-\infty}H_{\mu\arc\nu}(x)=-\infty$. So we will prove that if $\lim_{x\to-\infty}H_{\mu\arc\nu}(x)=-\infty,$ then 
$\lim_{x\to-\infty}C_\mu(x)>-1$. 
For future use, we prefer to 
prove now that $\lim_{x\to-\infty}C_\mu(x)>-1$ under the hypothesis $\lim_{x\to-\infty}H_{\mu\arc\nu}(x)=-\infty$  than under the stronger one of the proposition. 

Assume thus that $\lim_{x\to-\infty}H_{\mu\arc\nu}(x)=-\infty.$
Recall first equality (\ref{H}) : \begin{equation}\label{soulstorm-nile}\forall z\in(-\infty,0),\quad H_{\mu\arc\nu}(z)=H_\nu\left(\frac{H_{\nu\arc\mu}(z)}{T\left[M_{\mu\arc\nu}(z)-C_\mu(H_{\mu\arc\nu}(z))
\right]}\right),
\end{equation}
where $M_{\mu\arc\nu}(z)$ is 
the analytic extension of
 $U\left(\frac{H_{\mu\arc\nu}(z)}{z}-1\right)$
  which can be found in \ref{Remarks.about.H} (a)(and which allows us to claim that $\lim_{x\to-\infty}M_{\mu\arc\nu}(x)$ exists and is equal to $\mu\arc\nu(\{0\})-1$
). 

 Equality (\ref{soulstorm-nile})
together with 
the continuity of $H_\nu$ on 
$(-\infty, 0]$  and the hypothesis $\displaystyle\lim_{x\to-\infty}H_{\mu\arc\nu}(x)=-\infty$ imply that 
$$\lim_{x\to-\infty}\frac{H_{\mu\arc\nu}(x)}{T\left[M_{\mu\arc\nu}(x)-C_\mu(H_{\mu\arc\nu}(x))\ri]}=-\infty$$(indeed, for any sequence $x_n$ of negative numbers which tends to $-\infty$, since by \ref{Remarks.about.H}, for all $n$, $y_n:=\frac{H_{\mu\arc\nu}(x_n)}{T\left[M_{\mu\arc\nu}(x_n)-C_\mu(H_{\mu\arc\nu}(x_n))\ri]}\in(-\infty, 0)$, if $y_n$ doesn't tend to $-\infty$, a subsequence of $H_{\mu\arc\nu}(x_n)$ will have a finite  limit, which is impossible). So 
by \ref{Remarks.about.H} (f) we obtain
\begin{eqnarray*}
T(\nu(\{0\})-1) & = & \lambda\nu(\{0\})^2+(1-\lambda)\nu(\{0\})\\
 & = & \lim_{x\to-\infty} \f{H_\nu(x)}{x}\\ &=&
\lim_{x\to-\infty}\frac{T(M_{\mu\arc\nu}(x)-C_\mu(H_{\mu\arc\nu}(x)))}{H_{\mu\arc\nu}(x)}
H_\nu\left(\frac{H_{\mu\arc\nu}(x)}{T(M_{\mu\arc\nu}(x)-C_\mu(H_{\mu\arc\nu}(x)))}\right)\\
& = & \lim_{x\to-\infty}\frac{T(M_{\mu\arc\nu}(x)-C_\mu(H_{\mu\arc\nu}(x)))}{H_{\mu\arc\nu}(x)}H_{\mu\arc\nu}(x)\\
& = & \lim_{x\to-\infty}T(M_{\mu\arc\nu}(x)-C_\mu(H_{\mu\arc\nu}(x))).
\end{eqnarray*}
Thus, $\lim_{x\to-\infty}T(M_{\mu\arc\nu}(x)-C_\mu(H_{\mu\arc\nu}(x)))\in[0,1]$. Note that  \eqref{10.01.06.3} allows us to claim that  $\lim_{x\to-\infty}C_{\mu}(x)$ exists in $[-\infty, 0]$, hence  by above, $M_{\mu\arc\nu}(x)-C_\mu(H_{\mu\arc\nu}(x))$ has also a limit $l=(\mu\arc\nu)(\{0\})-1-\lim_{w\to -\infty}C_\mu(w)\geq -1$ as $x$ goes to $-\infty$. 
Since $T^{-1}([0,1])=[-\ff{\la}-1,-\ff{\la}]\cup[-1,0]$, one has $l\in [-1,0]$, hence $l=\nu(\{0\})-1$.
We conclude that 
$$\lim_{w\to-\infty}C_\mu(w)=
\lim_{x\to-\infty}C_\mu(H_{\mu\arc\nu}(x))=(\mu\arc\nu)(\{0\})-1-(\nu(\{0\})-1)\in
(-1,0],$$ as claimed.
Moreover, this equality together with Lemma \ref{rect-id} implies that 
$\mu(\{0\})-1=(\mu\arc\nu)(\{0\})-\nu(\{0\})$, which is equivalent to item 2.
 \end{pr}

\begin{propo}\label{hole}
Let $\mu$ be $\arc$-infinitely divisible and $\nu$ be arbitrary. Assume that $\nu(\{0\})+
\mu(\{0\})<1.$ Then ${\rm supp}(\mu\arc\nu)$ has a hole around the origin.
\end{propo}

\begin{pr}
Let us denote $r:=\lim_{x\to-\infty}H_{\mu\arc\nu}(x)$
(which exists and belongs to  $[-\infty,0)$ by  \ref{Remarks.about.H} (f)). The first step in our proof will be to show that
under our hypothesis, $r>-\infty.$ Then we will view $r$ as the Denjoy-Wolff point of
a certain self-map of the left half-plane $i\mathbb C^+$, and use this fact to see
$W_{\mu\arc\nu}(x)=H_{\mu\arc\nu}(1/x)$ as an implicit function which is defined on
a neighbourhood of zero and extends $H_{\mu\arc\nu}(1/x)$ from the left half-line.
Finally, we will argue that on a small 
enough interval, $H_{\mu\arc\nu}(1/x)
\in [-(1-\lambda)^2/4\la x
,0]$ for all $x>0$
small  enough, which is equivalent to the existence 
of an open  neighborhood
of the origin  which does not intersect 
 the support of $\mu\arc\nu$ (as can be checked by 
using lemma \ref{19.09.06.8} and remark  \ref{rmqhow}).

We shall prove that $r$ is the Denjoy-Wolff point of $f_1$
if 
 $$f_t(z)=H_\nu\left(\frac{z}{T[-t-C_\mu(z)]}\right),\, t\in[0,1].$$
 First, we claim that $f_t$ is defined on the left 
half-plane $i\mathbb C^+$, and moreover that $f_t(i\mathbb C^+)\subseteq
i\mathbb C^+$ for all $t\in[0,1]$. Indeed,  from Theorem 
\ref{10.01.06.2} it follows that $C_\mu(\mathbb C^+)\subseteq\mathbb C^+$ and
$\arg C_\mu(z)>\arg z$ for any $z\in\mathbb C^+$. Thus, since $0<\lambda<1,$
and $0\le t\le1$, 
$\pi>\arg \left(C_\mu(z)+t-1\right)\ge\arg C_\mu(z)>\arg z$ and
$\pi>\arg \left(C_\mu(z)+t-\frac1\lambda\right)>\arg C_\mu(z)>\arg z$, so that
$\arg T[-t-C_\mu(z)]\in(2\arg z,2\pi)$ for all $z\in\mathbb C^+\cap i\mathbb C^+
.$ 
We conclude that $\arg \left(\frac{T[-t-C_\mu(z)]}{z}\right)\in(\arg z,2\pi-\arg z)
\subset(\pi/2,3\pi/2)$ for any $z\in\mathbb C^+\cap i\mathbb C^+
,$ so that $\frac{T[-t-C_\mu(z)]}{z}$ maps $\mathbb C^+\cap i\mathbb C^+$
in $i\mathbb C^+$. Since $i\mathbb C^+$ is invariant under the maps 
$z\mapsto 1/z$
and $z\mapsto\bar{z}$,
and $\overline{\left(\frac{T[-t-C_\mu(z)]}{z}\right)}=\frac{T[-t-C_\mu(\bar z)]}{\bar z}$,
we conclude that $z\mapsto \frac{z}{T[-t-C_\mu(z)]}$ maps $i\mathbb C^+$
into itself. Since, by subsection \ref{Remarks.about.H} (d), $H_\nu(i\mathbb C^+)
\subseteq i\mathbb C^+$, our claim is proved.

By the last remark, we also have that
$f_t(\bar{z})=\overline{f_t(z)},$  for all
$t\in[0,1]$, and hence in particular $f_t((-\infty,0])\subset (-\infty,0]$.

We next show the existence and uniqueness of the
Denjoy Wolff point of $f_1$
 as a consequence of 
Theorem 
\ref{Denj}.
In fact, 
$f_1$ is not a conformal automorphism of $i\mathbb C^+$. Indeed, 
there are only two  conformal automorphisms of $i\C^+$ 
which fix  $ (-\infty,0]$ up to multiplication by positive
scalar; 
 the identity and $z\to1/z$. The case
$az=H_\nu\left(\frac{z}{T[-1-C_\mu(z)]}\right)$ can be discarded since as $z$
goes to zero 
 along the negative half-line, $C_\mu(z)/z$
converges to $\int (1+t^2)dG(t)$ by monotone convergence theorem
(with $G$ the L\'evy measure of $\mu$)
and so $y(z):=\frac{z}{T[-1-C_\mu(z)]}\in (-\infty,0]$ goes to the constant
\begin{equation}\label{eqr}
 \lim_{x\uparrow 0}\frac{x}{T[-1-C_\mu(x)]}=\lim_{x\uparrow0}\frac{1}{\lambda\frac{
C_\mu(x)^2}{x}+(\lambda-1)\frac{C_\mu(x)}{x}}=\left((\lambda-1)\int(t^2+1)\,dG(t)
\right)^{-1},
\end{equation}
which is null only if $\int (1+t^2)dG(t)$ is infinite. 
If this  constant does not vanish, we obtain a contradiction
since $H_\nu$ does not vanish on $(-\infty,0)$. 
If it vanishes, we write $aT[-1-C_\mu(z)]=H_\nu(y(z))/y(z)$ 
with $y(z)$ negative going to zero
as $z$ goes to zero. 
This is in contradition with the fact
that $H_\nu(z)/z$ goes to one (see Remark \ref{rmqhow}).
The case 
$a/z=H_\nu\left(\frac{z}{T[-1-C_\mu(z)]}\right)$ leads also to
a contradiction by letting $z$ going to zero.

The uniqueness of the Denjoy-Wolff
point  given by Theorem 
\ref{Denj} implies that this point can only belong 
to $[-\infty,0]$ since $f_t(\bar z)=\bar f_t(z)$.
 We shall first show that 
zero cannot be the Denjoy-Wolff
point of $f_1$, and secondly 
we  show that infinity
can be  the  Denjoy-Wolff
point of $f_1$ only when $\mu(\{0\})+\nu(\{0\})\ge 1$. 

For zero to be the Denjoy-Wolff
point of $f_1$, we would first need to 
have $\lim_{x\uparrow0}f_1(x)=0.$ 
Since $H_\nu$ vanishes on $(-\infty,0]$
only at the origin, we must have by \eqref{eqr}
that the L\'{e}vy measure of $\mu$ has infinite second moment.
The second requirement for zero to be the Denjoy-Wolff
point of $f_1$ is that $\lim_{x\uparrow0}f_1(x)/x\in(0,1].$ But
\begin{eqnarray*}
\lim_{x\uparrow0}\frac{f_1(x)}{x} & = & \lim_{x\uparrow0}
\frac{H_\nu\left(\frac{x}{T[-1-C_\mu(x)]}\right)}{\frac{x}{T[-1-C_\mu(x)]}}\cdot
\frac{1}{T[-1-C_\mu(x)]}\\
& = & \lim_{x\uparrow0}\frac{H_\nu(x)}{x}\cdot\lim_{x\uparrow0}\frac{1}{T[-1-C_\mu
(x)]} =  \infty,
\end{eqnarray*}
since $\lim_{x\uparrow0}\frac{H_\nu(x)}{x}=1$ 
and $\lim_{x\uparrow0}C_\mu(x)=0$.
We conclude that zero cannot be the Denjoy-Wolff point of $f_1$. 

Now we show under that under our
condition $\nu(\{0\})+
\mu(\{0\})<1$, 
$f_1$ cannot have
infinity as Denjoy-Wolff point (recall 
that $(\mu\arc\nu)(\{0\})=0$  by Proposition \ref{atom}
under this assumption.) The two
requirements that $f_1$
must  verify to have infinity as
Denjoy-Wolff point
are $\lim_{x\to-\infty}f_1(x)=-\infty$ and $
\lim_{x\to-\infty}f_1(x)/x\in[1,+\infty).$ The continuity of $H_\nu$ on 
$(-\infty,0]$ translates the first requirement into $\lim_{x\to-\infty}
\frac{x}{T[-1-C_\mu(x)]}=-\infty$ and $\lim_{x\to-\infty}H_\nu(x)=-\infty$. Applying
\ref{Remarks.about.H} (f) and the above, we obtain:

\begin{eqnarray*}
\lim_{x\to-\infty}\frac{f_1(x)}{x} & = & \lim_{x\to-\infty}
\frac{H_\nu\left(\frac{x}{T[-1-C_\mu(x)]}\right)}{\frac{x}{T[-1-C_\mu(x)]}}\cdot
\frac{1}{T[-1-C_\mu(x)]}\\
& = & \lim_{x\to-\infty}\frac{H_\nu(x)}{x}\cdot\lim_{x\to-\infty}\frac{1}{T[-1-C_\mu
(x)]}\\
& = & (\lambda\nu(\{0\})^2+(1-\lambda)\nu(\{0\}))\cdot\frac{1}{\lambda c^2+(\lambda-1)c}\\
& = & \frac{T(\nu(\{0\})-1)}{T(-c-1)},
\end{eqnarray*}
where $c:=\lim_{x\to-\infty}C_\mu(x)\in[-\infty,0).$ To have $
\lim_{x\to-\infty}f_1(x)/x\ge1$, we must have $\nu(\{0\})>0$ and $c>-\infty$.
Thus, we may write
$$T(\nu(\{0\})-1)\ge T(-c-1),$$
which implies,
since $T$ is increasing on $[-1,+\infty)$,
that   $1\ge \nu(\{0\})\ge-c.$ Using Lemma \ref{rect-id} (2)
we conclude that $\nu(\{0\})+\mu(\{0\})\ge1.$ 

Thus, $\nu(\{0\})+\mu(\{0\})<1$ implies that $f_1$ has a Denjoy-Wolff point
$s\in(-\infty,0).$ 
We claim that $s=r$. Indeed, by taking limit when $z\to-\infty$
in equation (\ref{H}), using Proposition \ref{atom} and the fact that $\lim_{x\to-\infty}
M_{\mu\arc\nu}(x)=(\mu\arc\nu)(\{0\})-1=-1$, we obtain that
\begin{equation}\label{leqr}
r=H_\nu\left(\frac{r}{T[-1-C_\mu(r)]}\right)=f_1(r),
\end{equation}
where, if $r=-\infty,$  the second term must be also understood as a limit.

We finally
show that $r$ cannot be infinite, which will imply with \eqref{leqr} that
$r=s$. 
By Theorem \ref{Denj}, 
 $f'_1(s)\in
(-1,1)$. By continuity
of $f'_1$, 
there exists $\delta>0$
so that  if $D=\{y:|x-s|<\delta\}$, 
  $\rho:=\sup_{\bar D}
|f'_1(x)|<1$ and therefore 
$f_1(\bar{D})\subset \{y:|y-s|\le \rho\delta\}\subset D.$ 
Since $f_t$ converges to $f_1$ as $t\to1$ uniformly on 
compact subsets of $i\mathbb C^+$,
 there exists $\varepsilon>0$ so that $f_t(\bar{D})
\subseteq D$
and moreover 
 $f_t$ is not an hyperbolic 
rotation  for all $1-\varepsilon\le t\le1$.
 Thus, by Theorem \ref{Denj},
$f_t$ has a unique Denjoy-Wolff point
and it must be in $\bar D$ (has can be seen by iterating $f$
from a point in $\bar D$). Thus, 
the Denjoy-Wolff points of $f_t$ converge to $s$
as $t\to1$. 

Now, since $\lim_{x\to-\infty}
M_{\mu\arc\nu}(x)=-1$ as
$\mu\arc\nu(\{0\})=0$, for $x$ large enough we have $M_{\mu\arc\nu}(x)\in
(-1,\varepsilon-1).$ From equation (\ref{H}) it follows that
$$H_{\mu\arc\nu}(x)=f_{-M_{\mu\arc\nu}(x)}(H_{\mu\arc\nu}(x))$$
and therefore 
 $H_{\mu\arc\nu}(x)$ is the
Denjoy-Wolff point of $f_{-M_{\mu\arc\nu}(x)}.$
Thus we conclude that $r=\lim_{x\to-\infty}H_{\mu\arc\nu}(x)=s$, which proves our 
claim. 

Let us define $W_{\mu\arc\nu}(z)=H_{\mu\arc\nu}(1/z)$ and 

$$g(x,w)=H_\nu\left(\frac{w}{T[U(xw-1)-C_\mu(w)]}\right)-w.$$
It is easy to observe that for $x<0$ close to zero, we have $g(x,W_{\mu\arc\nu}(x))=0
$, as the formula $M_{\mu\arc\nu}(z)=U(H_{\mu\arc\nu}(z)/z-1)$ must hold for all
$z\in\mathbb R^-.$ Moreover, there obviously exists a small enough 
interval $I$ centered at zero so that $g$ is actually defined on $I\times(I+r)$, and
of course, by equation (\ref{leqr}), $g(0,r)=0$. Let us differentiate $g$
with respect to $w$:
\begin{eqnarray*}
\partial_wg(x,w) & = & H_\nu' \left(\frac{w}{T[U(xw-1)-C_\mu(w)]}\right)\\
& &  \mbox{}\times
\frac{T[U(xw-1)-C_\mu(w)]-wT'[U(xw-1)-C_\mu(w)][xU'(xw-1)-C_\mu'(w)]}{
T[U(xw-1)-C_\mu(w)]^2}-1.
\end{eqnarray*}
Since $U$ is differentiable in $-1$, and $T[-1-C_\mu(r)]\neq0$,
we have
$$\partial_wg(0,r) = H_\nu' \left(\frac{r}{T[-1-C_\mu(r)]}\right)
\frac{T[-1-C_\mu(r)]-rT'[-1-C_\mu(r)][-C_\mu'(r)]}{
T[-1-C_\mu(r)]^2}-1=f_1'(r)-1.$$
Since we have shown that $|f_1'(r)|<1,$ we conclude that $\partial_wg(0,r)\neq0,$
so we can apply the implicit function theorem to it in the point $(0,r)$ to 
extend $W_{\mu\arc\nu}$ to a small neighborhood
of the origin. Then, $W_{\mu\arc\nu}(x)=H_{\mu\arc\nu}(1/x)
$ takes its values in a finite
neighborhood of $r\in (-\infty,0)$
which is included into $[-(1-\lambda)^2/4\lambda x,0]$
for sufficiently small $x$, and hence $\mu\arc\nu$ put no mass
in a open neighborhood of the origin. 
\end{pr}

\end{document}